\documentclass[12pt]{article}
\usepackage{a4wide}
\usepackage[french,english]{babel}
\usepackage[utf8]{inputenc}
\usepackage[T1]{fontenc}
\usepackage{amsthm}
\usepackage{amssymb}
\usepackage{amsmath}
\usepackage{stmaryrd}
\usepackage{makeidx}
\usepackage{dsfont}
\usepackage{graphicx}
\usepackage{caption}
\theoremstyle{plain}
\usepackage[colorlinks=true,breaklinks=true,linkcolor=black]{hyperref} %pour faire des liens hypertextes
\usepackage{amsfonts}

\newtheorem{theorem}{Theorem}[section] %[subsection sert à numéroter les théorèmes d'après la sous-section où ils se trouvent
\newtheorem{proposition}[theorem]{Proposition}
\newtheorem{corollary}[theorem]{Corollary}
\newtheorem{definition}[theorem]{Definition}
\newtheorem{lemma}[theorem]{Lemma}

\newtheorem{remark}[theorem]{Remark}
\newtheorem{definition/proposition}[theorem]{Definition/Proposition}
\date{}
\makeindex
\usepackage{mathrsfs}
\usepackage[all]{xy}
\theoremstyle{definition}
\usepackage[left,modulo]{lineno}
\usepackage{geometry}
\usepackage{comment}

\title{Principal series  representations of Iwahori-Hecke algebras for Kac-Moody groups over local fields}
\author{Auguste \textsc{Hébert} \\École normale supérieure de Lyon\\ UMR 5669 CNRS,
auguste.hebert@univ-lorraine.fr}

\hypersetup{
pdfauthor={Hebert, Auguste},
pdftitle={Principal series representations},
}

\theoremstyle{definition}\newtheorem{thm*}{Theorem}
\theoremstyle{definition}\newtheorem{proposition*}{Proposition}
\makeatletter \@addtoreset{figure}{section}\makeatother

\newcommand{\R}{\mathbb{R}}
\newcommand{\A}{\mathbb{A}}

\newcommand{\N}{{\mathbb{\Z}_{\geq 0}}}
\newcommand{\Z}{\mathbb{Z}}
\newcommand{\Q}{\mathbb{Q}}
\newcommand{\C}{\mathbb{C}}
\newcommand{\Ne}{{\mathbb{Z}_{\geq 1}}}

\newcommand{\I}{\mathcal{I}}
\newcommand{\T}{\mathcal{T}}
\newcommand{\Id}{\mathrm{Id}}

\newcommand{\conv}{\mathrm{conv}}

\newcommand{\HH}{\mathcal{H}}
\newcommand{\supp}{\mathrm{supp}}

\newcommand{\h}{\mathrm{ht}}

\newcommand{\Hom}{\mathrm{Hom}}

\newcommand{\AC}{{^{\mathrm{BL}}\mathcal{H}}}
\newcommand{\ATF}{\AC(T_\FC)}
\newcommand{\AF}{\AC_\FC}
\newcommand{\ATC}{\AC(T_\C)}
\newcommand{\BC}{\mathcal{B}}

\newcommand{\DC}{\mathcal{D}}
\newcommand{\FC}{\mathcal{F}}
\newcommand{\HC}{\mathcal{H}}
\newcommand{\KC}{\mathcal{K}}
\newcommand{\OC}{\mathcal{O}}

\newcommand{\RC}{\mathcal{R}}
\newcommand{\SC}{\mathcal{S}}
\newcommand{\UC}{\mathcal{U}}
\newcommand{\VC}{\mathcal{V}}

\newcommand{\HCW}{\mathcal{H}_{W^v,\C}}

\newcommand{\HWC}{\mathcal{H}_{W^v,\C}}
\newcommand{\HFW}{{\mathcal{\HC}_{W^v,\FC}}}
\newcommand{\ev}{\mathrm{ev}}

\newcommand{\CCC}{\mathscr{C}}
\newcommand{\RCC}{\mathscr{R}}
\newcommand{\SCC}{\mathscr{S}}

\newcommand{\Wt}{\mathrm{Wt}}

\newcommand{\Wta}{{W_{\!(\tau)\!}}}
\newcommand{\wb}{\mathbf{w}}
\newcommand{\dw}{d^{W^+}}
\newcommand{\CO}{\mathscr{C}_0}

\newcommand{\ve}{\mathrm{vert}}
\newcommand{\dv}{d^{Y^{++}}}
\newcommand{\dy}{d^{Y^+}}
\newcommand{\vb}{\mathbf{v}}

\begin{document}

\maketitle

\begin{abstract}
Recently, Iwahori-Hecke algebras were associated with Kac-Moody groups over non-Archimedean local fields. We introduce principal series representations for these algebras. We study these representations and partially generalize irreducibility criteria of Kato and Matsumoto.
\end{abstract}

\section{Introduction}

\subsection{The reductive case}

Let  $G$ be a split reductive group over a non-Archimedean local field $\KC$. Let $T$ be a maximal split torus of $G$ and $Y$ be the cocharacter lattice of $(G,T)$. Let $B$ be a Borel subgroup of $G$ containing $T$.  Let $T_\C=\Hom_{\mathrm{Gr}}(Y,\C^*)$. Then $\tau$ can be extended to a character $\tau:B\rightarrow \C^*$. If $\tau\in T_\C$,  the principal series representation $I(\tau)$  of $G$ is the induction of $\tau\delta^{1/2}$ from $B$ to $G$, where $\delta:B\rightarrow \R^*_+$ is the modulus character of $B$. More explicitly, this  is the space of locally constant functions $f:G\rightarrow \C$ such that $f(bg)=\tau\delta^{1/2}(b)f(g)$ for every $g\in G$ and $b\in B$. Then $G$ acts on $I(\tau)$ by right translation.

 To each open compact subgroup $K$ of $G$ is associated the Hecke algebra $\HC_K$. This is the algebra of functions from $G$ to $\C$ which have compact support and are $K$-bi-invariant.  There exists a strong link between the smooth representations of $G$ and the representations of the Hecke algebras of $G$. Let $K_I$ be the Iwahori subgroup of $G$. Then the Hecke algebra $\HC_\C$ associated with $K_I$ is called the Iwahori-Hecke algebra of $G$ and plays an important role in the representation theory of $G$.

 The algebra  $\HC_\C$ acts on $I_{\tau,G}:=I(\tau)^{K_I}$ by the formula \[\phi.f=\int_G \phi(g) g.fd\mu(g),\forall(\phi,f)\in \HC_\C\times I(\tau)^{K_I},\] where $\mu$ is a Haar measure on $G$.  This formula can actually be rewritten as \begin{equation}\label{eqFormula_action_IH_algebra}\phi.f=\mu(K_I)\sum_{g\in G/K_I} \phi(g) g.f,\forall(\phi,f)\in \HC_\C\times I(\tau)^{K_I}.\end{equation} Then $I(\tau)$ is irreducible as a representation of $G$ if and only $I_{\tau,G}$ is irreducible as a representation of $\HC_\C$.

Let  $W^v$ be the vectorial Weyl group of $(G,T)$. By the Bernstein-Lusztig relations, $\HC_\C$ admits a basis  $(Z^\lambda H_w)_{\lambda\in Y, w\in W^v}$ such that $\bigoplus_{\lambda\in Y}\C Z^\lambda$ is a subalgebra of $\HC_\C$ isomorphic to the group algebra $\C[Y]$ of $Y$. We identify $\bigoplus_{\lambda\in Y}\C Z^\lambda$ and $\C[Y]$. We regard $\tau$ as an algebra morphism $\tau:\C[Y]\rightarrow\C$. Then $I_{\tau,G}$ is isomorphic to the induced representation $I_\tau=\mathrm{Ind}_{\C[Y]}^{\HC_\C}(\tau)$ and we refer to \cite[Section 3.2]{solleveld2009periodic} for a survey on this subject.

Matsumoto and Kato gave criteria for the irreducibility of $I_\tau$.  The group  $W^v$ acts on $Y$ and thus it acts on $T_\C$. If $\tau\in T_\C$, we denote by $W_\tau$ the stabilizer  of $\tau$ in $W^v$. Let $\Phi^\vee$ be the coroot lattice of $G$. Let $q$ be the residue cardinal of $\KC$.  Let $\Wta$ be the subgroup of $W_\tau$ generated by the reflections $r_{\alpha^\vee}$, for $\alpha^\vee\in \Phi^\vee$ such that $\tau(\alpha^\vee)=1$. Then Kato proved the following theorem (see \cite[Theorem 2.4]{kato1982irreducibility}):
 
 \begin{thm*}\label{thm*Kato's theorem}
 Let $\tau\in T_\C$. Then $I_\tau$ is irreducible if and only if it satisfies the following conditions: \begin{enumerate}
\item\label{itWchi engendré par ses réflexions} $W_\tau=\Wta$,

\item for all $\alpha^\vee\in \Phi^\vee$, $\tau(\alpha^\vee)\neq q$. 
 \end{enumerate}
 \end{thm*}

When $\tau$ is \textbf{regular}, that is when $W_\tau=\{1\}$, condition~(\ref{itWchi engendré par ses réflexions}) is satisfied and this is a result by Matsumoto (see \cite[Th{\'e}or{\`e}me 4.3.5]{matsumoto77Analyse}).

\subsection{The Kac-Moody case}

Let $G$ be a split Kac-Moody group over a non-Archimedean local field $\KC$. We do not know which topology on $G$ could replace the usual topology on reductive groups over $\KC$.  There is up to now no definition of smoothness for the representations of $G$. However one can define certain Hecke algebras in this framework. In \cite{braverman2011spherical} and \cite{braverman2016iwahori}, Braverman, Kazhdan and Patnaik defined the spherical Hecke algebra and the Iwahori-Hecke $\HC_\C$ of $G$ when $G$ is affine. In \cite{gaussent2014spherical} and \cite{bardy2016iwahori}, Bardy-Panse, Gaussent and Rousseau generalized these constructions to the case where $G$ is a general Kac-Moody group. They achieved this construction by using masures (also known as hovels), which are analogous to Bruhat-Tits buildings (see \cite{gaussent2008kac}). Together with Abdellatif, we attached   Hecke algebras to subgroups slightly more general than the Iwahori subgroup (see \cite{abdellatif2019completed}).

Let $B$ be a positive Borel subgroup of $G$ and $T$ be a maximal split torus of $G$ contained in $B$.  Let $Y$ be the cocharacter lattice of $G$, $W^v$ be the Weyl group of $G$ and $Y^{++}$ be the set of dominant cocharacters of $Y$. The Bruhat decomposition does not hold on $G$: if $G$ is not reductive, \[G^+:=\bigsqcup_{\lambda\in Y^{++} }K_I \lambda K_I\subsetneq G.\]  The set $G^+$ is a sub-semi-group of $G$. Then $\HC_\C$ is defined to be the set of functions from $K_I\backslash G^+/K_I$ to $\C$ which have finite support.  The Iwahori-Hecke algebra $\HC_\C$ of $G$ admits a Bernstein-Lusztig presentation but it is no longer indexed by $Y$. Let $Y^+=\bigcup_{w\in W^v} w.Y^{++}\subset Y$. Then $Y^+$ is the \textbf{integral Tits cone} and we have $Y^+=Y$ if and only $G$ is reductive. The \textbf{Bernstein-Lusztig-Hecke algebra of }$G$ is the space $\AC_\C=\bigoplus_{w\in W^v} \C[Y] H_w$ subject to to some relations (see subsection~\ref{subIH algebras}). Then $\HC_\C$ is isomorphic to $\bigoplus_{w\in W^v} \C[Y^+] H_w$.

Let $B^+=B\cap G^+$. Let $T_\C^+=\Hom_{\mathrm{Mon}}(Y^+,\C)\setminus\{0\}$ and $T_\C=\Hom_{\mathrm{Gr}}(Y,\C^*)$. Let $\epsilon\in \{+,\emptyset\}$. If $\tau^\epsilon\in T_\C^\epsilon$  we define the space $\widehat{I(\tau^\epsilon)^\epsilon}$  of functions $f$ from $G^\epsilon$  to $\C$ such that for every $g\in G^\epsilon$  and $b\in B^\epsilon$, 
 $f(bg)=\tau\delta^{1/2}(b)f(g)$. As we do not know which condition could replace 
``locally constant'', we do not impose any regularity condition on the functions of $\widehat{I(\tau^\epsilon)^\epsilon}$. Then $G^\epsilon$  acts by right translation on $\widehat{I(\tau^\epsilon)^\epsilon}$.
Let $I_{\tau^\epsilon,G^\epsilon}$ be the subspace of $\widehat{I(\tau^\epsilon)^\epsilon}$   of  functions which are invariant under the action of $K_I$ and whose support satisfy some finiteness conditions (see \ref{subsubWell_definedness_action}).    Inspired by formula~(\ref{eqFormula_action_IH_algebra}), we define an action of $\HC_\C$ on $I_{\tau^\epsilon,G^\epsilon}$ by \[\phi.f=\sum_{g\in G/K_I} \phi(g) g.f,\forall(\phi,f)\in \HC_\C\times I_{\tau^\epsilon,G^\epsilon}.\] 
 As often in the Kac-Moody theory, the fact that this formula is well-defined is not obvious. We prove some finiteness results on $G$ to prove that the formula only involves finite sums and that $\phi.f$ is an element of $I_{\tau^\epsilon,G^\epsilon}$ (see Definition/Proposition~\ref{def_prop_action}). 

We regard $\tau^\epsilon$ as an algebra morphism $\C[Y^\epsilon]\rightarrow \C$. Let $I_{\tau^\epsilon}^\epsilon$ be the representation of $\AC_\C^\epsilon$ (where $\AC_\C^+=\HC_\C$) defined by induction of $\tau^\epsilon$ from $\C[Y^\epsilon]$ to $\AC_\C^\epsilon$.

We prove the following proposition, which seems to indicate that the representations of $\HC_\C$ correspond to representations of $G^+$ and that the  representations of $\AC_\C$ correspond to representations of $G$:

\begin{proposition*}(see Proposition~\ref{propExtension_I_tau_G_+toG})

Let $\tau^+\in T_\C^+$. 
\begin{enumerate}

\item Suppose that $\tau^+$ is not the restriction to $Y^+$ of an element of $T_\C$. 

For every $f\in \widehat{I(\tau^+)}\setminus\{0\}$, for every $G$-module $M$, the restriction of $M$ to $G^+$ is not isomorphic to $G^+.f$. 

For every $x\in I_{\tau^+}^+\setminus\{0\}$, for every $\AC_\C$-module $M$, the restriction of  $M$ to $\HC_\C$ is not isomorphic to $\HC_\C.x$.

\item Suppose that $\tau^+$ is the restriction to $Y^+$ of a (necessarily unique) element $\tau$ of $T_\C$. 

Every element $f^+$ of $\widehat{I(\tau^+)^+}$ can be extended uniquely to an element $f$ of $\widehat{I(\tau)}$. Then $f^+\mapsto  f$ is an isomorphism of $G^+$-modules.

 The action of $\HC_\C$ on  $I_{\tau^+}^+$ extends uniquely to an action of $\AC_\C$ on $I_{\tau^+}^+$. Then $I_{\tau^+}^+$ is naturally isomorphic to $I_\tau$ as a $\AC_\C$-module.
\end{enumerate}
\end{proposition*}

Note that the existence of  elements of $T_\C^+$ which do not extend to elements of $T_\C$ depends on $G$. We prove that in some cases (for example when $G$ is affine or associated with a size $2$ Kac-Moody matrix)  every element of $T_\C^+$ is the restriction of an element of $T_\C$. We also prove that for some size $3$ Kac-Moody matrices,  there exists $\tau\in T_\C^+$ which is not the restriction of an element of $T_\C$ (see Lemma~\ref{lemNonexistence_degenerate_representations} and Lemma~\ref{lemExistence_degenerate_representations}).

We then restrict our study to the elements $\tau^+$ of $T_\C^+$ which are the restriction of an element $\tau$ of $T_\C$. We prove that $I_{\tau^+}^+$ is irreducible if and only if $I_\tau$ is (see Proposition~\ref{propEquality_submodules_H+BLH}). We then study the irreducibility of $I_\tau$.  We prove the following theorem, generalizing Matsumoto's irreducibility criterion (see Corollary~\ref{corMatsumoto theorem}): 

\begin{thm*}\label{thm*Matsumoto's criterion}
Let $\tau$ be a regular character. Then $I_\tau$ is irreducible if and only if for all $\alpha^\vee\in \Phi^\vee$, \[\tau(\alpha^\vee)\neq q.\]
\end{thm*}

We also generalize one implication of Kato's criterion (see Lemma~\ref{lemIrreducibility implies isomorphisms} and Proposition~\ref{propKato's weak theorem}). Let $\Wta$ be the subgroup of $W_\tau$ generated by the reflections $r_{\alpha^\vee}$, for $\alpha^\vee\in \Phi^\vee$ such that $\tau(\alpha^\vee)=1$.  

\begin{thm*}\label{thm*Kato's criterion}
Let $\tau\in T_\C$. Assume that $I_\tau$ is irreducible. Then:\begin{enumerate}
\item\label{itWchi reflections KM} $W_\tau=\Wta$,

\item\label{itCondition on the values of chi} for all $\alpha^\vee\in \Phi^\vee$, $\tau(\alpha^\vee)\neq q$.
\end{enumerate}
\end{thm*}

We then obtain Kato's criterion when the Kac-Moody group $G$ is associated with a size $2$ Kac-Moody matrix (see Theorem~\ref{thmKato's_theorem_dim2}):  

\begin{thm*}
Assume that $G$ is associated with a size $2$ Kac-Moody matrix. Let $\tau\in T_\C$. Then $I_\tau$ is irreducible if and only if it satisfies the following conditions: \begin{enumerate}
\item\label{itWchi engendré par ses réflexions} $W_\tau=\Wta$,

\item for all $\alpha^\vee\in \Phi^\vee$, $\tau(\alpha^\vee)\neq q$. 
 \end{enumerate}
\end{thm*}

\medskip 

In order to prove these theorems, we first establish the following irreducibility criterion. For $\tau\in T_\C$ set $I_\tau(\tau)=\{x\in I_\tau | \theta.x=\tau(\theta).x\ \forall \theta\in \C[Y]\}$. Then:
 \begin{thm*}\label{thm*Irreducibility_criterion} (see Theorem~\ref{thmIrreducibility criterion})
$I_\tau$ is irreducible if and only if:\begin{itemize}
\item $\tau(\alpha^\vee)\neq q$ for all $\alpha^\vee\in \Phi^\vee$

\item $\dim I_\tau(\tau)=1$. 
\end{itemize}
\end{thm*}

\begin{remark}
Suppose that $G$ is an affine Kac-Moody group. Then by \cite[7]{bardy2016iwahori}, some extension $\widetilde{\AC_\C}$ of $\AC_\C$ contains the double affine Hecke algebra introduced in \cite{cherednik1992double}. It would therefore be interesting to find a link between the representations of $\AC_\C$ and those of this algebra.
\end{remark}

\paragraph{Framework}
Actually, following \cite{bardy2016iwahori} we study Iwahori-Hecke algebras associated with abstract masures. In particular our results also apply when $G$ is an almost-split Kac-Moody group over a non-Archimedean local field. The definition of $\Wta$ and the  statements given in this introduction are not necessarily valid in this case and we refer to Proposition~\ref{propKato's weak theorem}, Theorem~\ref{thmKato's_theorem_dim2} and Theorem~\ref{thmIrreducibility criterion} for statements valid in this frameworks.

\paragraph{Organization of the paper}

The paper is organized as follows. In a first part (sections~\ref{secIH algebras} to \ref{secTau_simple_reflections}) we consider ``abstract'' Iwahori-Hecke algebras. We define them using the Bernstein-Lusztig presentation and they are a priori not associated with a group. The techniques used are mainly algebraic, based on the Bernstein-Lusztig relations. In a second part (section~\ref{secTowards_principal_series_representations}), we introduce Kac-Moody groups, masures and Iwahori-Hecke algebras associated with groups, and we associate some principal series representations to these groups. The techniques involved are mainly building theoretic.

In section~\ref{secIH algebras}, we recall the definition of the Iwahori-Hecke algebras and of the Bernstein-Lusztig-Hecke algebras,  introduce principal series representations and define an algebra $\ATF$ containing  $\AF$, where $\FC$ is the field of coefficients of $\AF$.

In section~\ref{secCY module and intertwining}, we study the $\FC[Y]$-module  $I_\tau$  and we study the intertwining operators from $I_\tau$ to $I_{\tau'}$, for $\tau,\tau'\in T_\FC$.

In section~\ref{secStudy of irreducibility}, we establish Theorem~\ref{thm*Irreducibility_criterion}. We then apply it to obtain Theorem~\ref{thm*Matsumoto's criterion} and Theorem~\ref{thm*Kato's criterion}.

In section~\ref{secTau_simple_reflections} we consider the weight vectors of $I_\tau$ and use them to  prove Kato's irreducibility criterion for size $2$ Kac-Moody matrices.

In section~\ref{secTowards_principal_series_representations}, we introduce Kac-Moody groups over local fields, masures, and Iwahori-Hecke algebras of these groups. We introduce some principal series representations of these groups, study them and relate them to the principal series representations studied in the previous sections.

There is an index of notations at the end of the paper.

\paragraph{Acknowledgments} I  would like to thank Ramla Abdellatif, Stéphane Gaussent and Dinakar Muthiah for the discussions we had on this topic. I also would like to thank Anne-Marie Aubert for her advice concerning references and Olivier Ta{\"i}bi for correcting the statement of the main theorem and for discussing this subject with me. I am grateful to Maarten Solleveld for his helpful corrections and comments which enabled me to simplify and improve some statements. Finally, I would like to thank the referees for their valuable comments and suggestions.

\paragraph{Funding}
The author was supported by the ANR grant ANR-15-CE40-0012.

\tableofcontents

\section{Bernstein-Lusztig presentation of Iwahori-Hecke algebras}\label{secIH algebras}
Let $G$ be a Kac-Moody group over a non-Archimedean local field. Then Gaussent and Rousseau constructed a space $\I$, called a masure on which $G$ acts, generalizing the construction of the Bruhat-Tits buildings (see \cite{gaussent2008kac}, \cite{rousseau2016groupes} and \cite{rousseau2017almost}). Rousseau then gave in \cite{rousseau2011masures} an axiomatic definition of masures inspired by the axiomatic definition of Bruhat-Tits buildings.  We call a masure satisfying these axioms an abstract masure. It is a priori not associated with any group.

In \cite{bardy2016iwahori}, Bardy-Panse, Gaussent and Rousseau attached an Iwahori-Hecke algebra $\HC_\RC$ to each abstract masure satisfying certain conditions and to each ring $\RC$. The algebra $\HC_\RC$ is an algebra of functions defined on some pairs of chambers of the masure, equipped with a convolution product. Then they prove that under some additional hypothesis on the ring $\RC$ (which are satisfied by $\R$ and $\C$),  $\HC_\RC$ admits a Bernstein-Lusztig presentation. In this section, we will only introduce the Bernstein-Lusztig presentation of $\HC_\RC$ and we do not introduce masures (we introduce them in section~\ref{secTowards_principal_series_representations}). We however introduce the standard apartment of a masure. We restrict our study to the case where $\RC=\FC$ is a field. 

\subsection{Standard apartment of a masure}\label{subRootGenSyst}

\subsubsection{Root generating system}

A \textbf{ Kac-Moody matrix} (or { generalized Cartan matrix}) is a square matrix $A=(a_{i,j})_{i,j\in I}$ indexed by a finite set $I$, with integral coefficients, and such that :
\begin{enumerate}
\item[\tt $(i)$] $\forall \ i\in I,\ a_{i,i}=2$;

\item[\tt $(ii)$] $\forall \ (i,j)\in I^2, (i \neq j) \Rightarrow (a_{i,j}\leq 0)$;

\item[\tt $(iii)$] $\forall \ (i,j)\in I^2,\ (a_{i,j}=0) \Leftrightarrow (a_{j,i}=0$).
\end{enumerate}
A \textbf{root generating system} is a $5$-tuple $\mathcal{S}=(A,X,Y,(\alpha_i)_{i\in I},(\alpha_i^\vee)_{i\in I})$\index{$\mathcal{S}$}\index{$Y$} made of a Kac-Moody matrix $A$ indexed by the finite set $I$, of two dual free $\Z$-modules $X$ and $Y$ of finite rank, and of a free family $(\alpha_i)_{i\in I}$ (respectively $(\alpha_i^\vee)_{i\in I}$) of elements in $X$ (resp. $Y$) called \textbf{simple roots} (resp. \textbf{simple coroots}) that satisfy $a_{i,j}=\alpha_j(\alpha_i^\vee)$ for all $i,j$ in $I$. Elements of $X$ (respectively of $Y$) are called \textbf{characters} (resp. \textbf{cocharacters}).

Fix such a root generating system $\mathcal{S}=(A,X,Y,(\alpha_i)_{i\in I},(\alpha_i^\vee)_{i\in I})$ and set $\A:=Y\otimes \R$\index{$\A$}. Each element of $X$ induces a linear form on $\A$, hence $X$ can be seen as a subset of the dual $\A^*$. In particular, the $\alpha_{i}$'s (with $i \in I$) will be seen as linear forms on $\A$. This allows us to define, for any $i \in I$, an involution $r_{i}$ of $\A$ by setting $r_{i}(v) := v-\alpha_i(v)\alpha_i^\vee$ for any $v \in \A$. Let $\SCC=\{r_i|i\in I\}$\index{$\SCC$} be the (finite) set of \textbf{simple reflections}.  One defines the \textbf{Weyl group of $\mathcal{S}$} as the subgroup $W^{v}$\index{$W^v$} of $\mathrm{GL}(\A)$ generated by $\SCC$. The pair $(W^{v}, \SCC)$ is a Coxeter system, hence we can consider the length $\ell(w)$ with respect to $\SCC$ of any element $w$ of $W^{v}$. If $s\in \SCC$, $s=r_i$ for some unique $i\in I$. We set $\alpha_s=\alpha_i$ and $\alpha_s^\vee=\alpha_i^\vee$.

The following formula defines an action of the Weyl group $W^{v}$ on $\A^{*}$:  
\[\displaystyle \forall \ x \in \A , w \in W^{v} , \alpha \in \A^{*} , \ (w.\alpha)(x):= \alpha(w^{-1}.x).\]
Let $\Phi:= \{w.\alpha_i|(w,i)\in W^{v}\times I\}$\index{$\Phi,\Phi^\vee$} (resp. $\Phi^\vee=\{w.\alpha_i^\vee|(w,i)\in W^{v}\times I\}$) be the set of \textbf{real roots} (resp. \textbf{real coroots}): then $\Phi$ (resp. $\Phi^\vee$) is a subset of the \textbf{root lattice} $Q := \displaystyle \bigoplus_{i\in I}\Z\alpha_i$ (resp. \textbf{coroot lattice} $Q^\vee=\bigoplus_{i\in I}\Z\alpha_i^\vee$). By \cite[1.2.2 (2)]{kumar2002kac}, one has $\R \alpha^\vee\cap \Phi^\vee=\{\pm \alpha^\vee\}$ and $\R \alpha\cap \Phi=\{\pm \alpha\}$ for all $\alpha^\vee\in \Phi^\vee$ and $\alpha\in \Phi$.

\subsubsection{Fundamental chamber, Tits cone and vectorial faces}

As in the reductive case, define the \textbf{fundamental chamber} as $C_{f}^{v}:= \{v\in \A \ \vert \ \forall s \in \SCC,  \alpha_s(v)>0\}$\index{$C^f_v$}. 

 Let $\mathcal{T}:= \displaystyle \bigcup_{w\in W^{v}} w.\overline{C^{v}_{f}}$\index{$\T$} be the \textbf{Tits cone}. This is a convex cone (see \cite[1.4]{kumar2002kac}).
 
 For $J\subset \SCC$, set $F^v(J)=\{x\in \A| \alpha_j(x)=0\forall j\in J\text{ and }\alpha_j(x)>0 \forall j\in \SCC\setminus J\}$\index{$F^v(J)$}. A \textbf{positive vectorial face} (resp. \textbf{negative}) is a set of the form $w.F^v(J)$ ($-w.F^v(J)$) for some $w\in W^v$ and $J\subset \SCC$. Then by \cite[5.1 Th\'eor\`eme (ii)]{remy2002groupes}, the family of positive vectorial faces of $\A$ is a partition of $\T$ and the stabilizer of $F^v(J)$ is $W_J=\langle J\rangle$. 

One sets $Y^{++}=Y\cap\overline{C^v_f}$\index{$Y^+$, $Y^{++}$} and $Y^+=Y\cap \T$.

\begin{remark}
By \cite[§4.9]{kac1994infinite} and \cite[§ 5.8]{kac1994infinite} the following conditions are equivalent:\begin{enumerate}
\item the Kac-Moody matrix $A$ is of finite type (i.e. is a Cartan matrix),

\item $\A=\T$

\item $W^v$ is finite.
\end{enumerate}
\end{remark}

\subsection{Recollections on Coxeter groups }\label{subReflection_subgroups}

\subsubsection{Bruhat order}

Let $(W_0,\SCC_0)$ be a Coxeter system. We equip it with the Bruhat order $\leq_{W_0}$ (see \cite[Definition 2.1.1]{bjorner2005combinatorics}). We have the following characterization (see \cite[Corollary 2.2.3]{bjorner2005combinatorics}): let $u,w\in W_0$. Then $u\leq_{W_0} w$ if and only if every reduced expression for $w$ has a subword that is a reduced
expression for $u$. By \cite[Proposition 2.2.9]{bjorner2005combinatorics}, $(W_0,\leq_{W_0})$ is a \textbf{directed poset}, i.e for every finite set $E\subset W_0$, there exists $w\in W_0$ such that $v\leq_{W_0} w$ for all $v\in E$. 

We write $\leq$ instead of $\leq_{W^v}$. For $u,v\in W^v$, we denote by $[u,v]$, $[u,v)$, $\ldots$ the sets $\{w\in W^v|u\leq w\leq v\}$, $\{w\in W^v|u\leq w <v\}$, $\ldots$. 

\subsubsection{Reflections  and coroots}

Let $\RCC=\{wsw^{-1}|w\in W^v, s\in \SCC\}$\index{$\RCC$} be the set of \textbf{reflections} of $W^v$. Let  $r\in \RCC$. Write $r=wsw^{-1}$, where $w\in W^v$, $s\in \SCC$ and $ws>w$ (which is possible because if $ws<w$, then $r=(ws)s(ws)^{-1}$). Then one sets $\alpha_r=w.\alpha_s\in \Phi_+$\index{$\alpha_r,\alpha_r^\vee$} (resp. $\alpha_r^\vee=w.\alpha_s^\vee\in\Phi^\vee_+$). This is well-defined by the lemma below.

\begin{lemma}\label{lemDefinition_root_reflection}
Let $w,w'\in W^v$ and $s,s'\in \SCC$ be such that $wsw^{-1}=w's'w'^{-1}$ and $ws>w$, $w's'>w'$. Then $w.\alpha_s=w'.\alpha_{s'}\in \Phi_+$ and $w.\alpha_s^\vee=w'.\alpha_{s'}^\vee\in \Phi^\vee_+$.
\end{lemma}

\begin{proof}
One has $r(x)=x-w.\alpha_s(x)w.\alpha_s^\vee=x-w'.\alpha_{s'}(x)w'.\alpha_{s'}^\vee$ for all $x\in \A$ and thus $w.\alpha_s\in \R^* w'.\alpha_{s'}$ and $w.\alpha_{s}^\vee\in \R^* w'.\alpha_{s'}^\vee$. As $\Phi$ and $\Phi^\vee$ are reduced, $w.\alpha_s=\pm w'.\alpha_{s'}$ and $w.\alpha_s^\vee=\pm w'.\alpha_s^\vee$. By \cite[Lemma 1.3.13]{kumar2002kac}, $w.\alpha_s,w'.\alpha_{s'}\in \Phi_+$ and $w.\alpha_s^\vee, w'.\alpha_{s'}^\vee\in \Phi^\vee_+$, which proves the lemma.
\end{proof}

\begin{lemma}\label{lemKumar1.3.11}
Let $r,r'\in \RCC$ and $w\in W^v$ be such that $w.\alpha_r=\alpha_{r'}$ or $w.\alpha_r^\vee=\alpha_{r'}^\vee$. Then $wrw^{-1}=r'$.
\end{lemma}

\begin{proof}
Write $r=vsv^{-1}$ and $r'=v's'v'^{-1}$ for $s,s'\in \SCC$ and $v,v'\in W^v$. Then $v'^{-1}wv.\alpha_s=\alpha_{s'}$. Thus by \cite[Theorem 1.3.11 (b5)]{kumar2002kac}, $v'^{-1}wvsv^{-1}w^{-1}v'=s'$ and hence $wrw^{-1}=r'$.
\end{proof}

Let $r\in \RCC$. Then for all $x \in \A$, one has: \[r(x)=x-\alpha_r(x)\alpha_r^\vee.\]  Let $\alpha^\vee\in \Phi^\vee$. One sets $r_{\alpha^\vee}=wsw^{-1}$\index{$r_{\alpha^\vee}$} where $(w,s)\in W^v\times \SCC$ is such that $\alpha^\vee=w.\alpha_s^\vee$. This is well-defined, by Lemma~\ref{lemKumar1.3.11}. Thus $\alpha^\vee\mapsto r_{\alpha^\vee}$ and $r\mapsto \alpha_r^\vee$ induce bijections $\Phi^\vee_+\rightarrow \RCC$ and $\RCC\rightarrow \Phi^\vee_+$.  If $r\in \RCC$, $r=wsw^{-1}$, one sets $\sigma_r=\sigma_s$, which is well-defined by assumption on the $\sigma_t$, $t\in \SCC$ (see Subsection~\ref{subIH algebras}).

For $w\in W^v$, set $N_{\Phi^\vee}(w)=\{\alpha^\vee\in \Phi^\vee_+|w.\alpha^\vee\in \Phi^\vee_-\}$\index{$N_{\Phi^\vee}(w)$}.

\begin{lemma}\label{lemKumar_1.3.14}(\cite[Lemma 1.3.14]{kumar2002kac})
Let $w\in W^v$.  Then $|N_{\Phi^\vee}(w)|=\ell(w)$ and if $w=s_1\ldots s_r$ is a reduced expression, then $N_{\Phi^\vee}(w)=\{\alpha_{s_r}^\vee,s_r.\alpha_{s_{r-1}}^\vee,\ldots,s_r\ldots s_{2}.\alpha_{s_1}^\vee\}$. 
\end{lemma}

\subsubsection{Reflections subgroups of a Coxeter group}\label{subsubReflections_subgroups}

If $W_0$ is a Coxeter group, a \textbf{Coxeter generating set} is a set $\SCC_0$ such that $(W_0,\SCC_0)$ is a Coxeter system.   Let $(W_0,\SCC_0)$ be a Coxeter system and $\RCC_0=\{w.s.w^{-1}|w\in W_0,s\in \SCC_0\}$ be its set of reflections. A \textbf{reflection subgroup of} $W_0$ is a group of the form $W_1=\langle \RCC_1\rangle$ for some $\RCC_1\subset \RCC_0$. For $w\in W_0$, set $N_{\RCC_0}(w)=\{r\in \RCC_0|rw^{-1}<w^{-1} \}$\index{$N_{\RCC_0}(w)$}. By \cite[3.3]{dyer1990reflection} or \cite[1]{dyer1991bruhat}, if $\SCC(W_1)=\{r\in \RCC_0|N_{\RCC_0}(r)\cap W_1=\{r\}\}$, then $(W_1,\SCC(W_1))$ is a Coxeter system.

Let $(W_0,\SCC_0)$ be a Coxeter system. The \textbf{rank of }$(W_0,\SCC_0)$ is $|\SCC_0|$.

\begin{remark}\label{rkCoxeter_generators}
\begin{enumerate}
\item\label{itCoxeter_generators_infinite_dihedral_group}  The rank of a Coxeter group is not well-defined. For example, by \cite[3]{muhlherr2005isomorphism}, if $k\in \Ne$ and $n=4(2k+1)$ then the dihedral group of order $n$ admits Coxeter generating sets of order $2$ and $3$. However by \cite{radcliffe1999rigidity},  all the Coxeter generating sets of the infinite dihedral group have cardinal $2$.

\item  Using \cite[IV 1.8  Proposition 7]{bourbaki1981elements} we can prove that if $(W_0,\SCC_0)$ is a Coxeter system of infinite rank, then every Coxeter generating set of $W_0$ is infinite. 

\item Reflection subgroups of finite rank Coxeter groups are not necessarily of finite rank. Indeed, let $W_0$ be the Coxeter group generated by the involutions $s_1,s_2,s_3$, with $s_is_j$ of infinite order when $i\neq j\in \llbracket 1,3\rrbracket$. Let $W_0'=\langle s_1,s_2\rangle\subset W_0$ and $\RCC_1=\{ws_3w^{-1}|w\in W_0'\}\subset \RCC_0$. Then $W_1=\langle \RCC_1 \rangle$ has infinite rank. Indeed, let $\psi:W_0\rightarrow W_0'$ be the group morphism defined by $\psi_{|W'_0}=\Id_{W'_0}$ and $\psi(s_3)=1$. Then $\RCC_1\subset \ker \psi$. Thus $s_3$ appears in the reduced writing of every nontrivial element of $W_1$. By \cite[Corollary 1.4.4]{bjorner2005combinatorics} if  $r\in \RCC_1$, then the unique element of $N_{\RCC_0}(r)$ containing an $s_3$ in its reduced writing is $r$. Thus $\SCC(W_1)\supset \RCC_1$ is infinite.
\end{enumerate}

\end{remark}

\subsection{Iwahori-Hecke algebras}\label{subIH algebras}

In this subsection, we give the definition of  the Iwahori-Hecke algebra via its Bernstein-Lusztig presentation, as done in \cite[Section 6.6]{bardy2016iwahori}. 

Let  $\RC_1=\Z[(\sigma_s)_{s\in \SCC},(\sigma'_s)_{s\in \SCC}]$, where $(\sigma_s)_{s\in \SCC}, (\sigma'_s)_{s\in \SCC}$ are two families of indeterminates satisfying the following relations:
\begin{itemize}
\item if $\alpha_{s}(Y) = \Z$, then $\sigma_{s} = \sigma'_{s}$\index{$\sigma_s,\sigma_s'$};  
\item if $s,t \in \SCC$ are such that the order of $st$ is finite and odd (i.e if $\alpha_s(\alpha_t^\vee)=\alpha_t(\alpha_s^\vee)=-1$), then $\sigma_{s}=\sigma_{t}=\sigma'_{s}=\sigma'_{t}$.
\end{itemize}

To define the Iwahori-Hecke algebra $\HH_{\RC_1}$ associated with $\A$ and $(\sigma_s, \sigma'_s)_{s \in \SCC}$, we first introduce the Bernstein-Lusztig-Hecke algebra. Let $^{BL}\HH_{\RC_1}$ be the free $\RC_1$-vector space with basis $(Z^\lambda H_w)_{\lambda\in Y,w\in W^v}$. For short, one sets $H_{w} = Z^{0}H_{w}$\index{$H_w$} for $w \in W^{v}$ and $Z^{\lambda} = Z^{\lambda}H_{1}$\index{$Z^\lambda$} for $\lambda \in Y$. The \textbf{Bernstein-Lusztig-Hecke algebra} $^{BL}\HH_{\RC_1}$ is the module $^{BL}\HH_{\RC_1}$ equipped with the unique product $*$ that turns it into an associative algebra and satisfies the following relations (known as the \textbf{Bernstein-Lusztig relations}):
\begin{itemize}
\item (BL1) $\forall \ (\lambda, w) \in Y \times W^{v}$, $Z^{\lambda} * H_{w} = Z^{\lambda}H_{w}$;
\item (BL2) $\forall \ s \in \SCC, \forall \ w \in W^{v}$, $H_{s}*H_{w}=\left\{\begin{aligned} & H_{sw} &\mathrm{\ if\ }\ell(sw)=\ell(w)+1\\ & (\sigma_{s}-\sigma_{s}^{-1}) H_{w}+H_{s w} &\mathrm{\ if\ }\ell(sw)=\ell(w)-1 \end{aligned}\right . \ $;
\item (BL3) $\forall \ (\lambda, \mu) \in Y^{2}$, $Z^{\lambda} * Z^{\mu} = Z^{\lambda + \mu}$;
\item (BL4) $\forall \ \lambda \in Y, \ \forall \ i \in I$, $H_{s}*Z^{\lambda}-Z^{s.\lambda}*H_{s} =Q_s(Z)(Z^\lambda-Z^{s.\lambda})$, where $Q_s(Z)=\frac{(\sigma_s-\sigma_s^{-1})+(\sigma_s'-\sigma_s'^{-1})Z^{-\alpha_s^\vee}}{1-Z^{-2\alpha_s^\vee}}$\index{$Q_s(Z)$}. 
\end{itemize}
The existence and uniqueness of such a product $*$ comes from \cite[Theorem 6.2]{bardy2016iwahori}. 

\begin{definition}\label{defBernstein-Lusztig algebra}
Let $\FC$ be a field of characteristic $0$  and $f:\RC_1\rightarrow \FC$ be a ring morphism such that $f(\sigma_s)$ and $f(\sigma_s')$ are invertible in $\FC$ for all $s\in \SCC$. Then the \textbf{Bernstein-Lusztig-Hecke algebra of } $(\A,(\sigma_s)_{s\in \SCC},(\sigma'_s)_{s\in \SCC})$ over $\FC$ is the algebra $\AC_\FC=\AC_{\RC_1}\otimes_{\RC_1}\FC$\index{$\AC_\FC$}.  Following \cite[Section 6.6]{bardy2016iwahori}, the \textbf{Iwahori-Hecke algebra} $\HH_\FC$\index{$\HC_\FC$} associated with $\mathcal{S}$ and $(\sigma_s, \sigma'_s)_{s \in \SCC}$ is now defined as the $\FC$-subalgebra of $^{BL}\HH_\FC$ spanned by $(Z^{\lambda} H_{w})_{\lambda\in Y^{+}, w\in W^{v}}$ (recall that $Y^{+} = Y \cap \T$ with $\T$ being the Tits cone). Note that for $G$ reductive, we recover the usual Iwahori-Hecke algebra of $G$, since $Y\cap \T=Y$.
\end{definition}

In certain proofs, when $\FC=\C$, we will make additional assumptions on the $\sigma_s$ and $\sigma_s'$, $s\in \SCC$. To avoid these assumptions, we can assume that $\sigma_s,\sigma_s'\in \C$ and $|\sigma_s|>1, |\sigma_s'|>1$ for all $s\in \SCC$.

\begin{remark}\label{remIH algebre dans le cas KM deploye}

\begin{enumerate}
\item Let $s\in \SCC$. Then if $\sigma_s=\sigma'_s$, $Q_s(Z)=\frac{(\sigma_s-\sigma_s^{-1})}{1-Z^{-\alpha_s^\vee}}$.

\item\label{itPolynomiality_Bernstein_Lusztig} Let $s\in \SCC$ and $\lambda\in Y$. Then $Q_s(Z)(Z^\lambda-Z^{s.\lambda})\in \FC[Y]$. Indeed, $Q_s(Z)(Z^\lambda-Z^{s.\lambda})=Q_s(Z).Z^\lambda(1-Z^{-\alpha_s(\lambda)\alpha_s^\vee})$. Assume that $\sigma_s=\sigma_s'$. Then \[ \frac{1-Z^{-\alpha_s(\lambda)\alpha_s^\vee}}{1-Z^{-\alpha_s^\vee}}=\left\{\begin{aligned} &\sum_{j=0}^{\alpha_s(\lambda)-1}Z^{-j\alpha_s^\vee}  &\mathrm{\ if\ }\alpha_s(\lambda)\geq 0\\ & -Z^{\alpha_s^\vee}\sum_{j=0}^{-\alpha_s(\lambda)-1}Z^{j\alpha_s^\vee}  &\mathrm{\ if\ }\alpha_s(\lambda)\leq 0,\end{aligned}\right.\] and thus $Q_s(Z)(Z^\lambda-Z^{s.\lambda})\in \FC[Y]$. Assume $\sigma_s'\neq \sigma_s$. Then $\alpha_s(Y)=2\Z$ and a similar computation enables to conclude.

\item From (BL4) we deduce that for all $s\in \SCC$, $\lambda\in Y$, \[Z^\lambda* H_s-H_s*Z^{s.\lambda}=Q_s(Z)(Z^\lambda-Z^{s.\lambda}).\]

\item When $G$ is a split Kac-Moody group over a non-Archimedean local field $\mathcal{K}$ with residue cardinal $q$, we can choose $\FC$ to be a field containing $\Z[\sqrt{q}^{\pm 1}]$ and  take $f(\sigma_s)=f(\sigma'_s)=\sqrt{q}$ for all $s\in \SCC$.

\item By (BL4), the family $(H_w*Z^\lambda)_{w\in W^v, \lambda\in Y}$ is also a basis of $\AC_\FC$.

\item\label{itWell_definition_sigma_w} Let $w\in W^v$ and $w=s_1\ldots s_k$, with $k\in \N$ and $s_1,\ldots,s_k\in  \SCC$ be a reduced expression of $w$. We set $\sigma_w=\sigma_{s_1}\ldots \sigma_{s_k}$. This is well-defined, independently of the choice of a reduced expression of $w$ by the conditions imposed on the $\sigma_s$ and by \cite[Theorem 3.3.1 (ii)]{bjorner2005combinatorics}.

\end{enumerate}

\end{remark}

We equip $\FC[Y]$ with an action of $W^v$. For $\theta=\sum_{\lambda\in Y}{a_\lambda Z^\lambda}\in \FC[Y]$ and $w\in W^v$, set $\theta^w:=\sum_{\lambda\in Y}a_\lambda Z^{w.\lambda}$. 

\begin{lemma}\label{lemCommutation relation}
Let $\theta\in \FC[Y]$ and $w\in W^v$. Then $\theta* H_w -H_w* \theta^{w^{-1}}\in \AF^{<w}:=\bigoplus_{v<w} H_v\FC[Y]$.  In particular, $\AF^{\leq w}:=\bigoplus_{v\leq w}H_v\C[Y]$ is a left finitely generated $\FC[Y]$-submodule of $\AF$.
\end{lemma}

\begin{proof}
We do it by induction on $\ell(w)$. Let $\theta \in \FC[Y]$ and $w\in W^v$ be such that $u:=\theta H_w-H_w \theta^{w^{-1}}\in \AC(T_\FC)^{<w}$. Let $s\in \SCC$ and assume that $\ell(ws)=\ell(w)+1$. Then by (BL4): \[\theta* H_{ws}=(H_w \theta^{w^{-1}}+u)*H_s=H_{ws}\theta^{sw^{-1}}+aH_w+uH_s,\] for some $a\in \FC$. Moreover, by \cite[Corollary 1.3.19]{kumar2002kac} 	and (BL2),  $u*H_s\in \AC(T_\FC)^{< ws}$ and the lemma follows.
\end{proof}

\begin{definition}\label{defHecke algebra of W}
Let $\HC_{\FC,W^v}=\bigoplus_{w\in W^v} \FC H_w\subset \HC_{\FC}$\index{$\HC_{\FC,W^v}$}. Then $\HC_{\FC,W^v}$ is a subalgebra of  $\HC_\FC$. This is the Hecke algebra of the Coxeter group $(W^v,\SCC)$. 
\end{definition}

\subsection{Principal series representations}\label{subPrincipal series representations}
In this subsection, we introduce the principal series representations of $\AF$. 

We now fix $(\A,(\sigma_s)_{s\in \SCC},(\sigma'_s)_{s\in \SCC})$ as in Subsection~\ref{subIH algebras} and a field $\FC$ as in Definition~\ref{defBernstein-Lusztig algebra}. Let $\HC_\FC$ and $\AC_\FC$ be the Iwahori-Hecke and the Bernstein-Lusztig Hecke algebras of 

$(\A,(\sigma_s)_{s\in \SCC},(\sigma'_s)_{s\in \SCC})$  over $\FC$.

Let  $T_\FC= \Hom_{\mathrm{Gr}}(Y,\FC^\times)$\index{$T_\FC$} be the group of homomorphisms from $Y$ to $\FC^*$. Let $\tau\in T_\FC$. Then $\tau$ induces an algebra morphism $\tau:\FC[Y]\rightarrow \FC$ by the formula $\tau(\sum_{y\in Y} a_y e^y)=\sum_{y\in Y} a_y \tau(y)$, for $\sum a_y e^y\in \FC[Y]$. This equips $\FC$ with the structure of a  $\FC[Y]$-module.

Let $I_\tau=\mathrm{Ind}^{\AC_\FC}_{\FC[Y]}(\tau)=\AC_\FC\otimes_{\FC[Y]} \FC$\index{$I_\tau$}.  For example if $\lambda\in Y$, $w\in W^v$  and $s\in \SCC$, one has: \[Z^\lambda.1\otimes_\tau 1=\tau(\lambda)1\otimes_\tau 1, H_w*Z^\lambda \otimes_\tau 1=\tau(\lambda)H_w\otimes_\tau 1\text{ and }\]  \[ Z^\lambda.H_s\otimes_\tau 1= H_s*Z^{s.\lambda}\otimes_\tau 1+Q_s(Z)(Z^\lambda-Z^{s.\lambda})\otimes_\tau 1=\tau(s.\lambda)H_i\otimes_\tau 1+\tau\big(Q_s(Z)(Z^\lambda-Z^{s.\lambda})\big)\otimes_\tau 1.\]

 Let $h\in I_\tau$. Write $h=\sum_{\lambda\in Y, w\in W^v} h_{w,\lambda} H_w Z^\lambda\otimes_{\tau} c_{w,\lambda}$, where $(h_{w,\lambda}),(c_{w,\lambda})\in \FC^{(W^v\times Y)}$, which is possible by Remark~\ref{remIH algebre dans le cas KM deploye}. Thus \[h=\sum_{\lambda\in Y,w\in W^v} h_{w,\lambda} c_{w,\lambda} \tau(\lambda)H_w \otimes_{\tau} 1=\big(\sum_{\lambda\in Y,w\in W^v} h_{w,\lambda} c_{w,\lambda} \tau(\lambda)H_w\big)1\otimes_{\tau} 1.\] Thus $I_\tau$ is a principal $\AC_\FC$-module and $(H_w\otimes_\tau 1)_{w\in W^v}$ is a basis of $I_\tau$. Moreover $I_\tau=\HC_{W^v,\FC}.1\otimes_\tau 1$ (see Definition~\ref{defHecke algebra of W} for the definition of $\HC_{W^v,\FC}$).

 The definition of principal series representations of $\HC_\FC$ is very similar: we replace $T_\FC$ by $T_\FC^+=\Hom_{\mathrm{Mon}}(Y^+,\C)\setminus\{0\}$ and $\FC[Y]$ by $\FC[Y^+]$ in the definition above. If $\tau\in T_\FC^+$, we denote by $I_{\tau^+}^+$ the principal series representation of $\HC_\FC$ associated with $\tau^+$.

\begin{remark}\label{rkBruhat order and CY structure}
Let $\tau\in T_\FC$. By Lemma~\ref{lemCommutation relation}, $I_\tau^{\leq w}$ and $I_\tau^{\ngeq w}$ are $\FC[Y]$-submodules of $I_\tau$. In particular    $\FC[Y].x$ is finite dimensional for all $x\in I_\tau$.
\end{remark}

\begin{lemma}\label{lemEquality_submodules_H_BLH}
Let $\tau\in T_\FC$. Let $M\subset I_\tau$  be a finite dimensional $\FC[Y^+]$-submodule of $I_\tau$. Then $M$ is an $\FC[Y]$-submodule of $I_\tau$.
\end{lemma}

\begin{proof}
Let $\lambda\in Y^+$. Let $\phi_\lambda:M\rightarrow M$ be defined by $\phi_\lambda(x)=Z^\lambda.x$, for all $m\in M$. Let $x\in \mathrm{ker}(\phi_\lambda)$. Then $Z^{-\lambda}.Z^\lambda.x=0=x$ and thus $\phi_\lambda$ is an isomorphism. Moreover, $\phi_\lambda^{-1}(x)=Z^{-\lambda}.x$ for all $x\in M$ and thus $Z^{-\lambda}.x\in M$, for all $x\in M$. As $Y^+-Y^-=Y$, we deduce the lemma.
\end{proof}

\begin{proposition}\label{propEquality_submodules_H+BLH}
 Let $\tau\in T_\FC$ and $M\subset I_\tau$. Then $M$ is an $\HC_\FC$-submodule of $I_\tau$ if and only if $M$ is a ${^{BL}\HC}_\FC$-submodule of $I_\tau$. In particular, $I_\tau$ is irreducible as a $\AC_\FC$-module if and only if $I_\tau$ is irreducible as an $\HC_\FC$-module.
\end{proposition}

\begin{proof}
Let $M\subset I_\tau$ be a $\HC_\FC$-submodule. Then $M$ is an $\FC[Y^+]$ submodule of $I_\tau$. Let $x\in M$. Then by Remark~\ref{rkBruhat order and CY structure}, $\FC[Y^+].x\subset \FC[Y].x$ is finite dimensional. Thus $M=\sum_{x\in M} \FC[Y^+].x$ and by Lemma~\ref{lemEquality_submodules_H_BLH}, $M$ is an $\FC[Y]$-submodule of $I_\tau$. As $\AC_\FC$ is generated as an algebra by $\HC_\FC$ and $\FC[Y]$, we deduce the proposition.
\end{proof}

\subsection{The algebra \texorpdfstring{${^\mathrm{BL}\HC_\FC}(T_\FC)$}{(BL)H(TC)}}

In this subsection, we introduce an  algebra $\ATF$ containing $\AF$. This algebra will enable us to regard the elements of $I_\tau$ as specializations at $\tau$ of certain elements of $\ATF$. When $\FC=\C$, this will enable us to make $\tau\in T_\C$ vary and to use density arguments and basic algebraic geometry to study the $I_\tau$.

\subsubsection{Description of $\ATF$}

Let $\ATF$\index{$\ATF$} be the right $\FC(Y)$ vector space $\bigoplus_{w\in W^v} H_w\FC(Y)$.  We equip $\FC(Y)$ with an action of $W^v$. For $\theta=\frac{\sum_{\lambda\in Y}{a_\lambda Z^\lambda}}{\sum_{\lambda\in Y} b_\lambda Z^\lambda}\in \FC(Y)$ and $w\in W^v$, set $\theta^w:=\frac{\sum_{\lambda\in Y}a_\lambda Z^{w.\lambda}}{\sum_{\lambda\in Y}b_\lambda Z^{w.\lambda}}$.

\begin{proposition}\label{propDefinition_ATF}
There exists a unique multiplication $*$ on $\ATF$ which equips $\ATF$ with the structure of an associative algebra and   such that: \begin{itemize}
\item $\FC(Y)$ embeds into $\ATF$ as an algebra,

\item (BL2) is satisfied,

\item the following relation (BL4') is satisfied: \[\text{for all }\theta\in \FC(Y)\text{ and }s\in \SCC, \theta* H_s-H_s*\theta^s=Q_s(Z)(\theta-\theta^s).\]
\end{itemize}
\end{proposition}

The proof of this proposition is postponed to~\ref{subsubConstruction_of_ATF}.

We regard the elements of $\FC[Y]$ as polynomial functions on $T_\FC$ by setting: \[\tau(\sum_{\lambda\in Y}a_\lambda Z^\lambda)=\sum_{\lambda\in Y}a_\lambda\tau(\lambda),\] for all $(a_\lambda)\in \FC^{(Y)}$. The ring $\FC[Y]$ is a unique factorization domain. Let $\theta\in \FC(Y)$ and $(f,g)\in \FC[Y]\times \FC[Y]^*$ be such that $\theta=\frac{f}{g}$ and $f$ and $g$ are coprime. Set $\DC(\theta)=\{\tau\in T_\FC|\theta(g)\neq 0\}$\index{$\DC(\theta)$}. Then we regard  $\theta$ as a map from $\DC(\theta)$ to $\FC$ by setting $\theta(\tau)=\frac{f(\tau)}{g(\tau)}$ for all $\tau\in \DC(\theta)$.

 For $w\in W^v$, let $\pi^H_w:\ATF\rightarrow \FC(Y)$\index{$\pi^H_w$} be defined by $\pi^H_w(\sum_{v\in W^v} H_v\theta_v)=\theta_w$, for $(\theta_v)\in (\HFW)^{W^v}$ with finite support. If $\tau\in T_\FC$, let $\FC(Y)_\tau=\{\frac{f}{g}|f,g\in \C[Y]\text{ and } g(\tau)\neq 0\}\subset \FC(Y)$\index{$\FC(Y)_\tau$}. Let $\ATF_\tau=\bigoplus_{w\in W^v} H_w \FC(Y)_\tau \subset \ATF$. This is a not a subalgebra of $\ATF$ (consider for example $\frac{1}{Z^\lambda-1}*H_s=H_s*\frac{1}{Z^{s.\lambda}-1}+\ldots$ for some well chosen $\lambda\in Y$, $s\in \SCC$ and $\tau\in T_\C$). It is however an $\HFW-\FC(Y)_\tau$ bimodule. For $\tau\in T_\FC$, we define $\ev_\tau:\ATF_\tau\rightarrow \HFW$\index{$ev_\tau$} by $\ev_\tau(h)=h(\tau)=\sum_{w\in W^v} H_w\theta_w(\tau)$ if $h=\sum_{w\in W^v}H_w \theta_w\in \HC(Y)_\tau$. This is a morphism of $\HFW-\FC(Y)_\tau$-bimodule.

\subsubsection{Construction of $\ATF$}\label{subsubConstruction_of_ATF}

We now prove the existence of $\ATF$. For this we use the theory of Asano and Ore of rings of fractions: $\ATF$ will be the ring $\AF*(\FC[Y]\setminus\{0\})^{-1}$.

Let $V=\AC_\FC\otimes_{\FC[Y]} \FC(Y)\supset \AC_\FC$, where $\AC_\FC$ is equipped with its structure of a right $\FC[Y]$-module. As a right $\FC(Y)$-vector space, $V=\bigoplus_{w\in W^v} H_w\FC(Y)$. The left action of $\FC[Y]$ on $\AC_\FC$ extends  to an action of $\FC[Y]$ on $V$ by setting $\theta.\sum_{w\in W^v} H_w f_w=\sum_{w\in W^v} (\theta.H_w)f_w$, for $\theta\in\FC[Y]$ and $(f_w)\in \FC(Y)^{W^v}$ with finite support. This equips $V$ with the structure of a $(\FC[Y]-\FC(Y))$-bimodule.

\begin{lemma}\label{lemDefinition_action_F(Y)_ATF}
The left action of $\FC[Y]$ on $V$ extends uniquely to a left action of $\FC(Y)$ on $V$. This equips $V$ with the structure of a $(\FC(Y)$-$\FC(Y))$-bimodule.  \end{lemma}

\begin{proof}
Let $w\in W^v$ and $P\in \FC[Y]\setminus\{0\}$.  Let $V^{\leq w}=\bigoplus_{v\in [1, w]} H_v\FC(Y)$.  By Lemma~\ref{lemCommutation relation}, the map $m_P:V^{\leq w}\rightarrow V^{\leq w}$ defined by $m_P(h)=P.h$ is well-defined. Thus the left action of $\FC[Y]$ on $V^{\leq w}$ induces a ring morphism $\phi_w:\FC[Y]\rightarrow \mathrm{End}_{v.s}(V^{\leq w})$, where $\mathrm{End}_{v.s}(V^{\leq w})$ is the  space of endomorphisms of the  $\FC(Y)$-vector space $V^{\leq w}$. 

Let us prove that $\phi_w(P)$ is injective. Let $h\in V^{\leq w}$. Write $h=\sum_{v\in [1,w]} H_v \theta_v$, with $\theta_v\in \FC(Y)$ for all $v\in [1,w]$. Suppose that $h\neq 0$. Let $v\in [1,w]$ be such that $\theta_v\neq 0$ and such that $v$ is maximal for this property for the Bruhat order. By Lemma~\ref{lemCommutation relation}, $P*h\neq 0$ and thus $\phi_w(P)$ is injective. As $V^{\leq w}$ is finite dimensional over $\FC(Y)$, we deduce that $\phi_w(P)$ is invertible for all $P\in \FC[Y]$. Thus $\phi_w$ extends uniquely to a ring morphism $\widetilde{\phi_w}:\FC(Y)\rightarrow V^{\leq w}$. As $(W^v,\leq)$ is a directed poset, there exists an increasing sequence $(w_n)_{n\in \N}$ (for the Bruhat order) such that $\bigcup_{n\in \N} [1,w_n]=W^v$. Let $m,n\in \N$ be such that $m\leq n$. Let $P\in \FC[Y]$ and $f^{(m)}=\widetilde{\phi_{w_m}}(P)$ and $f^{(n)}=\widetilde{\phi_{w_n}}(P)$. Then $f^{(n)}_{|V^{\leq w_m}}=f^{(m)}$ and thus for all $\theta\in \FC(Y)$ and $x\in \AC(T_\FC)$, $\theta.x:= \widetilde{\phi}_{w_k}(\theta)(x)$ is well-defined, independently of $k\in \N$ such that $x\in V^{\leq w_k}$. This defines an action of $\FC(Y)$ on $V$. 

Let $h\in V$, $\theta\in \FC(Y)$ and $P\in \FC[Y]\setminus\{0\}$. Let $x=\frac{1}{P}.h$. Then as $V$ is a $(\FC[Y]$-$\FC(Y))$-bimodule, $(P*x)*\theta=h*\theta=P*(x*\theta)$ and thus $x*\theta=\frac{1}{P}*(h*\theta)=(\frac{1}{P}*h)*\theta$. Thus $V$ is a $(\FC(Y)-\FC(Y))$-bimodule. 
\end{proof}

\begin{lemma}\label{lemOreCondition}
The set $\FC[Y]\subset \AF$ satisfies the right Ore condition: for all $P\in \FC[Y]\setminus\{0\}$ and $h\in \AF\setminus\{0\}$, $P*\AF\cap h*\FC[Y]\neq \{0\}$.
\end{lemma}

\begin{proof}
Let $P\in \FC[Y]\setminus\{0\}$ and $h\in \AF\setminus\{0\}$. Then by definition, $P*(\frac{1}{P}*h)=h\in V$. Moreover,  $V=\bigoplus_{w\in W^v}H_w\FC(Y)$ and thus  there exists $\theta\in \FC[Y]\setminus\{0\}$ such that $ \frac{1}{P}*h*\theta\in \AF\setminus\{0\}$. Then $P*\frac{1}{P}*h*\theta=h*\theta\in P*\AF\cap  h*\FC[Y]$, which proves the lemma.
\end{proof}

\begin{definition}
Let $R$ be a ring and $r$ in $R$. Then $r$ is said to be \textbf{regular} if for all $r'\in R\setminus \{0\}$, $rr'\neq 0$ and $r'r\neq 0$.

Let $R$ be a ring and  $X \subset R$ a multiplicative set of regular
elements. A \textbf{right ring of fractions for $R$ with
respect to $X$} is any overring $S \supset R$ such that:\begin{itemize}
\item Every element of $X$ is invertible in $S$.

\item  Every element of $S$ can be expressed in the form $ax^{-1}$ for some $a \in R$
and $x\in X$.
\end{itemize}

\end{definition}

 We can now prove Proposition~\ref{propDefinition_ATF}. The uniqueness of such a product follows from (BL4'). By Lemma~\ref{lemCommutation relation}, the elements of $\FC[Y]\setminus\{0\}$ are regular. By Lemma~\ref{lemOreCondition} and \cite[Theorem 6.2]{goodearl2004introduction}, there exists a right ring of fractions $\ATF$ for $\AF$  with respect to $\FC[Y]\setminus\{0\}$. Then $\ATF$ is an algebra over $\FC$ and as a vector space, $\ATF=\bigoplus_{w\in W^v} (H_w \FC[Y])(\FC[Y]\setminus\{0\})^{-1}=\bigoplus_{w\in W^v} H_w\FC(Y)$.
 
 Let $(f,g)\in \FC[Y]\times(\FC[Y]\setminus\{0\})$. Then it is easy to check that $g*\big(H_s*\frac{1}{g^s}+Q_s(Z)\big)(\frac{1}{g}-\frac{1}{g^s})\big)=H_s$ and thus $\frac{1}{g}*H_s=(H_s*\frac{1}{g^s}+Q_s(Z)(\frac{1}{g}-\frac{1}{g^s})$. Let $f\in \FC[Y]$. A straightforward computation yields the formula $\frac{f}{g}*H_s=H_s*(\frac{f}{g})^s+Q_s(Z)(\frac{f}{g}-(\frac{f}{g})^s)$ which finishes the proof of  Proposition~\ref{propDefinition_ATF}.

\begin{remark}
\begin{itemize}
\item Inspired by the proof of \cite[Theorem 6.2]{bardy2016iwahori} we  could try to define $*$ on $V$ as follows. Let $\theta_1,\theta_2\in \FC[Y]$ and $w_1,w_2\in W^v$. Write $\theta_1*H_{w_2}=\sum_{w\in W^v} H_w \theta_w$, with $(\theta_w)\in \FC(Y)^{(W^v)}$. Then $(H_{w_1} *\theta_1)*(H_{w_2}*\theta_2)=\sum_{w\in W} (H_{w_1}*H_w)*(\theta_2\theta_w)$. However it is not clear a priori that the so defined law is associative.

\item Suppose that $\HC_\FC$ is the Iwahori-Hecke algebra associated with some masure defined in \cite[Definition 2.5]{bardy2016iwahori}. Using the same procedure as above (by taking $S=\{Y^\lambda|\lambda\in Y^+\}$), we can construct the algebra $\AC_\FC$ from the algebra $\HC_\FC$. In this particular case, this gives an alternative proof of \cite[Theorem 6.2]{bardy2016iwahori}.
\end{itemize}
\end{remark}

\section{Weight decompositions and intertwining operators}\label{secCY module and intertwining}

Let $\tau \in T_\FC$. In this section, we study the structure of $I_\tau$ as a $\FC[Y]$-module and the set $\Hom_{\AC_\FC-\mathrm{mod}}(I_\tau,I_{\tau'})$ for $\tau'\in T_\FC$. 

In Subsection~\ref{subGeneralized_weight_spaces}, we study the weights of $I_\tau$ and decompose every $\AC_\FC$-submodule of $I_\tau$ as a sum of generalized weight spaces (see Lemma~\ref{lemDecomposition_submodules_I_tau}).

In Subsection~\ref{subIntertwining_operators_weight_space}, we relate intertwining operators and weight spaces. We then prove the existence of nontrivial intertwining operators $I_\tau\rightarrow I_{w.\tau}$ for all $w\in W^v$.

In Subsection~\ref{subNonTrivial_submodules_infinite_dimensional}, we prove that when $W^v$ is infinite, then every nontrivial submodule of $I_\tau$ is infinite dimensional. We deduce that contrary to the reductive case, there exist irreducible representations of $\AC_\FC$ which does not embed in any $I_\tau$.

\subsection{Generalized weight spaces of  \texorpdfstring{$I_\tau$}{Itau}}\label{subGeneralized_weight_spaces}

Let $\tau\in T_\FC$. Let $x\in I_\tau$. Write $x=\sum_{w\in W^v} x_w H_w\otimes_\tau 1$, with $(x_w)\in \FC^{(W^v)}$. Set $\supp(x)=\{w\in W^v|\ x_w\neq 0\}$.  Equip $W^v$ with the Bruhat order. If $E$ is a finite subset of $W^v$, $\max (E)$ is the set of elements of $E$ that are maximal for the Bruhat order.  Let $R$ be a binary relation on $W^v$  (for example $R=$``$\leq$'', $R=$``$\ngeq$'', ...) and $w\in W^v$. 
One sets \[I_\tau^{Rw}=\bigoplus_{v\in W^v|  v Rw} \FC H_v\otimes_\tau 1, \HFW^{Rw}=\bigoplus_{vRw} \FC H_v,\ \ATF^{Rw}=\bigoplus_{v Rw} H_v \FC(Y)\]\index{$I_\tau^{\leq w}, I_\tau^{<w}, \HFW^{\leq w},\ldots$} and $\AF^{Rw}=\ATF^{Rw}\cap \AF=\bigoplus_{vR w} H_v \FC[Y].$

\medskip

Let $V$ be a vector space over $\FC$ and $E\subset \mathrm{End}(V)$. For $\tau\in \FC^E$ set  $V(\tau)=\{v\in V|e.v=\tau(e).v\forall e\in E\}$ and  $V(\tau,\mathrm{gen})=\{v\in V|\exists k\in \N | (e-\tau(e)\Id)^k.v=0, \forall e\in E\}$. Let $\mathrm{Wt}(E)=\{\tau\in \FC^E|V(\tau)\neq \{0\}\}$. 

The following lemma is well known.

\begin{lemma}\label{lemFrobenius_thm}
Let $V$ be a finite dimensional vector space over $\FC$. Let $E\subset \mathrm{End}(V)$ be a subset such that for all $e,e'\in E$, \begin{enumerate}
\item $e$ is triangularizable

\item $ee'=e'e$.
\end{enumerate}

Then $V=\bigoplus_{\tau\in \Wt(E)}V(\tau,\mathrm{gen})$ and in particular $\Wt(E)\neq \emptyset$.

\end{lemma}

For $\tau\in T_\FC$, set $W_\tau=\{w\in W^v|\ w.\tau =\tau\}$.

Let $M$ be a $\AC_\FC$-module. For $\tau\in T_\FC$, set \[M(\tau)=\{m\in M|P.m=\tau(P).m\ \forall P\in \FC[Y]\}\]\index{$M(\tau)$, $I_\tau(\tau)$} and  \[M(\tau,\text{gen})=\{m\in M|\exists k\in \N|\forall P\in \FC[Y], (P-\tau(P))^k.m=0\}\supset M(\tau).\]\index{$M(\tau,\mathrm{gen})$, $I_\tau(\tau,\mathrm{gen})$} Let $\Wt(M)=\{\tau\in T_\FC| M(\tau)\neq\{0\}\}$\index{$\mathrm{Wt}(M)$} and $\Wt(M,\text{gen})=\{\tau\in T_\FC|M(\tau,\text{gen})\neq \{0\}\}$.

\begin{lemma}\label{lemDecomposition_submodules_I_tau}
\begin{enumerate}
\item\label{itMaximal_support} Let $\tau,\tau'\in T_\FC$. Let $x\in I_{\tau}(\tau',\text{gen})$. Then if $x\neq 0$, \[\max\supp(x)\subset \{w\in W^v|\ w.\tau=\tau'\}.\] In particular, if $I_{\tau}(\tau',\text{gen})\neq \{0\}$, then $\tau'\in W^v.\tau$ and thus \[\Wt(I_\tau) \subset W^v.\tau.\]

\item\label{itDecomposition_generalized_weight_spaces}  Let $\tau\in T_\FC$. Let $M\subset I_\tau$ be a $\FC[Y]$-submodule of $I_\tau$.   Then  $\Wt(M)=\Wt(M,\mathrm{gen})\subset W^v.\tau$ and $M=\bigoplus_{\chi\in \Wt(M)} M(\chi,\mathrm{gen})$. In particular, $\Wt(M)\neq \emptyset$.

\end{enumerate}

\end{lemma}

\begin{proof}
(\ref{itMaximal_support}) Let $x\in I_{\tau}(\tau',\text{gen})\setminus\{0\}$. Let $w\in \max \supp(x)$. Write $x=a_wH_w\otimes_{\tau}1+y$, where $a_w\in \FC\setminus\{0\}$ and $y\in I_\tau^{\ngeq w}$. Then by Lemma~\ref{lemCommutation relation}, \[Z^\lambda.x=a_wH_wZ^{w^{-1}.\lambda}\otimes_{\tau}1+y'=\tau(w^{-1}.\lambda)a_wH_w\otimes_{\tau} 1+y'=\tau'(\lambda)a_wH_w\otimes_{\tau} 1+\tau'(\lambda)y,\] where $y'\in I_\tau^{\ngeq w}$.  Therefore $w.\tau=\tau'$.

(\ref{itDecomposition_generalized_weight_spaces}) Let $w\in W^v$. Let  $P\in \FC[Y]$ and $m_P:I_\tau^{\leq w}\rightarrow I_\tau^{\leq w}$ be defined by $m_P(x)=P.x$ for all $x\in I_\tau^{\leq w}$. Then by Lemma~\ref{lemCommutation relation}, $(m_P-w.\tau(P)\Id)(I_\tau^{\leq w})\subset I_\tau^{<w}$. By induction on $\ell(w)$ we deduce that $m_P$ is triangularizable on $I_\tau^{\leq w}$ and $\Wt(I_\tau^{\leq w})\subset [1,w].\tau\subset W^v.\tau$. 

Let $x\in M$ and $M_x=\FC[Y].x$. By the fact that $(W^v,\leq)$ is  a directed poset and by  Lemma~\ref{lemCommutation relation}, there exists $w\in W^v$ such that $M_x\subset I_\tau^{\leq w}$.  Therefore, for all $P\in \FC[Y]$, $m_P:M_x\rightarrow M_x$ is triangularizable. Thus by Lemma~\ref{lemFrobenius_thm}, \[\FC[Y].x=\bigoplus_{\chi\in \Wt(M_x,\text{gen})}M_x(\chi,\text{gen})=\bigoplus_{\chi\in W^v.\tau}M_x(\chi,\text{gen}).\] Consequently, $M=\sum_{x\in M} M_x=\bigoplus_{\chi\in \Wt(M,\text{gen})}M(\chi,\text{gen})$ and $\Wt(M)\subset \bigcup_{w\in W^v}\Wt(I_\tau^{\leq w})\subset W^v.\tau$.

Let $\chi\in \Wt(M,\text{gen})$. Let $x\in M(\chi,\text{gen})\setminus \{0\}$ and $N=\FC[Y].x$. Then by Lemma~\ref{lemCommutation relation}, $N$ is a finite dimensional submodule of $I_\tau$. By Lemma~\ref{lemFrobenius_thm}, $\Wt(N)\neq \emptyset$. As $\Wt(N)\subset \{\chi\}$, $\chi\in \Wt(M)$. Thus $\Wt(M,\text{gen})\subset \Wt(M)$ and as the other inclusion is clear, we get the lemma. \end{proof}

\begin{proposition}\label{propNecessary condition for existence of intertwinners}(see \cite[4.3.3 Théorème (iii)]{matsumoto77Analyse})
Let $\tau,\tau'\in T_\FC$ and $M$ (resp. $M'$) be a $\AC_\FC$-submodule of $I_\tau $ (resp. $I_{\tau'}$).  Assume that $ \Hom_{\AC_\FC-\mathrm{mod}}(M,M')\setminus\{0\}$. Then  $\tau'\in W^v.\tau$.

\end{proposition}

\begin{proof}
Let $f\in \Hom_{\AC_\FC}(M,M')\setminus\{0\}$.  Then by Lemma~\ref{lemDecomposition_submodules_I_tau}~(\ref{itDecomposition_generalized_weight_spaces}), there exists $w\in W^v/W_\tau$ such that $f\big(M(w.\tau,\mathrm{gen})\big)\neq \{0\}$. Then $w.\tau\in \mathrm{Wt}(I_{\tau'})$ and by Lemma~\ref{lemDecomposition_submodules_I_tau}~(\ref{itMaximal_support}) the proposition follows.
\end{proof}

An element $\tau\in T_\FC$ is said to be \textbf{regular} if $w.\tau\neq \tau$ for all $w\in W^v\setminus\{1\}$. We denote by $T^{\mathrm{reg}}_\FC$\index{$T^{\mathrm{reg}}_\FC$} the set of regular elements of $T_\FC$. 

\begin{proposition}\label{propDecomposition_generalized_weight_spaces_Itau}(see \cite[Proposition 1.17]{kato1982irreducibility})
Let $\tau\in T_\FC$. \begin{enumerate}
\item\label{itGeneralized_weight_vectors} There exists a basis $(\xi_w)_{w\in W^v}$ of $I_\tau$ such that for all $w\in W^v$:\begin{itemize}
\item $\xi_w\in I_\tau^{\leq w}$ and $\pi^H_w(\xi_w)=1$

\item $\xi_w\in I_\tau(w.\tau,\mathrm{gen})$.
\end{itemize} Moreover, if $w\in W^v$ is minimal for $\leq$  among $\{v\in W^v|v.\tau=w.\tau\}$, then $\xi_w\in I_\tau(w.\tau)$. In particular, $\Wt(I_\tau)=W^v.\tau$. 

\item\label{itDecomposition_Itau_regular_case} If $\tau$ is regular, then $I_\tau(w.\tau,\mathrm{gen})=I_\tau(w.\tau)$ is one dimensional for all $w\in W^v$ and $I_\tau=\bigoplus_{w\in W^v}I_\tau(w.\tau)$. 

\end{enumerate}

\end{proposition}

\begin{proof}
(\ref{itGeneralized_weight_vectors}) Let $w\in W^v$. Then by Lemma~\ref{lemCommutation relation}, Lemma~\ref{lemFrobenius_thm} and Lemma~\ref{lemDecomposition_submodules_I_tau}, \[I_\tau^{\leq w}=\bigoplus_{\overline{v}\in W^v/W_\tau}I_\tau^{\leq w}(\overline{v}.\tau,\mathrm{gen}).\] Write $H_w\otimes_\tau 1= \sum_{\overline{v}\in W^v/W_\tau} x_{\overline{v}}$, where $x_{\overline{v}}\in I_\tau^{\leq w}(v.\tau,\mathrm{gen})$ for all $\overline{v}\in W^v/W_\tau$.  Let $\overline{v}\in W^v/W_\tau$ be such that $\pi^H_w(x_{\overline{v}})\neq 0$. Then $\max \supp (x_{\overline{v}})=\{w\}$ and by Lemma~\ref{lemDecomposition_submodules_I_tau}, $w.\tau=\overline{v}.\tau$. Set $\xi_w=\frac{1}{\pi^H_w(x_{\overline{v}})} x_{\overline{v}}$. Then $(\xi_u)_{u\in W^v}$ is a basis of $I_\tau$ and has the desired properties. Let  $w\in W^v$ be minimal for $\leq$  among $\{v\in W^v|v.\tau=w.\tau\}$. Let $\lambda\in Y$.  Then by Lemma~\ref{lemCommutation relation}, $(Z^\lambda-w.\tau(\lambda).\xi_w)\in I_\tau(w.\tau,\mathrm{gen})\cap I_\tau^{<w}$. By Lemma~\ref{lemDecomposition_submodules_I_tau}, we deduce that $(Z^\lambda-w.\tau(\lambda)).\xi_w=0$ and thus that $\xi_w\in I_\tau(w.\tau)$. Thus $w.\tau\in \Wt(I_\tau)$ and by Lemma~\ref{lemDecomposition_submodules_I_tau}, $\Wt(I_\tau)=I_\tau$. 

(\ref{itDecomposition_Itau_regular_case}) Suppose that $\tau$ is regular. Let $w\in W^v$, $\lambda\in Y$ and $x\in I_\tau(\tau,\mathrm{gen})$.  Then by Lemma~\ref{lemDecomposition_submodules_I_tau}~(\ref{itMaximal_support}), $x-\pi^H_w(x)\xi_w\in I_\tau(\tau,\mathrm{gen})\cap I_\tau^{<w}=\{0\}$. By~(\ref{itGeneralized_weight_vectors}), $\xi_w\in I_\tau(w.\tau)$ and thus $I_\tau(\tau)=I_\tau(\tau,\mathrm{gen})$ is one dimensional.  By Lemma~\ref{lemDecomposition_submodules_I_tau}, we deduce that $I_\tau=\bigoplus_{w\in W^v}I_\tau(w.\tau)$. 
\end{proof}

\subsection{Intertwining operators and weight spaces}\label{subIntertwining_operators_weight_space}

In this subsection, we relate  intertwining operators and weight spaces and study some consequences. Let $\tau\in T_\FC$. Using Subsection~\ref{subGeneralized_weight_spaces}, we prove the existence of nonzero morphisms $I_\tau\rightarrow I_{w.\tau}$ for all $w\in W^v$. We will give a more precise construction of such morphisms in Subsection~\ref{subWeight_vectors_regarded_as_rational_functions}.

Let $M$ be a $\AC_\FC$-module and $\tau\in T_\FC$. For $x\in M(\tau)$ define $\Upsilon_x:I_\tau\rightarrow M$\index{$\Upsilon$} by $\Upsilon_x(u.1\otimes_\tau 1)=u.x$, for all $u\in \AC_\FC$. Then $\Upsilon_x$ is well-defined. Indeed, let $u\in \AC_\FC$ be such that $u.1\otimes_\tau 1=0$. Then $u\in \FC[Y]$ and $\tau(u)=0$. Therefore $u.x=0$ and hence $\Upsilon_x$ is well-defined. The following lemma is then  easy to prove.

\begin{lemma}\label{lemFrobenius_reciprocity} (Frobenius reciprocity, see \cite[Proposition 1.10]{kato1982irreducibility})
Let $M$ be a $\AC_\FC$-module, $\tau\in T_\FC$ and $x\in M(\tau)$. Then the map $\Upsilon:M(\tau)\rightarrow \Hom_{\AC_\FC-mod}(I_\tau,M)$  mapping each $x\in M(\tau)$ to $\Upsilon_x$ is a vector space isomorphism and $\Upsilon^{-1}(f)=f(1\otimes_\tau 1)$ for all $f\in  \Hom_{\AC_\FC-mod}(I_\tau,M)$. 
\end{lemma}

\begin{proposition}\label{propLink Matsumoto def with BL def}(see \cite[(4.1.10)]{matsumoto77Analyse})
Let $M$ be a  $\AC_\FC$-module such that there exists $\xi\in M$ satisfying: \begin{enumerate}
\item there exists $\tau\in T_\FC$ such that $\xi\in M(\tau)$,

\item  $M=\AC_\FC.\xi$.

\end{enumerate}

Then there exists a surjective morphism $\phi:\I_\tau\twoheadrightarrow M$ of $\AC_\FC$-modules.
\end{proposition}

\begin{proof}
One can take $\phi=\Upsilon_\xi$, where $\Upsilon$ is as in Lemma~\ref{lemFrobenius_reciprocity}. \end{proof}

\begin{proposition} (see \cite[Th{\'e}or{\`e}me 4.2.4]{matsumoto77Analyse})
 Let $M$ be an irreducible representation of $\AC_{\FC}$ containing a finite
 dimensional $\FC[Y]$-submodule $M'\neq\{0\}$. Then there exists $\tau\in T_{\FC}$ such that there exists a surjective morphism of $\AC_{\FC}$-modules $\phi:I_\tau\twoheadrightarrow M$.
\end{proposition}

\begin{proof}
By Lemma~\ref{lemFrobenius_thm}, there exists $\xi\in M'\setminus\{0\}$ such that $Z^\mu.\xi\in \FC .\xi$ for all $\mu\in Y$. Let $\tau\in T_{\FC}$ be such that $\xi\in M(\tau)$. Then we conclude with Proposition~\ref{propLink Matsumoto def with BL def}.\end{proof}

\begin{remark}
Let $\mathcal{Z}(\AF)$ be the center of $\AF$. When $W^v$ is finite, it is well known that $\AC_\FC$ is a finitely generated $\mathcal{Z}(\AC_\FC)$ module and thus every irreducible representation of $\AC_\FC$ is finite dimensional. Assume that $W^v$ is infinite.  Using the same reasoning as in \cite[Remark 4.32]{abdellatif2019completed} we can prove that $\AC_\FC$ is not a finitely generated $\mathcal{Z}(\AF)$-module. As we shall see (see Remark~\ref{rkMeasure_irreducible_representations}), when $\FC=\C$, there exist irreducible  infinite dimensional representations of $\AF$. However we do not know if there exist an irreducible representation $V$ of $\AF$ such that for all $x\in V\setminus \{0\}$, $\FC[Y].x$ is infinite dimensional or equivalently, a representation which is not a quotient of a principal series representation.
\end{remark}

\begin{proposition}\label{propExistence_nontrivial_morphisms} (see \cite[(1.21)]{kato1982irreducibility})
Let $\tau\in T_\FC$ and $w\in W^v$. Then \[\Hom_{\AF-\mathrm{mod}}(I_\tau,I_{w.\tau})\neq \{0\}.\]
\end{proposition}

\begin{proof}
By Proposition~\ref{propDecomposition_generalized_weight_spaces_Itau} $w.\tau\in \Wt(I_\tau)$ and we conclude with Lemma~\ref{lemFrobenius_reciprocity}.
\end{proof}

\subsection{Nontrivial submodules of $I_\tau$ are infinite dimensional}\label{subNonTrivial_submodules_infinite_dimensional}

In this subsection, we prove that when $W^v$ is infinite, then every submodule of $I_\tau$ is infinite dimensional. We then deduce that there can exist an irreducible representation of $\AC_\C$ such that $V$ does not embed in any $I_\tau$, for $\tau\in T_\C$.

\begin{lemma}\label{lemNon_existence_left_maximal_elements}
Assume that $W^v$ is infinite. Let $w\in W^v$. Then there exists $s\in \SCC$ such that $sw>w$. 
\end{lemma}

\begin{proof}
Let $D_L(w)=\{s\in \SCC|sw<w\}$. By the proof of \cite[Lemma~3.2.3]{bjorner2005combinatorics}, $\SCC\nsubseteq D_L(w)$, which proves the lemma.
\end{proof}

\begin{proposition}\label{propSubmodules_infinite_dimensional} (compare \cite[4.2.4]{matsumoto77Analyse})
Let $\tau\in T_\FC$. Let $M\subset I_\tau$ be a nonzero $\HC_{W^v,\FC}$-submodule. Then the dimension of $M$ is infinite. In particular, if $V$ is a finite dimensional irreducible representation of $\AC_\FC$, then $\Hom_{\AC_\FC-\mathrm{mod}}(V,I_\tau)=\{0\}$ for all $\tau\in T_\FC$.
\end{proposition}

\begin{proof}
Let $m\in M\setminus\{0\}$. Let $\ell(m)=\max\{ \ell(v)|v\in \supp(m)\}$.  Let $w\in \supp(m)$ be such that $\ell(w)=\ell(m)$. By Lemma~\ref{lemNon_existence_left_maximal_elements} there exists $(s_n)\in \SCC^\Ne$ such that if $w_1=w$ and $w_{n+1}=s_n w_{n}$ for all $n\in \Ne$, one has $\ell(w_{n+1})=\ell(w_n)+1$ for all $n\in \Ne$. Let $m_1=m$ and $m_{n+1}=H_{s_n}.m_n$ for all $n\in \Ne$. Then for all $n\in \Ne$, $w_n\in\max\big(\supp(m_n)\big)$, which proves that $M$ is infinite dimensional.
\end{proof}

As we shall see in Appendix~\ref{appExistence_one_dimensional_representations}, there can exist finite dimensional representations of $\AC_\C$.

\section{Study of the irreducibility of $I_\tau$}\label{secStudy of irreducibility}

In this section, we study the irreducibility of $I_\tau$. 

In Subsection~\ref{subintertwining_operators_simple_reflections}, we describe certain intertwining operators between $I_\tau$ and $I_{s.\tau}$, for $s\in \SCC$ and $\tau\in T_\FC$. For this, we   introduce elements $F_s\in \ATF$ such that $F_s(\chi)\otimes_\chi 1 \in I_\chi(s.\chi)$ for all $\chi\in \T_\FC$ for which this is well-defined.

In Subsection~\ref{subNecessary condition for irreducibility}, we establish that the condition~(\ref{itCondition on the values of chi}) appearing in Theorems~\ref{thm*Kato's theorem}, \ref{thm*Matsumoto's criterion} and~\ref{thm*Kato's criterion} is  a necessary condition for the irreducibility of $I_\tau$. This conditions comes from the fact that when $I_\tau$ is irreducible, certain intertwinners have to be isomorphisms. 

In Subsection~\ref{subIrreducibility criterion}, we prove an irreducibility criterion for $I_\tau$ involving the dimension of $I_\tau$ and the values of $\tau$ (see Theorem~\ref{thmIrreducibility criterion}). We then deduce Matsumoto criterion.

In Subsection~\ref{subWeight_vectors_regarded_as_rational_functions} we introduce and study, for every $w\in W^v$, an element $F_w\in \ATF$ such that $F_w(\chi)\otimes_\chi 1 \in I_\chi(w.\chi)$ for every $\chi\in T_\C$ for which this is well-defined.

In Subsection~\ref{subOne implication of Kato's theorem} we prove one implication of Kato's criterion (see Proposition~\ref{propKato's weak theorem}).

The definition we gave for $I_\tau$ is different from the definition of Matsumoto (see~\cite[(4.1.5)]{matsumoto77Analyse}). It seems to be well known that these definitions are equivalent. We justify this equivalence in Subsection~\ref{subLink with Matsumoto and Kato}. We also explain why it seems difficult to adapt Kato's proof in our framework.

\subsection{Intertwining operators associated with simple reflections}\label{subintertwining_operators_simple_reflections}

 Let $s\in \SCC$. In this subsection we define and study an element $F_s\in \ATF$ such that $F_s(\chi)\otimes_\chi 1\in I_\chi(s.\chi)$ for all $\chi$ such that $F_s(\chi)$ is well-defined. 

Let $s\in \SCC$ and $T_s=\sigma_s H_s$. Let $w\in W^v$ and $w=s_1\ldots s_k$ be a reduced writing. Set $T_w=T_{s_1}\ldots T_{s_k}$\index{$T_w$}. This is independent of the choice of the reduced writing by \cite[6.5.2]{bardy2016iwahori}.

 Set $B_s=T_s-\sigma_s^2\in \HC_{W^v,\FC}$\index{$B_s$}. One has $B_s^2=-(1+\sigma_s^2)B_s$. Let  $\zeta_{s}=-\sigma_sQ_s(Z)+\sigma_s^2\in \FC(Y)\subset \AC(T_\FC)$\index{$\zeta_s$, $\zeta_{\alpha^\vee}$}.  When $\sigma_s=\sigma_s'=\sqrt{q}$ for all $s\in \SCC$, we have $\zeta_{s}=\frac{1-qZ^{-\alpha_s^\vee}}{1-Z^{-\alpha^\vee_s}}\in \FC(Y)$. Let $F_s=B_s+\zeta_{s}\in \AC(T_\FC)$\index{$F_s$}.

Let  $\alpha^\vee\in \Phi^\vee$. Write $\alpha^\vee=w.\alpha_s^\vee$ for $w\in W^v$ and $s\in \SCC$. We set $\zeta_{\alpha^\vee}=(\zeta_s)^w$.  

Let $\alpha^\vee\in\Phi^\vee$. Write $\alpha=w.\alpha_s^\vee$, with $w\in W^v$ and $s\in \SCC$. We set $\sigma_{\alpha^\vee}=\sigma_s$\index{$\sigma_{\alpha^\vee}$, $\sigma'_{\alpha^\vee}$} and $\sigma_{\alpha^\vee}'=w.\sigma_s'$. This is well-defined by Lemma~\ref{lemKumar_1.3.14} and by the relations on the $\sigma_t$, $t\in \SCC$ (see Subsection~\ref{subIH algebras}).

The ring $\FC[Y]$ is a  unique factorization domain. For $\alpha^\vee$, write $\zeta_{\alpha^\vee}=\frac{\zeta_{\alpha^\vee}^{\mathrm{num}}}{\zeta_{\alpha^\vee}^{\mathrm{den}}}$\index{$\zeta_{\alpha^\vee}^{\mathrm{den}}$, $\zeta_{\alpha^\vee}^{\mathrm{num}}$} where $\zeta_{\alpha^\vee}^{\mathrm{num}},\zeta_{\alpha^\vee}^{\mathrm{den}}\in \FC[Y]$ are pairwise coprime. For example if $\alpha^\vee\in \Phi^\vee$ is such that $\sigma_{\alpha^\vee}=\sigma_{\alpha^\vee}'$  we can take $\zeta_{\alpha^\vee}^{\mathrm{den}}=1-Z^{-\alpha^\vee}$ and in any case we will choose $\zeta_{\alpha^\vee}^{\mathrm{den}}$ among $\{1-Z^{-\alpha^\vee},1+Z^{-\alpha^\vee},1-Z^{-2\alpha^\vee}\}$.

\begin{remark}\label{rkNecessary_condition_stau_neq_tau}
Let $\tau\in T_\FC$ and $r=r_{\alpha^\vee}\in \RCC$. Suppose that $r.\tau\neq \tau$. Then $\zeta_{\alpha^\vee}^{\mathrm{den}}(\tau)\neq 0$. Indeed, let $\lambda\in Y$ be such that $\tau(r.\lambda)\neq \tau(\lambda)$. Then $\tau(r.\lambda-\lambda)=\tau(\alpha_r^\vee)^{\alpha_r(\lambda)}\neq 1$. Suppose  $\sigma_{\alpha^\vee}=\sigma_{\alpha^\vee}'$, then $\zeta_{\alpha^\vee}^{\mathrm{den}}=1-Z^{-\alpha_r^\vee}$ and thus $\tau(\zeta_{\alpha^\vee}^{\mathrm{den}})\neq0$. Suppose $\sigma_r=\sigma_r'$. Then $\alpha_r(\lambda)\in 2\Z$ thus $\tau(\alpha_r^\vee)\notin \{-1,1\}$ and hence $\tau(\zeta_{\alpha^\vee}^{\mathrm{den}})\neq 0$. 
\end{remark}

\begin{lemma}\label{lem1.10 of Reeder}
Let $s\in \SCC$ and $\theta\in \FC(Y)$. Then \[\theta*F_s=F_s*\theta^s.\] In particular, for all $\tau\in T_\FC$ such that $\tau(\zeta_{s}^\mathrm{den})\neq 0$, $F_s(\tau)\otimes_\tau 1\in I_\tau(s.\tau)$ and $F_s(\tau)\otimes_{s.\tau} 1\in I_{s.\tau}(\tau)$. 
\end{lemma}

\begin{proof}
Let $\lambda\in Y$. Then \[\begin{aligned} Z^\lambda*B_s-B_s*Z^{s.\lambda} = & \sigma_s(Z^\lambda*H_s-H_s*Z^{s.\lambda})+\sigma_s^2(Z^{s.\lambda}-Z^\lambda) \\  =  & -\sigma_s Q_s(Z)(Z^{s\lambda}-Z^{\lambda})+\sigma_s^2(Z^{s.\lambda}-Z^\lambda)\\  =  &\zeta_s (Z^{s.\lambda}-Z^\lambda).\end{aligned}\]

 Thus $Z^\lambda*F_s=Z^{\lambda}*(B_s+\zeta_s)=F_s*Z^{s.\lambda}$ and hence $\theta*F_s=F_s*\theta^s$ for all $\theta\in \FC[Y]$.

Let $\theta\in \FC[Y]\setminus \{0\}$. Then $\theta*(F_s*\frac{1}{\theta^s})=F_s$ and thus $\frac{1}{\theta}*F_s=F_s*\frac{1}{\theta^s}$.  Lemma follows.
\end{proof}

\begin{lemma}\label{lemFs2}
Let $s\in \SCC$. Then $F_s^2=\zeta_s\zeta_s^s\in \FC(Y)\subset \AC(T_\FC)$.

\end{lemma}

\begin{proof}
By Lemma~\ref{lem1.10 of Reeder}, one has: \[\begin{aligned} F_s^2&=(B_s+\zeta_s)*F_s\\ &=B_s*F_s+F_s*\zeta_s^s\\ &=B_s^2+B_s\zeta_s+B_s\zeta_s^s+\zeta_s\zeta_s^s \\ &=B_s(-1-\sigma_s^2+\zeta_s+\zeta_s^s)+\zeta_s\zeta_s^s\\ &=\zeta_s\zeta_s^s.\end{aligned}\]
\end{proof}

\subsection{A necessary condition for irreducibility}\label{subNecessary condition for irreducibility}

In this subsection, we establish that the condition~(\ref{itCondition on the values of chi}) appearing in Theorems~\ref{thm*Kato's theorem}, \ref{thm*Matsumoto's criterion} and~\ref{thm*Kato's criterion} is  a necessary condition for the irreducibility of $I_\tau$.
 
Recall the definition of $\Upsilon$ from Subsection~\ref{subIntertwining_operators_weight_space}. 
 
\begin{lemma}\label{lemCondition on values for isomorphisms}
Let $\tau\in T_\FC$ and $s\in \SCC$ be such that $\tau(\zeta_s^{\mathrm{den}})\tau((\zeta_s^{\mathrm{den}})^s)\neq 0$.  Let $\phi(\tau,s.\tau)=\Upsilon_{F_s(\tau)\otimes_{s.\tau} 1}:I_\tau\rightarrow I_{s.\tau}$ and $\phi(s.\tau,\tau)=\Upsilon_{F_s(\tau)\otimes_\tau 1}:I_{s.\tau}\rightarrow I_\tau$. Then \[\phi(s.\tau,\tau)\circ\phi(\tau,s.\tau)=\tau(\zeta_s\zeta^s_s)\Id_{I_\tau} \text{ and }\phi(\tau,s.\tau)\circ\phi(s.\tau,\tau)=\tau(\zeta_s\zeta^s_s)\Id_{I_{s.\tau}}.\]
\end{lemma}

\begin{proof}
By Lemma~\ref{lem1.10 of Reeder} and  Lemma~\ref{lemFrobenius_reciprocity}, $\phi(s.\tau,\tau)$ and $\phi(\tau,s.\tau)$ are well-defined. Let $f=\phi(s.\tau,\tau)\circ\phi(\tau,s.\tau)\in \mathrm{End}_{\AC_\FC-\mathrm{mod}}(I_\tau)$. Then by Lemma~\ref{lem1.10 of Reeder} and  Lemma~\ref{lemFs2}: \[f(1\otimes_\tau 1)=\phi(s.\tau,\tau)\big(F_s(\tau)\otimes_{s.\tau} 1\big)=F_s(\tau).\phi(s.\tau,\tau)\big(1\otimes_{s.\tau} 1\big)= F_s(\tau)^2\otimes_\tau 1=\tau(\zeta_s\zeta^s_s)\otimes_\tau 1.\] By symmetry, we get the lemma.
\end{proof}

Let $\UC_\FC$\index{$\UC_\FC$} be the set of $\tau\in T_\FC$ such that for all $\alpha^\vee\in \Phi^\vee$,  $\tau(\zeta_{\alpha^\vee}^{\mathrm{num}})\neq 0$. When $\sigma_s=\sigma_s'=\sqrt{q}$ for all $s\in \SCC$, then $\UC_\FC=\{\tau\in T_\FC|\tau(\alpha^\vee)\neq q,\ \forall \alpha^\vee\in \Phi^\vee\}$. 

We assume that for all $s\in \SCC$, $\sigma_s'\notin \{\sigma_s^{-1}, -\sigma_s,-\sigma_s^{-1}\}$. Under this condition, if $\alpha^\vee\in \Phi^\vee$ and $\tau\in T_\FC$ are such that $\tau(\zeta_{\alpha^\vee}^{\mathrm{den}})=0$, then $\tau(\zeta_{\alpha^\vee}^{\mathrm{num}})\neq 0$.

\begin{lemma}\label{lemIrreducibility implies isomorphisms}
\begin{enumerate}
\item\label{itU implies isomorphism} Let $\tau\in \UC_\FC$. Then for all $w\in W^v$, $I_\tau$ and $I_{w.\tau}$ are isomorphic as $\AC_\FC$-modules.

\item\label{itIrreducible implies U} Let $\tau\in T_\FC$ be such that $I_\tau$ is irreducible. Then $\tau\in \UC_\FC$.

\end{enumerate}
\end{lemma}

\begin{proof}
Let $\tau\in \UC_\FC$. Let $w\in W^v$ and $\tilde{\tau}=w.\tau$. Let $s\in \SCC$. Assume that $s.\tilde{\tau}\neq \tilde{\tau}$. Then by Remark~\ref{rkNecessary_condition_stau_neq_tau}, $\zeta_s^{\mathrm{den}}(\tilde{\tau})\neq 0$ and $\zeta_s^{\mathrm{den}}(s.\tilde{\tau})\neq 0 $. Therefore $\zeta_s(\tau)$, $\zeta_s(s.\tilde{\tau})$ are well-defined and hence $F_s(\tilde{\tau})$, $F_s(\tilde{\tau})$ are well-defined. Let $\phi(\tilde{\tau},s.\tilde{\tau})=\Upsilon_{F_s(\tilde{\tau})\otimes_{s.\tilde{\tau}} 1}:I_{\tilde{\tau}}\rightarrow I_{s.\tilde{\tau}}$ and $\phi(s.\tilde{\tau},\tilde{\tau})=\Upsilon_{F_s(\tilde{\tau})\otimes_{\tilde{\tau}} 1}:I_{s.\tilde{\tau}}\rightarrow I_{\tilde{\tau}}$. Then by Lemma~\ref{lemCondition on values for isomorphisms}, \[\phi(s.\tilde{\tau},\tilde{\tau})\circ\phi(\tilde{\tau},s.\tilde{\tau})=\tilde{\tau}(\zeta_s\zeta^s_s)\Id_{I_{\tilde{\tau}}} \text{ and }\phi(\tilde{\tau},s.\tilde{\tau})\circ\phi(s.\tilde{\tau},\tilde{\tau})=\tilde{\tau}(\zeta_s\zeta^s_s)\Id_{I_{s.\tilde{\tau}}}.\]  By definition of $\UC_\FC$, $\tilde{\tau}(\zeta_s\zeta^s_s)=\tilde{\tau}(\zeta_s)\tilde{\tau}(\zeta^s_s)\neq 0$ and thus $\phi(\tilde{\tau},s.\tilde{\tau})$ and $\phi(s.\tilde{\tau},\tilde{\tau})$ are isomorphisms. Consequently $I_{\tilde{\tau}}$ is isomorphic to $I_{s.\tilde{\tau}}$ and~(\ref{itU implies isomorphism}) follows by induction.

Let $\tau\in T_\FC$ be such that $I_\tau$ is irreducible. Let $s\in \SCC$. 

Suppose $\tau(\zeta_s^{\mathrm{den}})=0$. Then by assumption, $\tau(\zeta_s^{\mathrm{num}})\neq 0$. Moreover by Remark~\ref{rkNecessary_condition_stau_neq_tau}, $I_{s.\tau}=I_\tau$. 

Suppose now $\tau(\zeta_s^{\mathrm{den}})\neq 0$. Then (with the same notations as in Lemma~\ref{lemCondition on values for isomorphisms}), $\phi(s.\tau,\tau)\neq 0$ and $\mathrm{Im}\big(\phi(s.\tau,\tau)\big)$ is a $\AC_\FC$-submodule of $I_\tau$: $\mathrm{Im}\big(\phi(s.\tau,\tau)\big)=I_\tau$. Therefore $\phi(\tau,s.\tau)\circ \phi(s.\tau,\tau)\neq 0$.  Thus by  Lemma~\ref{lemCondition on values for isomorphisms}, $\phi(\tau,s.\tau)$ is an isomorphism and  $\tau(\zeta_s\zeta_s^s)\neq 0$. In particular, $\tau(\zeta_s^{\mathrm{num}})\neq 0$.

 Therefore in any cases, $I_\tau$ is isomorphic to $I_{s.\tau}$ and $\tau(\zeta_s^{\mathrm{num}})\neq 0$.  By induction we deduce that $I_{w.\tau}$ is isomorphic to $I_\tau$. Thus $I_{w.\tau}$ is irreducible for all $w\in W^v$. Thus  $w.\tau(\zeta_s^{\mathrm{num}})\neq 0$ for all $w\in W^v$ and $s\in \SCC$, which proves that $\tau\in \UC_\FC$.
\end{proof}

\begin{lemma}\label{lemIsomorphism between weight spaces}
Let $\tau\in T_\FC$ be such that $I_{w.\tau}\simeq I_\tau$ (as a $\AC_\FC$-module) for all $w\in W^v$. Then for all $w\in W^v$, there exists a vector space isomorphism $I_\tau(\tau)\simeq I_\tau (w.\tau)$.
\end{lemma}

\begin{proof}
Let $w\in W^v$. Then by hypothesis, $\mathrm{Hom}_{\AC_{\FC}-\mathrm{mod}}(I_\tau,I_\tau)\simeq \mathrm{Hom}_{\AC_\FC-\mathrm{mod}}(I_{w.\tau},I_{w.\tau})$. Let $\phi:I_\tau\rightarrow I_{w.\tau}$ be a $\AC_\FC$-module isomorphism. Then $\phi$ induces an isomorphism of vector spaces $I_\tau(w.\tau)\simeq I_{w.\tau}(w.\tau)$. By Lemma~\ref{lemFrobenius_reciprocity}, 

$I_\tau(\tau)\simeq \mathrm{Hom}_{\AC_\FC-\mathrm{mod}}(I_\tau,I_\tau)\simeq \mathrm{Hom}_{\AC_\FC-\mathrm{mod}}(I_{w.\tau},I_{w.\tau})\simeq I_{w.\tau}(w.\tau)\simeq I_\tau(w.\tau).$\end{proof}

\subsection{An irreducibility criterion for $I_\tau$}\label{subIrreducibility criterion}

In this subsection, we give a characterization of irreducibility for $I_\tau$, for $\tau\in T_\C$.

If $\mathcal{B}$ is a $\C$-algebra with unity $e$ and $a\in \mathcal{B}$, one sets \[\mathrm{Spec}(a)=\{\lambda\in \C|\ a-\lambda e\mathrm{\ is\ not\ invertible}\}.\] Recall the following theorem of Amitsur (see Th{\'e}or{\`e}me B.I of \cite{renard2010representations}): 

\begin{theorem}\label{thmAmitsur}
Let $\BC$ be a $\C$-algebra with   unity $e$. Assume that the dimension of $\BC$ over $\C$ is countable. Then for all $a\in \BC$, $\mathrm{Spec}(a)\neq \emptyset$.
\end{theorem}

Recall that $\UC_\C$ is the set of $\tau\in T_\C$ such that for all $\alpha^\vee\in \Phi^\vee$,  $\tau(\zeta_{\alpha^\vee}^{\mathrm{num}})\neq 0$. 

\begin{theorem}\label{thmIrreducibility criterion}
Let $\tau\in T_\C$. Then the following are equivalent: \begin{enumerate}
\item\label{itIrreducibility} $I_\tau$ is irreducible,

\item\label{itWeight vectors criterion}  $I_\tau(\tau)=\C. 1\otimes_\tau 1$ and $\tau\in \UC_\C$,

\item\label{itEndomorphisms criterion} $\mathrm{End}_{\AC_\C-\mathrm{mod}}(I_\tau)=\C.  \mathrm{Id}$ and $\tau\in \UC_\C$.
\end{enumerate}

\end{theorem}

\begin{proof}
Assume that $\BC=\mathrm{End}_{\AC_\C-\mathrm{mod}}(I_\tau)\neq \C \mathrm{Id}$. By Lemma~\ref{lemFrobenius_reciprocity} and the fact that $I_\tau$ has countable dimension, $\BC$ has countable dimension. Let $\phi\in \BC\setminus \C \mathrm{Id}$. Then by Amitsur Theorem, there exists $\gamma\in \mathrm{Spec}(\phi)$. Then $\phi-\gamma \mathrm{Id}$ is non-injective or non-surjective and therefore $\mathrm{Ker}(\phi-\gamma \mathrm{Id})$ or $\mathrm{Im} (\phi-\gamma \mathrm{Id})$ is a non-trivial $\AC_\C$-module, which proves that $I_\tau$ is reducible. Using Lemma~\ref{lemIrreducibility implies isomorphisms} we deduce that~(\ref{itIrreducibility}) implies~(\ref{itEndomorphisms criterion}). 

By Lemma~\ref{lemFrobenius_reciprocity}, (\ref{itWeight vectors criterion}) is equivalent to~(\ref{itEndomorphisms criterion}).

Let $\tau\in T_\C$ satisfying~(\ref{itWeight vectors criterion}). Then by Lemma~\ref{lemIrreducibility implies isomorphisms} and  Lemma~\ref{lemIsomorphism between weight spaces}, $\dim I_\tau(w.\tau)=1$ for all $w\in W^v$. By Lemma~\ref{lemIrreducibility implies isomorphisms}, for all $w\in W^v$, there exists an  isomorphism of $\AC_\C$-modules $f_w:I_{w.\tau}\rightarrow I_\tau$. As $\C.f_w(1\otimes_{w.\tau}1)\subset I_\tau(w.\tau)$ we deduce that $I_\tau(w.\tau)=\C.f_w(1\otimes_{w.\tau} 1)$ for all $w\in W^v$. 

Let $M\neq \{0\}$ be a $\AC_\C$-submodule of $I_\tau$. Let $x\in M\setminus\{0\}$. Then $M'=\C[Y].x$ is a finite dimensional $\C[Y]$-module. Thus by Lemma~\ref{lemFrobenius_thm}), there exists $\xi\in M'\setminus \{0\}$ such that $Z^\lambda.\xi\in \C.\xi$ for all $\lambda\in Y$. Then $\xi\in I_\tau(\tau')$ for some $\tau'\in T_\C$. By Lemma~\ref{lemDecomposition_submodules_I_tau}, $\tau'=w.\tau$, for some $w\in W^v$. Thus $\xi\in \C^*f_w(1\otimes_{w.\tau}1)$. One has \[\AC_\C.\xi=f_w(\AC_\C.1\otimes_{w.\tau} 1)=f_w(I_{w.\tau})=I_\tau\subset M.\] Hence $I_\tau$ is irreducible, which finishes the proof of the theorem.
\end{proof}

\begin{remark}\label{rkAlgebraically closed field}
Actually, our proof of the equivalence between~(\ref{itWeight vectors criterion}) and~(\ref{itEndomorphisms criterion}), and of the fact that~(\ref{itWeight vectors criterion}) implies~(\ref{itIrreducibility}) is valid when $\FC$ is a field, without assuming   $\FC=\C$.
\end{remark}

Recall that  an element $\tau\in T_{\FC}$ is called regular if $w.\tau\neq \tau$ for all $w\in W^v$.

\begin{corollary}\label{corMatsumoto theorem}(see \cite[Th{\'e}or{\`e}me 4.3.5]{matsumoto77Analyse} )
Let $\tau\in T_{\FC}$ be regular. Then $I_\tau$ is irreducible if and only if $\tau\in \UC_{\FC}$.
\end{corollary}

\begin{proof}
By Lemma~\ref{lemIrreducibility implies isomorphisms}, if $I_\tau$ is irreducible, then $\tau\in \UC_{\FC}$. 

Assume that $\tau\in \UC_{\FC}$. Then by Proposition~\ref{propDecomposition_generalized_weight_spaces_Itau}~(\ref{itDecomposition_Itau_regular_case}), $\dim I_\tau(\tau)=1$ and we conclude with Theorem~\ref{thmIrreducibility criterion} and Remark~\ref{rkAlgebraically closed field}.
\end{proof}

\begin{remark}\label{rkMeasure_irreducible_representations}
Assume that $\FC=\C$ and that $\sigma_s=\sigma_s'=\sqrt{q}$ for all $s\in \SCC$, for some $q\in \Z_{\geq 2}$. Let $(y_j)_{j\in J}$ be a $\Z$-basis of $Y$. Then the map $T_\C\rightarrow (\C^*)^J$ defined by $\tau\in T_\C\mapsto (\tau(y_j))_{j\in J}$ is a group isomorphism. We equip $T_\C$ with a Lebesgue measure  through this isomorphism. Then the set of measurable subsets of $T_\C$ having full measure does not depend on the choice of the $\Z$-basis of $Y$.  Then $\UC_\C=\bigcap_{\alpha^ \vee\in \Phi^\vee} \{\tau\in T_\C|\tau(\alpha^\vee)\neq q\}$ has full measure in $T_\C$. Moreover $T^{\mathrm{reg}}_\C\supset \bigcap_{\lambda\in Y\setminus\{0\}} \{\tau\in T_\C|\tau(\lambda)\neq 1\}$ has full measure in $T_\C$ and thus $\{\tau\in T_\C|I_\tau\text{ is irreducible}\}$ has full measure in $T_\C$.
\end{remark}

Recall that $\RCC=\{wsw^{-1}|w\in W^v,s\in \SCC\}$ is the set of reflections of $W^v$. For $\tau\in T_\C$, set $W_\tau=\{w\in W^v|\ w.\tau=\tau\}$\index{$W_\tau$}, $\Phi^\vee_{(\tau)}=\{\alpha^\vee\in \Phi^\vee_+| \zeta_{\alpha^\vee}^{\mathrm{den}}(\tau)=0\}$\index{$\Phi^\vee_{(\tau)}$},  $\RCC_{(\tau)}=\{r=r_{\alpha^\vee}\in \RCC|\alpha^\vee\in \Phi^\vee_{(\tau)}\}$\index{$\RCC_{(\tau)}$} and \[\Wta=\langle \RCC_{(\tau)}\rangle=\langle \{r=r_{\alpha^\vee}\in \RCC| \zeta_{\alpha^\vee}^{\mathrm{den}}(\tau)=0\}\rangle\subset W^v.\]\index{$\Wta$} By Remark~\ref{rkNecessary_condition_stau_neq_tau}, $\Wta\subset W_\tau$. It is moreover normal in $W_\tau$. When $\alpha_s(Y)=\Z$ for all $s\in \SCC$, then $\Wta=\langle W_\tau\cap \RCC\rangle$.

\begin{corollary}\label{corIrreducibility when Wchi=s}
Let $\tau\in T_\FC$ be such that $W_\tau=\Wta=\{1,t\}$ for some reflection $t$. Then $I_\tau$ is irreducible if and only if $\tau\in \UC_\FC$.
\end{corollary}

\begin{proof}
By Lemma~\ref{lemIrreducibility implies isomorphisms}, if $I_\tau$ is irreducible, then $\tau \in \UC_\FC$. Conversely, let $\tau\in \UC_\FC$ be such that $W_{\tau}=\Wta=\{1,t\}$, for some $t\in \RCC$.  Write $t=v^{-1}sv$ for $s\in \SCC$ and $v\in W^v$.  Let  $\tilde{\tau}=v.\tau$. One has $s.\tilde{\tau}=\tilde{\tau}$ and $W_{\tilde{\tau}}=\{1,s\}$. By Lemma~\ref{lemDecomposition_submodules_I_tau}, $I_{\tilde{\tau}}(\tilde{\tau})\subset I_{\tilde{\tau}}^{\leq s}$. 

Let $\lambda\in Y$. Then $Z^\lambda.H_s\otimes_{\tilde{\tau}}1=\tilde{\tau}(\lambda)H_s\otimes_{\tilde{\tau}} 1+\tilde{\tau}(Q_s(Z)(Z^\lambda-Z^{s.\lambda}))1\otimes_{\tilde{\tau}}1$.

 Suppose $\sigma_s=\sigma_s'$. Then as $W_{(\tilde\tau)}=v.\Wta.v^{-1}=\{1,s\}$, one has  $\tilde{\tau}(\alpha_s^\vee)=1$.  By Remark~\ref{remIH algebre dans le cas KM deploye}, $\tilde{\tau}((Q_s(Z)(Z^\lambda-Z^{s.\lambda}))=(\sigma_s-\sigma_s^{-1})\alpha_s(\lambda)$. As there exists $\lambda\in Y$ such that $\alpha_s(\lambda)\neq 0$, we deduce that $H_s\otimes_{\tilde{\tau}} 1\notin I_{\tilde{\tau}}(\tilde{\tau})$ and thus $ I_{\tilde{\tau}}(\tilde{\tau})=\FC. 1\otimes_{\tilde{\tau}}1$. Similarly, if  $\sigma_s\neq \sigma_s'$ then  $I_{\tilde{\tau}}(\tilde{\tau})=\FC. 1\otimes_{\tilde{\tau}}1$. By Theorem~\ref{thmIrreducibility criterion} and Remark~\ref{rkAlgebraically closed field}, we deduce that $I_{\tilde{\tau}}$  is irreducible. By Lemma~\ref{lemIrreducibility implies isomorphisms} we deduce that $I_\tau$ is isomorphic to $I_{\tilde{\tau}}$ and thus $I_\tau$ is irreducible.
\end{proof}

\subsection{Weight vectors regarded as rational functions}\label{subWeight_vectors_regarded_as_rational_functions}

In this subsection, we introduce and study elements $F_w\in \ATF$, $w\in W^v$, such that for all $\chi\in T_\FC$ such that $F_w(\chi)$ is well-defined, $F_w(\chi)\otimes_\chi 1 \in I_\chi(w.\chi)$. 

 For $w\in W^v$, let $ \pi^T_w:\ATF\rightarrow \FC(Y)$\index{$\pi^T_w$} be the right $\FC(Y)$-module morphism defined by $\pi^T_w(T_v)=\delta_{v,w}$ for all $v\in W^v$.

\begin{lemma}\label{lemDensity_regulars_FC}
Let $\FC'$ be a uncountable field containing $\FC$. Let $P\in \FC[Y]$ be such that $P(\tau)=0$ for all $\tau\in T^{\mathrm{reg}}_{\FC'}$. Then $P=0$.
\end{lemma}

\begin{proof}
Let $\FC_0\subset \FC$ be a countable field (one can take $\FC_0=\Q$ or $\FC_0=\mathbb{F}_\ell$ for some prime power $\ell$). Write $P=\sum_{\lambda\in Y} a_\lambda Z^\lambda$, with $a_\lambda\in \FC$ for all $\lambda \in Y$. Let $(y_j)_{j\in J}$ be a $\Z$-basis of $Y$ and  $X_j=Z^{y_j}$ for all $j\in J$. Let $\FC_1=\FC(a_\lambda|\lambda\in Y)$. Let $(x_j)_{j\in J}\in (\FC')^J$ be algebraically independent over $\FC_1$. Let $\tau\in T_{\FC'}$ be defined by $\tau(y_j)=x_j$ for all $j\in J$. 

 Let us prove that $\tau\in T_\FC^{\mathrm{reg}}$.  Let $w\in W^v\setminus \{1\}$. Let $\lambda\in Y$ be such that $w^{-1}.\lambda-\lambda\neq 0$. Write $w^{-1}.\lambda-\lambda=\sum_{j\in J} n_j y_j$ with $n_j\in \Z$ for all $j\in J$. Let $Q=\prod_{j\in J}Z_j^{n_j}\in \FC_1[Z_j,j\in J]$. Then $Q\neq 1$ and thus $\tau(w^{-1}.\lambda-\lambda)=Q((x_j)_{j\in J})\neq 1$. Thus $w.\tau \neq \tau$ and $\tau\in T_{\FC'}^{\mathrm{reg}}$. Thus $P(\tau)=0$ and by choice of $(x_j)_{j\in J}$ this implies $P=0$.
\end{proof}

Let $w\in W^v$. Let $w=s_1\ldots s_r$ be a reduced expression of $w$. Set $F_w=F_{s_r}\ldots F_{s_1}=(B_{s_r}+\zeta_{s_r})\ldots (B_{s_1}+\zeta_{s_1})\in \AC(T_\FC)$\index{$F_w$}. By the lemma below, this does not depend on the choice of the reduced expression of $w$.

\begin{lemma}\label{lem4.3 de reeder}(see \cite[Lemma 4.3]{reeder1997nonstandard}) Let $w\in W^v$. \begin{enumerate}

\item\label{itWell_definedness_Fw} The element $F_w\in \AC(T_\FC)$ is well-defined, i.e it does not depend on the choice of a reduced expression for $w$.

\item\label{itLeading_coefficient_Fw} One has $F_w-T_w\in \ATF^{<w}$.

\item If $\theta\in \FC(Y)$, then $\theta*F_w=F_w*\theta^{w^{-1}}$.

\item\label{itDomain_Fw} If $\tau\in T_\FC$ is such that $\zeta_{\beta^\vee}\in \FC(Y)_\tau$ for all $\beta^\vee\in N_{\Phi^\vee}(w)$, then $F_w\in \AC(T_\FC)_\tau$ and $F_w(\tau).1\otimes_\tau 1\in I_\tau(w.\tau)$.

\item\label{itF_w_well_defined_at_regular} Let $\tau\in T^{\mathrm{reg}}_\FC$. Then $F_w\in \AC(T_\FC)_\tau$.

\end{enumerate}
\end{lemma}

\begin{proof}
Let us prove~(\ref{itDomain_Fw}) by induction on $\ell(w)$. By Lemma~\ref{lem1.10 of Reeder}, $\theta*F_w=F_w*\theta^{w^{-1}}$ for all $\theta\in \FC(Y)$. Let $n\in \N$ and assume that (\ref{itDomain_Fw}) is true for all $w\in W^v$ such that $\ell(w)\leq n$. Let $w\in W^{v}$ be such that $\ell(w)\leq n+1$. Write $w=sv$, with $s\in \SCC$ and $\ell(v)\leq n$. By Lemma~\ref{lemKumar_1.3.14}, $N_{\Phi^\vee}(sv)=N_{\Phi^\vee}(v)\cup \{v^{-1}.\alpha_s^\vee\}$. Let $\tau\in T_\FC$ be such that  be such that $\zeta_{\alpha^\vee}\in \FC(Y)_\tau$ for all $\alpha^\vee\in N_{\Phi^\vee}(w)$. One has $F_w=(B_{s}+\zeta_s)*F_v$. As $F_v\in \AC(T_\FC)_\tau$ and $\AC(T_\FC)_\tau$ is a left $\HFW$-submodule of $\AC(T_\FC)$, $B_s*F_v\in \AC(T_\FC)_\tau$. One has $\zeta_s*F_v=F_v*\zeta_s^{v^{-1}}\in \AC(T_\FC)_\tau$ and hence $F_w\in \AC(T_\FC)_\tau$. 

Let $\tau\in T_\FC$ be such that $\zeta_{\alpha^\vee}\in \FC(Y)_\tau$ for all $\alpha^\vee\in N_{\Phi^\vee}(w)$. Let $\theta\in \FC[Y]$. Then \[(\theta*F_w)(\tau)=(F_w*\theta^{w^{-1}})(\tau)=\tau(\theta^{w^{-1}})\tau(F_w(\tau)),\] which finishes the proof of~(\ref{itDomain_Fw}).

Let $\tau\in T_{\FC}^{\mathrm{reg}}$ and $\alpha^\vee\in \Phi^\vee$. Write $\alpha^\vee=w.\alpha_s^\vee$ for $w\in W^v$ and $s\in \SCC$. Then $s.w^{-1}.\tau\neq w^{-1}.\tau$ and  by Remark~\ref{rkNecessary_condition_stau_neq_tau}, $w^{-1}.\tau(\zeta_s^{\mathrm{den}})\neq 0$ or equivalently $\tau(\zeta_{\alpha^\vee}^{\mathrm{den}})\neq 0$. By~(\ref{itDomain_Fw}) we deduce that $F_w\in \AC(T_\FC)_\tau$ for all $\tau\in T_{\FC}^ {\mathrm{reg}}$, which proves~(\ref{itF_w_well_defined_at_regular}).

Let us prove~(\ref{itLeading_coefficient_Fw}). Let $v\in W^v$ be such that $h:=F_v-T_v\in \ATF^{<v}$ and $s\in \SCC$ be such that $sv>v$. Then\[F_{sv}=(T_{s}-\sigma_s^2+\zeta_s)*(T_{v}+h)=T_{sv}+(-\sigma_s^2+\zeta_s)*T_v+(-\sigma_s^2+\zeta_s)*h+T_s*h.\] By Lemma~\ref{lemCommutation relation}, \[(-\sigma_s^2+\zeta_s)*T_v,(-\sigma_s^2+\zeta_s)*h\in \AC(T_\FC)^{\leq v}.\] By \cite[Corollary 1.3.19]{kumar2002kac}, $s.[1,v)\subset [1,sv)$ and thus $T_s*h\in \AC(T_\FC)^{< sw}$ thus $F_{sv}-T_{sv}\in \ATF^{<sv}$.  By induction we deduce~(\ref{itLeading_coefficient_Fw}). 

 Let $w=s_1\ldots s_r=s_1'\ldots s_r'$ be reduced expressions of $w$. Let $F_w$ be associated to $s_1\ldots s_r$ and $F_w'$ be associated to $s_1'\ldots s_r'$. Let $\FC'$ be a uncountable field containing $\FC$.  Then by Proposition~\ref{propDecomposition_generalized_weight_spaces_Itau}~(\ref{itDecomposition_Itau_regular_case}), for all $\tau\in T^{\mathrm{reg}}_{\FC'}$ there exists $\theta(\tau)\in \FC'^*$ such that $F_w(\tau)=\theta(\tau)F_w'(\tau)$. Let $v\in W^v$ be such that $\pi^v(F_w')\neq 0$ and $\theta_v=\frac{\pi^H_v(F_w)}{\pi^H_v(F_w')}\in \FC(Y)$. Then $\theta_v(\tau)=\theta(\tau)$ for all $\tau\in T^{\mathrm{reg}}_{\FC'}$.  But by (\ref{itLeading_coefficient_Fw}), $\theta(\tau)=1$ for all $\tau\in T^{\mathrm{reg}}_{\FC'}$. Thus by Lemma~\ref{lemDensity_regulars_FC}, $\theta=1=\theta_v$ and $F_w'=F_w$.
\end{proof}

\begin{remark}\label{rkSupport_F_w}
\begin{enumerate} 
\item When $\sigma_s=\sigma_s'$ for all $s\in \SCC$, the condition~(\ref{itDomain_Fw}) is equivalent to $\tau(\beta^\vee)\neq 1$ for all $\beta^\vee\in N_{\Phi^\vee}(w)$.

\end{enumerate}
\end{remark}

\subsection{One implication of Kato's criterion}\label{subOne implication of Kato's theorem}

Recall the definition of $\Wta$ from Subsection~\ref{subIrreducibility criterion}.

In this subsection, we prove that if $I_\tau$ is irreducible, then $W_\tau=\Wta$.

\begin{lemma}\label{lemKatos theorem reciprocal}
Let $\tau\in T_\C$ be such that $W_\tau\neq\Wta$.  Let $w\in W_\tau\setminus\Wta$ be of minimal length. Then $F_w\in \AC(T_\FC)_\tau$. 
\end{lemma}

\begin{proof}
Write $w=s_k\ldots s_1$, where $k=\ell(w)$ and $s_1,\ldots,s_k\in \SCC$. Let  $j\in \llbracket 0,k-1\rrbracket$. Set $w_j=s_j\ldots s_1$. Suppose that $w_j.\zeta_{s_{j+1}}^{\mathrm{den}}(\tau)=0$. Then $r_{w_j.\alpha_{s_{j+1}}^\vee}=s_1\ldots s_j s_{j+1} s_j\ldots s_1 \in \Wta$. Moreover as $\Wta\subset W_\tau$, we have $s_{j+1}\ldots s_1.\tau= s_{j}\ldots s_1.\tau$.  Therefore \[\tau=w.\tau=s_k\ldots s_j\ldots s_1 .\tau=s_k\ldots \hat{s}_{j+1} \ldots s_1 .\tau,\] and $w'=s_k \ldots \hat{s}_{j+1} \ldots s_1\in W_\tau$. By definition of $w$, $w'\in \Wta$. Consequently \[w=s_k\ldots \hat{s}_{j+1} \ldots s_1.s_1\ldots s_j.s_{j+1}.s_j\ldots s_1=w'r_{w_j.\alpha_{s_{j+1}}^\vee}\in \Wta:\] a contradiction. Therefore $w_j.\zeta_{s_{j+1}}^{\mathrm{den}}(\tau)\neq 0$ and by Lemma~\ref{lemKumar_1.3.14} and Lemma~\ref{lem4.3 de reeder}, $F_w\in \AC(T_\FC)_\tau$.
\end{proof}

\begin{proposition}\label{propKato's weak theorem}
Let $\tau\in T_\C$ be such that  $W_\tau\neq \Wta$. Then $I_\tau$ is reducible. 
\end{proposition}

\begin{proof}
Let $w\in W_\tau\setminus\Wta$ be of minimal length. Then by Lemma~\ref{lemKatos theorem reciprocal} and Lemma~\ref{lem4.3 de reeder}, $F_{w}(\tau)\otimes_\tau 1\in I_\tau(\tau)$. Moreover, $\pi^T_w\big(F_{w}(\tau)\otimes_\tau 1\big)=1$ and thus $F_{w}(\tau)\otimes_\tau 1\notin \C 1\otimes_\tau 1$. We conclude with Theorem~\ref{thmIrreducibility criterion}.
\end{proof}

\subsection{Link with the works of Matsumoto and Kato}\label{subLink with Matsumoto and Kato}

Assume that $W^v$ is finite. Then $\HC_\C=\AC_\C$. Let $\tau\in T_\C$. Then by Subsection~\ref{subPrincipal series representations}, $\dim_\C I_\tau=|W^v|$.  One has $Z^\lambda.1\otimes_\tau 1=\tau(\lambda)1\otimes_\tau 1$ for all $\lambda\in Y$ and $\HC_\C.1\otimes_\tau 1=I_\tau$. Thus by \cite[Th{\'e}or{\`e}me 4.1.10]{matsumoto77Analyse} the definition we used is equivalent to Matsumoto's one.

\medskip

Assume that $\HC_\C$ is associated with a split reductive group over a field with residue cardinal $q$.  Then by (BL2), one has: \[\forall \ s \in \SCC, \forall \ w \in W^{v},\  T_{s}*T_{w}=\left\{\begin{aligned} & T_{sw} &\mathrm{\ if\ }\ell(sw)=\ell(w)+1\\ & (q-1)T_{w}+qT_{s w} &\mathrm{\ if\ }\ell(sw)=\ell(w)-1. \end{aligned}\right . \ \] 

Set $1'_\tau=\sum_{w\in W^v} T_w\otimes_\tau 1$. Then if $s\in \SCC$, $T_s.1_\tau'=q1_\tau'$. Then by \cite[(1.19)]{kato1982irreducibility}, $1'_\tau$ is proportional to the vector $1_\tau$ defined in \cite{kato1982irreducibility}. Kato proves Theorem~\ref{thm*Kato's theorem} by studying whether the following property is satisfied: ``for all $w\in W^v$, $\HC_\C.1_{w.\tau}'=I_{w.\tau}$'' (see \cite[Lemma 2.3]{kato1982irreducibility}). When $W^v$ is infinite, we do not know how to define an analogue of $1'_\tau$ and thus we do not know how to adapt Kato's proof.

\section{Description of generalized weight spaces}\label{secTau_simple_reflections}
In this section, we describe $I_\tau(\tau,\mathrm{gen})$, when $\tau\in T_\C$ is such that $\Wta=W_\tau$. We then deduce Kato's criterion for size $2$ matrices. 

Let us  sketch our proof of this criterion. By Theorem~\ref{thmIrreducibility criterion} and Proposition~\ref{propKato's weak theorem}, it suffices to study $I_\tau(\tau)$ when $\tau\in \UC_\C$ is such that $W_\tau=\Wta$. For this, we begin by describing $I_\tau(\tau,\mathrm{gen})$. Let $\tau\in T_\C$ satisfying the above condition.   By  Dyer's theorem, $(\Wta,\SCC_\tau)$ is a Coxeter system, for some $\SCC_\tau\subset \Wta$. Let $r\in \SCC_\tau$. We study the singularity of $F_r$ at $\tau$, that is, we determine an (explicit) element $\theta\in \C(Y)$ such that $F_r-\theta$ is  defined at $\tau$ (see Lemma~\ref{lemReeder14.3}). Using this, we then describe $I_\tau(\tau,\mathrm{gen})$. We then deduce that when $W_\tau=\Wta$ is the infinite dihedral group then $I_\tau(\tau)$ is irreducible. After classifying the subgroups of the infinite dihedral group (see Lemma~\ref{lemClassification_subgroups_W_dim2}), we deduce Kato's criterion for size $2$ matrices.

\medskip 

In Subsection~\ref{subComplex_torus}, we study the torus $T_\C$.

In Subsection~\ref{subNew_basis_HCW}, we introduce a new basis of $\HCW$ which enables us to have information on the poles of the coefficients of the $F_w$. 

In Subsection~\ref{subExpression_coefficients_Fw_H_v}, we give a recursive formula which enables us to have information on the poles of the coefficients of the $F_w$.

In Subsection~\ref{subTau_simple reflections}, we study the singularity of $F_r$ at $\tau$, for $r \in \SCC_\tau$.

In Subsection~\ref{subDescription_generalized_weight_spaces}, we give a description of $I_\tau(\tau,\mathrm{gen})$, when $W_\tau=\Wta$.

In Subsection~\ref{subIrreducibility_Itau_infinite_dihedral_group}, we prove that when $W_\tau=\Wta$ is the infinite dihedral group and $\tau\in \UC_\C$, then $I_\tau$ is irreducible.

In Subsection~\ref{subKatos_irreducibility_criterion}, we prove Kato's criterion for size $2$ Kac-Moody matrices.

\medskip

 This section is strongly inspired by \cite{reeder1997nonstandard}. 
 
\medskip

 In certain proofs, when $\FC=\C$, we will make additional assumptions on the $\sigma_s$ and $\sigma_s'$, $s\in \SCC$. To avoid these assumptions, we can assume that $\sigma_s,\sigma_s'\in \C$ and $|\sigma_s|>1, |\sigma_s'|>1$ for all $s\in \SCC$. 

\subsection{The complex torus $T_\C$}\label{subComplex_torus}

We assume that $|\sigma_s|\in \R_{>1}$ for all $s\in \SCC$. 
Let $(y_j)_{j\in J}$ be a $\Z$-basis of $Y$. The map $T_\C\rightarrow (\C^*)^{J}$  mapping each $\tau\in T_\C$ on $(\tau(y_j))_{j\in J}$ is a bijection. We identify $T_\C$ and $(\C^*)^J$. We equip $T_\C$ with the usual topology on $(\C^*)^J$. This does not depend on the choice of a basis $(y_j)_{j\in J}$.

\begin{lemma}\label{lemDensity_regular_elements}
The set  $\{\tau\in T_\C|\forall (w,\lambda)\in W^v\setminus\{1\}\times (C^v_f\cap Y),\  w.\tau(\lambda)\neq \tau(\lambda)\}$ is dense in $T_\C$. In particular,  $T_\C^{\mathrm{reg}}$ is dense in $T_\C$.
\end{lemma}

\begin{proof}
Let $\lambda\in C^v_f\cap Y$. By \cite[V.Chap 4 §6 Proposition 5]{bourbaki1981elements}, for all $w\in W^v\setminus \{1\}$, $w.\lambda\neq \lambda$. Let $(\gamma_j)_{j\in J}\in (\C^*)^J$ be algebraically independent over $\Q$ and  $\tau_\gamma\in T_\C$ be defined by $\tau_\gamma(y_j)=\gamma_j$ for all $j\in J$. Then $w.\tau_\gamma(\lambda)\neq \tau_\gamma(\lambda)$ for all $w\in W^v\setminus\{1\}$. Let $\tau\in T_\C$. Let $(\gamma^{(n)})\in \big((\C^*)^J\big)^\N$ be such that $\gamma^{(n)}$ is algebraically independent over $\Q$ for all $n\in \N$ and such that $\gamma^{(n)}\rightarrow (\tau(y_j))_{j\in J}$. Then $\tau_{\gamma^{(n)}}\rightarrow \tau$ and we get the lemma. \end{proof}

Let $A\subset \R$ be a ring. We set $Q^\vee_A=\bigoplus_{s\in \SCC} A\alpha^\vee_s\subset \A$\index{$Q^\vee_A$}.

\begin{lemma}\label{lemExistence_character_values_alphavee}
Let $(\gamma_s)\in (\C^*)^\SCC$. Then there exists $\tau\in T_\C$ such that $\tau(\alpha_s^\vee)=\gamma_s$ for all $s\in \SCC$.
\end{lemma}

\begin{proof}
Let us prove that there exists $n\in \Ne$ such that $\frac1nQ^\vee_\Z\supset Y\cap Q^\vee_\Q$. The module $Y\cap Q^\vee_\Q$ is a $\Z$-submodule of the free module $Y$. Thus it is a free module and its rank is lower or equal to the rank of $Y$. Let $(y_j)_{j\in J}$ be a $\Z$-basis of $Y\cap Q^\vee_\Q$. As $\alpha_s^\vee\in Y\cap Q^\vee_\Q$ for all $s\in \SCC$, we have we have $\mathrm{vect}_\Q(Y\cap Q^\vee_\Q)=Q^\vee_\Q$. Therefore for all $j\in J$, there exists $(m_{j,s})\in \Q^{\SCC}$ such that $y_j=\sum_{j\in J} m_{j,s}\alpha_s^\vee$ and thus there exists $n\in \Ne$ such that $\frac{1}{n}Q^\vee_\Z\supset Y\cap Q^\vee_\Q$. 

Let $S$ be a complement of $Y\cap Q^\vee_\Q$ in $Y\otimes \Q$.  For $s\in \SCC$, choose $\gamma_s^{\frac{1}{n}}\in \C^*$ such that $(\gamma_s^{\frac{1}{n}})^n=\gamma_s$. Let $\tilde{\tau}:\frac{1}{n}Q^\vee_\Z \oplus S\rightarrow \C^*$ be defined by $\tilde{\tau}(\sum_{s\in \SCC} \frac{a_s}{n} \alpha_s^\vee+x)=\prod_{s\in \SCC} (\gamma_s^{\frac{1}{n}})^{a_s}$ for all $(a_s)\in \Z^\SCC$ and $x\in S$. Let $\tau=\tilde{\tau}_{|Y}$. Then $\tau\in T_\C$ and $\tau(\alpha_s^\vee)=\gamma_s$ for all $s\in\SCC$.
\end{proof}

\subsection{A new basis of $\HCW$}\label{subNew_basis_HCW}

In \cite{kazhdan1979representations}, Kazhdan and Lusztig defined the Kazhdan-Lusztig basis $(C_w)_{w\in W^v}$ of $\HCW$ in the case where $\sigma_s=\sigma$ for all $s\in \SCC$. This basis is defined by its properties with respect to some involution of $\HCW$ and by the fact that $C_w-T_w\in\bigoplus_{v<w} \C T_v$ , for $w\in W^v$ (see \cite[Theorem 1.1]{kazhdan1979representations} for a precise statement). This basis was then  defined in the general case (where the $\sigma_s$, $s\in S$ need not be all equal) see \cite[6]{lusztig1983left} for example. We now define a basis $(B_w)_{w\in W^v}$ of $\HCW$ from the Kazhdan-Lusztig basis $(C_w)_{w\in W^v}$ and then compute the coefficient   in front of $B_1$ of the expansion of  $F_w$ in the  basis $(B_v)_{v\in W^v}$, for $w\in W^v$ (see Lemma~\ref{lemConstant_term_intertwining_operators}). This will enable us to have information on the coefficient $\pi^H_1(F_w)\in \C(Y)$, for $w\in W^v$ (see Lemma~\ref{lemConstant_term_intertwining_operators} and Lemma~\ref{lemReeder14.3}). Our computation relies on certain multiplicative properties of $(B_w)$ (see Lemma~\ref{lem5.1 de Reeder}) and we will not need the precise definition of the Kazhdan-Lusztig basis.

Let $(C_w)_{w\in W^v}$ be the basis introduced in \cite[6]{lusztig1983left}. For $w\in W^v$, we set $B_w=(-1)^{\ell(w)} \sigma_w C_w$, where $\sigma_w$ is defined in Remark \ref{remIH algebre dans le cas KM deploye}~(\ref{itWell_definition_sigma_w}). Then for $s\in \SCC$, one has $B_s=T_s-\sigma_s^2$ and thus this notation is coherent with the notation $B_s$ introduced in Subsection~\ref{subintertwining_operators_simple_reflections}.

\begin{lemma}\label{lem5.1 de Reeder}
The  basis $(B_w)_{w\in W^v}$\index{$B_w$} satisfies:\begin{enumerate}

\item\label{itDefBs} $B_s=T_s-\sigma_s^2$ for all $s\in \SCC$,

\item\label{itTriangularity} $B_w-T_w\in \HWC^{< w}$ for all $w\in W^v$,

\item\label{itRelation_product} For all $w\in W^v$ and $s\in \SCC$ we have:\[ B_wB_s=\left\{\begin{aligned} &-(1+\sigma_s^2)B_w &\mathrm{if\ }ws<w\\ & B_{ws}+\sum_{vs<v<w}b(v,w)B_v &\mathrm{if\ }ws>w,\end{aligned}\right.\] for some $b(v,w)\in \C$.
\end{enumerate}
\end{lemma}

\begin{proof}
(2) is a consequence of \cite[2. Proposition]{lusztig1983left}.

(3) Let $w\in W^v$ and $s\in \SCC$ be such that $ws<w$. By \cite[6.4]{lusztig1983left}, $C_w(H_s+\sigma_s^{-1})=0$, thus $(-1)^{\ell(w)}\sigma_wC_w(T_s+1)=0$ and hence $B_w(T_s+1-\sigma_s^2-1)=B_wB_s=-(\sigma_s^2+1)B_w$.

Let $w\in W^v$ and $s\in \SCC$ be such that $ws>w$. Then by \cite[6.3]{lusztig1983left}, one has $C_w(-C_s)\in C_{ws}+\bigoplus_{vs<v<w}\C C_v$ and thus \[(-1)^{\ell(w)}\sigma_w C_w(-\sigma_s C_s)=B_wB_s \in (-1)^{\ell(w)+1}\sigma_{ws} C_{ws}+\bigoplus_{vs<v<w}\C B_v =B_{ws}+\bigoplus_{vs<v<w}\C B_v,\] which proves the lemma.
\end{proof}

As $(B_w)_{w\in W^v}$ is a $\C$-basis of $\HCW$, $(B_w)_{w\in W^v}$ is a $\C(Y)$-basis of the right module $\ATC$.

Let $w\in W^v$. Write $F_w=\sum_{v\in W^v} B_v p_{v,w}$, where $(p_{v,w})\in \C(Y)^{(W^v)}$. By an induction on $\ell(w)$ using Lemma~\ref{lem5.1 de Reeder}~(\ref{itTriangularity}) we have $\bigoplus_{v\leq w}H_v\C(Y)=\bigoplus_{v\leq w} B_v \C(Y)$ for all $w\in W^v$. Thus for all $v\in W^v$ such that $v\nleq w$, one has $p_{v,w}=0$. In \cite[5.3]{reeder1997nonstandard}, Reeder gives recursive formulae for the $p_{v,w}$. The following lemma is a particular case of them.

For $v\in W^v$, define $\pi^B_v:\ATC\rightarrow \C(Y)$\index{$\pi^B_w$} by $\pi^B_v(\sum_{u\in W^v} B_u f_u)=f_v$ for all $(f_u)\in \C(Y)^{(W^v)}$.

\begin{lemma}\label{lemConstant_term_intertwining_operators}
Let $w\in W^v$. Then $p_{1,w}=\zeta_w:=\prod_{\beta^\vee\in N_{\Phi^\vee}(w)}\zeta_{\beta^\vee}$\index{$\zeta_w$}.
\end{lemma}

\begin{proof}
We prove it by induction on $\ell(w)$. 

 Let $v\in W^v$ and assume that $p_{1,v}=\zeta_v$.  Let $s\in \SCC$ be such that $vs>v$. By Lemma~\ref{lem1.10 of Reeder} one has \[\begin{aligned} F_{vs}= & F_v*F_s\\ = & (\sum_{u\in W^v} B_u p_{u,v})*F_s\\ =& \sum_{u\in W^v} B_u*F_s p_{u,v}^s= & \sum_{u\in W^v} B_u *B_s p_{u,v}^s+\sum_{u\in W^v} B_up_{u,v}^s\zeta_s .\end{aligned}\] By Lemma~\ref{lem5.1 de Reeder}, we have $\pi^B_1(\sum_{u\in W^v} B_u *B_s p_{u,v}^s)=0$ and $\pi^B_1(\sum_{u\in W^v}B_up_{u,v}^s\zeta_s)=p_{1,v}^s\zeta_s$. By Lemma~\ref{lemKumar_1.3.14}, $N_{\Phi^\vee}(vs)=s.N_{\Phi^\vee}(v)\sqcup\{\alpha_s^\vee\}$ and thus $\pi^B_1(F_{vs})=p_{1,vs}=p_{1,v}^s\zeta_s=\zeta_{vs}$ which proves the lemma.
\end{proof}

\begin{remark}
In the proof of Lemma~\ref{lemConstant_term_intertwining_operators}, we only used the properties of $(B_w)_{w\in W^v}$ described in Lemma~\ref{lem5.1 de Reeder} and not its precise definition. In \cite[Lemma 5.2]{reeder1997nonstandard}, Reeder gives an elementary proof of the existence of a basis $(B_w)_{w\in W^v}$ satisfying Lemma~\ref{lem5.1 de Reeder}. Its proof can be adapted to our framework to construct a basis $(B_w)$ without using Kazhdan-Lusztig basis.
\end{remark}

\subsection{An expression for the coefficients of the $F_w$ in the basis $(T_v)$}\label{subExpression_coefficients_Fw_H_v}

In this subsection, we give a recursive formula for the coefficients of the $F_w$ in the basis $(T_v)_{v\in W^v}$ (see formula~(\ref{eqDefinitionQv,w}) below and Lemma~\ref{lemExpression_Fw}). We will deduce information concerning the elements $v\in W^v$ such that  $\pi^T_v(F_w)$ is well-defined at $\tau$, for a given $\tau\in T_\C$ (see Lemma~\ref{lemPoints_where_Fw_regular}).

\medskip

 Let $\lambda\in Y$ and $w\in W^v$. By (BL4), Remark~\ref{remIH algebre dans le cas KM deploye}~(\ref{itPolynomiality_Bernstein_Lusztig}) and an induction on $\ell(w)$, there exists $(P_{v,w,\lambda}(Z))_{v\in W^v}\in \C[Y]^{(W^v)}$ such that $Z^\lambda*T_w=\sum_{v\in W^v} T_v*P_{v,w,\lambda}(Z)$. Moreover $P_{w,w,\lambda}=Z^{w^{-1}.\lambda}$ and for all $v\in W^v\setminus[1,w]$,  $P_{v,w,\lambda}=0$.

Let $\lambda\in C^v_f\cap Y$. Then by  \cite[V.Chap 4 §6 Proposition 5]{bourbaki1981elements},  for all $v,w\in W^v$ such that $v\neq w$, one has $v.\lambda\neq w.\lambda$. Let $w\in W^v$. Let $w=s_1\ldots s_k$ be a reduced expression. Set $Q_{w,w,\lambda}(Z)=1\in \C(Y)$. For $v\in W^v\setminus[1,w]$, set $Q_{v,w,\lambda}(Z)=0$. Define $(Q_{v,w,\lambda}(Z))_{v\in [1,w]}$ by decreasing induction by setting: \begin{equation}\label{eqDefinitionQv,w} Q_{v,w,\lambda}(Z)=\frac{1}{Z^{w^{-1}.\lambda}-Z^{v^{-1}.\lambda}}\sum_{w\geq u> v}Q_{u,w,\lambda} P_{v,u,\lambda}\in \C(Y).\end{equation}

\begin{lemma}\label{lemCaracterization_weight_vector_regular_case}
Let $\lambda\in C^v_f\cap Y$, $w\in W^v$ and $\tau\in T_\C^\mathrm{reg}$ be such that $v.\tau(\lambda)\neq \tau(\lambda)$ for all $v\in W^v\setminus\{1\}$. Let $x\in I_\tau$ be such that $Z^\lambda.x=w.\tau(\lambda).x$. Then $x\in I_\tau(w.\tau)$.
\end{lemma}

\begin{proof}
By Proposition~\ref{propDecomposition_generalized_weight_spaces_Itau}~(\ref{itDecomposition_Itau_regular_case}), we can write $x=\sum_{v\in W^v} x_v$ where $x_v\in I_\tau(v.\tau)$ for all $v\in W^v$. One has $Z^\lambda.x-w.\tau(\lambda).x=0=\sum_{v\in W^v}(v.\tau(\lambda)-w.\tau(\lambda))x_v$. As $v.\tau(\lambda)\neq w.\tau(\lambda)$ for all $v\neq w$, we deduce that $x=x_w$.
\end{proof}

\begin{lemma}\label{lemExpression_Fw}
Let $v,w\in W^v$. Then $\pi^T_{v}(F_w)=Q_{v,w,\lambda}$, for any $\lambda\in C^v_f\cap Y$. In particular, $Q_{v,w,\lambda}$ does not depend on the choice of $\lambda$.
\end{lemma}

\begin{proof}
Let $\lambda\in C^v_f$ and $h=\sum_{v\in W^v}T_vQ_{v,w,\lambda}\in \ATC$. One has: \[ \begin{aligned} Z^\lambda*h &=  Z^\lambda *\sum_{v\in W^v} T_v Q_{v,w,\lambda} \\ & = \sum_{u,v\in W^v} T_u P_{u,v,\lambda} Q_{v,w,\lambda} \\ & =\sum_{u\in W^v } T_u \sum_{v\in W^v} P_{u,v,\lambda}Q_{v,w,\lambda}.\end{aligned}\]

Let $u\in W^v$. Then: \[\begin{aligned}\sum_{v\in W^v}P_{u,v,\lambda}Q_{v,w,\lambda}&=P_{u,u,\lambda}Q_{u,w,\lambda}+\sum_{v>u} P_{u,v,\lambda}Q_{v,w,\lambda}\\ &=Z^{u^{-1}.\lambda}+(Z^{w^{-1}.\lambda}-Z^{u^{-1}.\lambda})Q_{u,w,\lambda}\\ &=Z^{w^{-1}.\lambda}Q_{u,w,\lambda},\end{aligned}\] and therefore $Z^{\lambda}.h=h. Z^{w^{-1}.\lambda}$.

Let $\lambda\in C^v_f\cap Y$ and  $\tau\in T^{\mathrm{reg}}_\C$ be such that $u.\tau(\lambda)\neq \tau(\lambda)$ for all $u\in W^v\setminus\{1\}$. Then $\ev_\tau(Z^\lambda*h)=\ev_\tau(h*Z^{v^{-1}.\lambda})=w.\tau(\lambda).h(\tau)$. By Lemma~\ref{lemCaracterization_weight_vector_regular_case} we deduce that $h(\tau)\in I_\tau(w.\tau)$. By Proposition~\ref{propDecomposition_generalized_weight_spaces_Itau}~(\ref{itDecomposition_Itau_regular_case}) and Lemma~\ref{lem4.3 de reeder} we deduce that $h(\tau)=F_w(\tau)$. By Lemma~\ref{lemDensity_regular_elements}, we deduce that $h=F_w$, which proves the lemma.
\end{proof}

\begin{lemma}\label{lemPoints_where_Fw_regular}
Let $w\in W^v$, $\tau\in T_\C$ and $v\in [1,w]$. Assume that for all $u\in [v,w)$, $u.\tau \neq w.\tau$. Then for all $u\in [v,w]$, $\pi^T_u(F_w)\in \C(Y)_\tau$.
\end{lemma}

\begin{proof}
We do it by decreasing induction on $v$. Suppose that for all $u\in (v,w)$, $\pi^T_{u}(F_w)\in \C(Y)_\tau$. Let $\lambda\in C^v_f\cap Y$ be such that $v.\tau(\lambda)\neq w.\tau(\lambda)$, which exists because $C^v_f\cap Y$ generates $Y$. By Lemma~\ref{lemExpression_Fw} we have \[\pi^T_v(F_w)=Q_{v,w,\lambda}=\frac{1}{Z^{w^{-1}.\lambda}-Z^{v^{-1}.\lambda}}\sum_{w\geq u>v} Q_{u,w,\lambda}P_{v,u,\lambda}.\] We deduce that $\pi^T_v(F_w)\in \C(Y)_\tau$ because by assumption $Q_{u,w,\lambda}\in\C(Y)_\tau$ for all $u\in (v,w]$. Lemma follows.
\end{proof}

\subsection{$\tau$-simple reflections and intertwining operators}\label{subTau_simple reflections}

Let $\tau\in T_\C$. Following \cite[14]{reeder1997nonstandard}, we introduce $\tau$-simple reflections (see Definition~\ref{defTau_simple_reflections}). If $\SCC_\tau$ is the set of $\tau$-simple reflections,  then $(\Wta,\SCC_\tau)$ is a Coxeter system. We study, for such a reflection $r$, the singularity of $F_r$ at $\tau$: we prove that $F_r-\zeta_r$ is in $\ATC_\tau$ (see Lemma~\ref{lemReeder14.3}). This enables us to define $K_r(\tau)=(F_r-\zeta_r)(\tau)\in \HCW$.  This will be useful to describe $I_\tau(\tau,\mathrm{gen})$.  

We now define $\tau$-simple reflections. Our definition slightly differs from \cite[Definition 14.2]{reeder1997nonstandard}. These definitions are equivalent (see Lemma~\ref{lemComparison_definition_tau_simples}).

\begin{definition}\label{defTau_simple_reflections}
Let $\tau\in T_\C$. A coroot $\beta^\vee\in \Phi^\vee_\tau$ and its corresponding reflection $r_{\beta^\vee}$ are said to be \textbf{$\tau$-simple} if $N_{\RCC}(r_{\beta^\vee})\cap \Wta=\{r_{\beta^\vee}\}$. We denote by $\SCC_\tau$\index{$\SCC_\tau$} the set of $\tau$-simple reflections.
\end{definition}

Recall that $\Phi^\vee_{(\tau)}=\{\alpha^\vee\in \Phi^\vee_+|\zeta^{\mathrm{den}}_{\alpha^\vee}(\tau)=0\}$ and $\RCC_{(\tau)}=\{r_{\alpha^\vee}|\alpha^\vee\in \Phi^\vee_{(\tau)}\}$.

\subsubsection{Coxeter structure of $W_{(\tau)}$ and comparison of the definitions of $\tau$-simplicity}

We use the same notation as in~\ref{subsubReflections_subgroups}. Then $\SCC_\tau=\SCC(\Wta)$ and thus $(\Wta,\SCC_\tau)$ is a Coxeter system.

Let $\leq_\tau$\index{$\leq_\tau$} and $\ell_\tau$\index{$\ell_\tau$} be the Bruhat order  and the length on $(W_{(\tau)},\SCC_\tau)$. 

\begin{lemma}\label{lemComparison_Bruhat_orders}
Let $x,y\in W_{(\tau)}$ be such that $x\leq_\tau y$. Then $x\leq y$.
\end{lemma}

\begin{proof}
By definition, if $x,y\in W_{(\tau)}$, then $x\leq_\tau y$ (resp. $x\leq y$) if there exist $n\in \N$ and  $x_0=x,x_1, \ldots ,x_n=y\in W_{(\tau)}$ (resp . $W^v$) such that $(x_i,x_{i+1})$ is an arrow of the graph of \cite[Definition 1.1]{dyer1991bruhat} for all $i\in \llbracket 0,n-1\rrbracket$. We conclude with \cite[Theorem 1.4]{dyer1991bruhat} 
\end{proof}

\begin{remark}
The orders $\leq$ and $\leq_\tau$ can be different on $\Wta$: there can exist $v,w\in \Wta$ such that $v$ and $w$ are not comparable for $\leq_\tau$ and $v< w$. For example if $W^v=\{s_1,s_2\}$ is the infinite dihedral group, $r_1=s_1$ and $r_2=s_2s_1s_2$ (see Lemma~\ref{lemTau_fixer_2_generators}), then $r_1<r_2$ but $r_1$ and $r_2$ are not comparable for $<_\tau$.
\end{remark}

Set $\Phi^\vee_{(\tau),+}=\Phi^\vee_{(\tau)}\cap \Phi^\vee_+$ and $\Phi^\vee_{(\tau),-}=\Phi^\vee_{(\tau)}\cap\Phi^\vee_-$. For $w\in W_{(\tau)}$, set $N_{\Phi^\vee_{(\tau)}}(w)=N_{\Phi^\vee}(w)\cap \Phi^\vee_{(\tau)}$.

\begin{lemma}\label{lemPhiveetau_stable_under_Wta}
Let $w\in \Wta$. Then $w.\Phi^\vee_{(\tau)}=\Phi^\vee_{(\tau)}$ and $w.\RCC_{(\tau)}.w^{-1}=\RCC_{(\tau)}$.
\end{lemma}

\begin{proof}
Let $\alpha^\vee\in \Phi^\vee_{(\tau)}$. One has $\zeta_{w.\alpha^\vee}^{\mathrm{den}}=(\zeta_{\alpha^\vee}^{\mathrm{den}})^w$ and hence \[\zeta_{\alpha^\vee}^{\mathrm{den}}(\tau)=(\zeta_{\alpha^\vee}^{\mathrm{den}})^w(\tau)=(\zeta_{\alpha^\vee}^{\mathrm{den}})(w^{-1}.\tau)=0\] because $w\in \Wta\subset W_\tau$. Thus $w.\alpha^\vee\in \Phi^\vee_{(\tau)}$ and $r_{v.\alpha^\vee}=vr_{\alpha^\vee}v^{-1}\in \RCC_{(\tau)}$, which proves the lemma.
\end{proof}

We now prove that our definition of $\tau$-simplicity is equivalent to the definition of \cite[14.2]{reeder1997nonstandard}. This equivalence will be useful in our study of the weight spaces of $I_\tau$ and thus in the study of the irreducibility of $I_\tau$. Indeed, our definition of $\tau$-simplicity is well adapted to the study of the Coxeter structure of $W_{(\tau)}$ whereas Reeder’s one is well adapted to the study of the singularity $F_r$ at $\tau$.

\begin{lemma}\label{lemComparison_definition_tau_simples}
\begin{enumerate}
\item One has $\SCC_\tau\subset\RCC\cap \Wta=\RCC_{(\tau)}$.

\item Let $r=r_\beta^\vee\in \RCC$. Then $r\in \SCC_\tau$ if and only if $N_{\Phi^\vee}(r_{\beta^\vee})\cap \Phi^\vee_{(\tau)}=\{\beta^\vee\}$.

\item Let $w\in \Wta$. Let $w=r_1\ldots r_k$ be  a reduced writing of $\Wta$, with $k=\ell_\tau(w)$ and $r_1,\ldots,r_k\in \SCC_\tau$. Then $|N_{\Phi^\vee_{(\tau)}}(w)|=\{\alpha^\vee_{r_k},r_k.\alpha^\vee_{r_{k-1}},\ldots,r_k\ldots r_2.\alpha_{r_1}^\vee\}$ and $|N_{\Phi^\vee}(w)\cap \Phi^\vee_{(\tau)}|=k=\ell_\tau(w)$.

\end{enumerate}

\end{lemma}

\begin{proof}

We begin by proving a part of (3).  By Lemma~\ref{lemComparison_Bruhat_orders} and \cite[Lemma 1.3.13]{kumar2002kac}, for $v\in \Wta$ and $r\in \SCC_\tau$, one has $\ell_\tau(vr)>\ell_\tau(v)$ if and only if $vr>_\tau v$ if and only if $vr>v$ if and only if  $v.\alpha_r^\vee\in \Phi^\vee_+$ if and only if $v.\alpha_r^\vee\in \Phi^\vee_{(\tau),+}$.

One has $N_{\Phi^\vee_{(\tau)}}(w)=\{\alpha^\vee\in \Phi^\vee_{(\tau),+}|w.\alpha^\vee\in \Phi^\vee_{(\tau),-}\}$. Then using the same proof as in \cite[Lemma 1.3.14]{kumar2002kac}, one has $N_{\Phi^\vee_{(\tau)}}(w)\supset \{\alpha^\vee_{r_k},r_k.\alpha^\vee_{r_{k-1}},\ldots,r_k\ldots r_2.\alpha_{r_1}^\vee\}$ and 

$|\{\alpha^\vee_{r_k},r_k.\alpha^\vee_{r_{k-1}},\ldots,r_k\ldots r_2.\alpha_{r_1}^\vee\}|=k=\ell_\tau(w)$.

We now prove (1) and (2). Let $f:\Phi^\vee_+\rightarrow \RCC$ be the map defined by $f(\alpha^\vee)=r_{\alpha^\vee}$ for  $\alpha^\vee\in \Phi^\vee_+$. Then by Subsection~\ref{subReflection_subgroups}, $f$ is a bijection. Let $r=r_{\beta^\vee}\in \SCC_\tau$. One has $f\big(N_{\Phi^\vee}(r)\cap \Phi^\vee_{(\tau)}\big)=N_{\RCC}(r)\cap \RCC_{(\tau)}$. Moreover, $\RCC_{(\tau)}\subset \Wta\cap \RCC$. Thus \[f^{-1}\big(N_\RCC(r)\cap \Wta\big)=\{\beta^\vee\}\supset f^{-1}\big(N_{\RCC}(r)\cap \RCC_{(\tau)}\big)=N_{\Phi^\vee}(r)\cap \Phi^\vee_{(\tau)}.\] Moreover, $|N_{\Phi^\vee}(r)\cap \Phi^\vee_{(\tau)}|\geq 1$ and thus $|N_{\Phi^\vee}(r)\cap \Phi^\vee_{(\tau)}|=\{\beta^\vee\}$. In particular, $\beta^\vee\in \Phi^\vee_{(\tau)}$ and $r\in \RCC_{(\tau)}$. Thus $\SCC_\tau\subset \RCC_{(\tau)}$.

By \cite[Theorem 3.3 (i)]{dyer1990reflection}, $\RCC\cap \Wta=\bigcup_{w\in \Wta} w\SCC_\tau w^{-1}$ and thus by Lemma~\ref{lemPhiveetau_stable_under_Wta},  $\RCC\cap \Wta\subset \bigcup_{w\in \Wta} w.\RCC_{(\tau)}.w^{-1}=\RCC_{(\tau)}$. As by definition, $\RCC_{(\tau)}\subset \Wta\cap \RCC$, we deduce that $\RCC_{(\tau)}=\Wta\cap \RCC$, which proves (1). 

Let $r=r_\beta^\vee\in \RCC$. Suppose that $N_{\Phi^\vee}(r_{\beta^\vee})\cap \Phi^\vee_{(\tau)}=\{\beta^\vee\}$. Then \[f\big(N_{\Phi^\vee}(r_{\beta^\vee})\cap \Phi^\vee_{(\tau)}\big)=\{r_{\beta^\vee}\}=N_{\RCC}(r_{\beta^\vee})\cap \RCC_{(\tau)}=N_{\RCC}(r_{\beta^\vee})\cap \Wta ,\] which proves (2). 

Let $\alpha^\vee\in N_{\Phi^\vee_{(\tau)}}(w)$. Then there exists $j\in \llbracket 2,k\rrbracket$ such that $r_j\ldots r_k.\alpha^\vee \in \Phi^\vee_{(\tau),+}$ and $r_{j-1}\ldots r_k.\alpha^\vee\in \Phi^\vee_{(\tau),-}$. Thus $r_{j-1}\ldots r_k.\alpha^\vee\in N_{\Phi^\vee_{(\tau)}}(r_j)=\{\alpha_{r_j}^\vee\}$ and hence $\alpha^\vee=r_k\ldots r_{j-1}.\alpha_{r_j}^\vee$, which concludes the proof of the lemma.
\end{proof}

\subsubsection{Singularity of $F_r$ at $\tau$ for a $\tau$-simple reflection}

\begin{lemma}\label{lemFactorization_Fr_simple_reflection}
Let $\tau\in T_\C$ and $r_{\beta^\vee}\in \SCC_\tau$. Then there exists $h'\in \ATC_\tau$ such that $F_{r_{\beta^\vee}}=h'.\zeta_{\beta^\vee}^{\mathrm{den}}$.
\end{lemma}

\begin{proof}
Using \cite[1. Exercise 10]{bjorner2005combinatorics}, we write $r_{\beta^\vee}=wsw^{-1}$ with $w\in W^v$, $s\in \SCC$ and $\ell(wsw^{-1})=2\ell(w)+1$. One has $\beta^\vee=w.\alpha_s^\vee$. Let $r_{\beta^\vee}=s_m\ldots s_1$  be a reduced  expression of $r_{\beta^\vee}$, with $m\in\N$ and  $s_1,\ldots,s_m\in \SCC$. Let $k\in \llbracket 0,m-1\rrbracket$ and $v=s_{k}\ldots s_1$. Suppose that  $F_v=h'_k.(\zeta_{\beta^\vee}^{\mathrm{den}})^{\eta(k)}$ where $h'_k\in \ATC_\tau$ and $\eta(k)\in \N$. Then $F_{s_{k+1}v}=F_{s_{k+1}}*F_{v}=(B_{s_{k+1}}+\zeta_{s_{k+1}})*F_{v}$. One has  $\zeta_{s_{k+1}}*F_{v}=F_v.\zeta_{s_{k+1}}^{v^{-1}}$ by Lemma~\ref{lem4.3 de reeder}. 

By Lemma~\ref{lemComparison_definition_tau_simples}  if $\zeta_{s_{k+1}}^{v^{-1}}$ is not defined in $\tau$  then $k=\ell(w)$. As $B_{s_{k+1}}\in \HCW$ and $\ATC_\tau$ is a left $\HCW$-module, we can write $F_{s_{k+1}v}=h'_{k+1}.(\zeta_{\beta^\vee}^{\mathrm{den}})^{\eta(k+1)}$ where $h'_{k+1}\in \ATC_\tau$ and $\eta(k+1)\leq \eta(k)$ if $k\neq \ell(w)$ and $\eta(k+1)\leq \eta(k)+1$ if $k=\ell(w)$, which proves the lemma.
\end{proof}

\begin{lemma}\label{lemComparison_poles_H_B}
Let $h\in \ATC$ and $\tau\in T_\C$. Then \[\max\{ u\in W^v|\pi^H_u(h)\notin \C(Y)_\tau\}=\max\{ u\in W^v|\pi^B_u(h)\notin \C(Y)_\tau\}.\]
\end{lemma}

\begin{proof}
Let $v\in \max\{ u\in W^v|\pi^H_u(h)\notin \C(Y)_\tau\}$.  By~\ref{lem5.1 de Reeder}~(\ref{itTriangularity}), \[\pi^B_v(h)=\sum_{u\geq v} \pi^{B}_v(H_u)\pi^H_u(h)=\pi^B_v(H_v)\pi^H_v(h)+\sum_{u>v}\pi^{B}_v(H_u)\pi^H_u(h).\] Moreover, by Lemma~\ref{lem5.1 de Reeder}~(\ref{itDefBs}) $\pi^B_v(H_v)\in \C^*$. Thus $\pi^B_v(h)\notin \C(Y)_\tau$. Similarly if $v'\in \max\{u\in W^v, u\geq v| \pi^B_u(h)\notin \C(Y)_\tau\}$, then $\pi^H_{v'}(h)\notin \C(Y)_\tau$.  Hence $v\in\max\{ u\in W^v|\pi^B_u(h)\notin \C(Y)_\tau\}$ and consequently $\max\{ u\in W^v|\pi^H_u(h)\notin \C(Y)_\tau\}\subset \max\{ u\in W^v|\pi^B_u(h)\notin \C(Y)_\tau\}$. By a similar reasoning we get the other inclusion.
\end{proof}

\begin{lemma}\label{lemCaracterisation_simple_reflections} 
Let $w\in W^v$. Suppose that for some $s\in \SCC$, we have $w.\lambda-\lambda\in \R \alpha_s^\vee$ for all $\lambda\in Y$. Then $w\in\{\Id, s\}$.
\end{lemma}

\begin{proof}
Let $\beta^\vee\in N_{\Phi^\vee}(w)$. Write $\beta^\vee=\sum_{t\in \SCC}n_t\alpha_t^\vee$, with $n_t\in\N$ for all $t\in \SCC$. Then $w.\beta^\vee\in \Phi^\vee_-$ and by assumption, $n_t=0$ for all $t\in \SCC\setminus\{s\}$. Therefore $\beta^\vee\in \N \alpha_s^\vee \cap \Phi^\vee=\{\alpha_s ^\vee\}$. We conclude with Lemma~\ref{lemKumar_1.3.14}.
\end{proof}

\begin{lemma}\label{lemDensity_characters}
Let $\chi\in T_\C$. Assume that there exists $\beta^\vee\in \Phi^\vee_+$ such that $r_{\beta^{\vee}}\in W_\chi$. Then there exists $(\chi_n)\in (T_\C)^\N$ such that:\begin{itemize}
\item $\chi_n\rightarrow \chi$,

\item $W_{\chi_n}=\langle r_{\beta^\vee}\rangle$ for all $n\in \N$,

\item $\chi_n(\beta^\vee)=\chi(\beta^\vee)$ for all $\in \N$. 
\end{itemize}
\end{lemma}

\begin{proof}
We first assume that $\beta^\vee=\alpha_s^\vee$, for some $s\in\SCC$. Let $(y_j)_{j\in J}$ be a $\Z$-basis of $Y$. For all $j\in J$, choose $z_j\in \C$ such that $\chi(y_j)=\exp(z_j)$. Let $g:\A\rightarrow \C$ be the linear map such that $g(y_j)=z_j$ for all $j\in J$.  Let $V$ be a complement of $Q^\vee_\R$ in $\A$. Let $n\in \Ne$.
 Let $b_s^{(n)}=g(\alpha_s^\vee)$ and $(b_t^{(n)})\in \C^{\SCC\setminus\{s\}}$ be such that $|b_t^{(n)}-g(\alpha_t^\vee)|<\frac{1}{n}$ and such that the $\exp(b_t^{(n)})$, $t\in \SCC\setminus\{s\}$ are algebraically independent over $\Q$.  Let $g_n:\A\rightarrow \C$ be the linear map such that $g_n(\alpha_t^\vee)=b_t^{(n)}$ for all $t\in \SCC$ and $g_n(v)=g(v)$ for all $v\in V$. For $n\in \N$ set $\chi_n=(\exp\circ g_n)_{|Y}\in T_\C$. For all $x\in \A$, $g_n(x)\rightarrow g(x)$ and thus $\chi_n\rightarrow \chi$. 

Let $n\in \Ne$. Then $\chi(\alpha_s^\vee)=\chi_n(\alpha_s^\vee)$ and thus $s\in W_{\chi_n}$. Let $w\in W_{\chi_n}$. Then $w^{-1}.\lambda-\lambda\in \Z\alpha_s^\vee$ for all $\lambda\in Y$. By Lemma~\ref{lemCaracterisation_simple_reflections} we deduce  that $w\in \{\Id, s\}$. Therefore $W_{\chi_n}=\{\Id, s\}$.

We no more assume that  $\beta^\vee=\alpha_s^\vee$ for some $s\in \SCC$. Write $\beta^\vee=w.\alpha_s^\vee$ for some $w\in W^v$ and $s \in \SCC$. Let $\tilde{\chi}=w^{-1}.\chi$. Then $s\in W_{\tilde{\chi}}$.  Thus there exists $(\tilde{\chi}_n)\in (T_\C)^\N$ such that $\tilde{\chi}_n \rightarrow \tilde{\chi}$ and $W_{\tilde{\chi}_n}=\{\Id, s\}$ for all $n\in \N$. Let $(\chi_n)=(w.\tilde{\chi}_n)$. Then $\chi_n\rightarrow \chi$ and $W_{\chi_n}=\{1,r_{\beta^\vee}\}$ for all $n\in \N$. 

Moreover, $\chi(\beta^\vee)\in \{-1,1\}$ and $\chi_{n}(\beta^\vee)\in \{-1,1\}$ for all $n\in \N$. Maybe considering a subsequence of $(\chi_n)$, we may assume that there exists $\epsilon\in \{-1,1\}$ such that $\chi_n(\beta^\vee)=\epsilon$ for all $n\in \N$. As $\chi_n\rightarrow \chi$, $\chi_n(\beta^\vee)=\epsilon\rightarrow \chi(\beta^\vee)$, which proves the lemma. 
\end{proof}

Let $\C[Q^\vee_\Z]=\bigoplus_{\lambda\in Q^\vee_\Z}\C Z^\lambda \subset \C[Y]$. This is the group algebra of $Q^\vee_\Z$. Let $\C(Q^\vee_\Z)\subset \C(Y)$ be the field of fractions of $\C[Q^\vee_\Z]$ and $\HC(Q^\vee_\Z)=\bigoplus_{w\in W^v} H_w \C(Q^\vee_\Z)\subset \ATC$. This is a ($\HCW-\C(Q^\vee_\Z)$)-bimodule of $\ATC$ and a left $\C(Q^\vee_\Z)$-submodule of $\ATC$. Consequently $F_w\in \HC(Q^\vee_\Z)$ for all $w\in W^v$.

Let $A=\C[Z^{\alpha_s^\vee}|s\in \SCC]\subset \C[Q^\vee_\Z]$. This is a unique factorization domain and $\C(Q^\vee_\Z)$ is the  field of fractions of $A$.

\begin{lemma}\label{lemIrreducibility_Zbeta-1}
Let $\beta^\vee\in \Phi^\vee$. Then $Z^{\beta^\vee}-1$ and $Z^{\beta^\vee}+1$ are irreducible in $A$.
\end{lemma}

\begin{proof}
Write $\beta^\vee=w.\alpha_s^\vee$, where $w\in W^v$ and $s\in \SCC$. Then $Z^{\beta^\vee}=(Z^{\alpha_s^\vee})^w$.
\end{proof} 

\begin{lemma}\label{lemReeder14.3}(see \cite[Proposition 14.3]{reeder1997nonstandard})
Let  $\tau\in T_\C$ and $r=r_{\beta^\vee}\in \SCC_\tau$. Then $F_{r_\beta^\vee}-\zeta_{\beta^\vee}\in \ATC_\tau$.
\end{lemma}

\begin{proof}
One has $ F_{r_\beta^\vee}-\zeta_{\beta^\vee}\in \HC(Q^\vee_\Z)$. Write $F_{r_\beta^\vee}-\zeta_{\beta^\vee}=\sum_{u\in W^v} H_u \frac{f_u}{g_u}$, with $f_u,g_u\in A$  and $f_u\wedge g_u=1$ for all $u\in W^v$. Let $u\in (1,r_{\beta^\vee})$. Let us prove that $\zeta_{\beta^\vee}^{\mathrm{den}}\wedge g_u=1$. Suppose that $\zeta_{\beta^\vee}^{\mathrm{den}}\wedge g_u\neq 1$.   Then there exists $\eta\in \{-1,1\}$ such that  $Z^{\beta^\vee}+\eta$ divides $g_u$. 

Let $\chi\in T_\C$ be such that $\chi(\beta^\vee)=-\eta$. By Remark~\ref{rkNecessary_condition_stau_neq_tau}, $r_{\beta^\vee}\in W_\chi$.   Let $(\chi_n)\in (T_\C)^\N$ be such that $\chi_n\rightarrow \chi$ and $W_{\chi_n}=\{1,r_{\beta^\vee}\}$ for all $n\in \N$, and $\chi_n(\beta^\vee)=-\eta$ for all $n\in \N$. 
 whose existence is provided by Lemma~\ref{lemDensity_characters}. 
One has $g_u(\chi_n)=0$ for all $n\in \N$. Moreover by Lemma~\ref{lemPoints_where_Fw_regular}, $ \pi^H_u(F_{r_\beta^\vee})=\frac{f_u}{g_u}\in \C(Y)_{\chi_n}$ for all $n\in \N$. Therefore, $f_u(\chi_n)=0$ for all $n\in \N$ and thus $f_u(\chi)=0$. 

 By the Nullstellensatz (see~\cite[IX, Theorem~1.5]{lang2002algebra} for example), there exists $n\in \N$ such that $Z^{\beta^\vee}+\eta$ divides $f_u^n$ in $A$. By Lemma~\ref{lemIrreducibility_Zbeta-1}, $Z^{\beta^\vee}+\eta$ is irreducible in $A$ and thus $Z^{\beta^\vee}+\eta$ divides $f_u$: a contradiction. Therefore   $\zeta_{\beta^\vee}^{\mathrm{den}}\wedge g_u=1$. By Lemma~\ref{lemFactorization_Fr_simple_reflection}, $g_u(\tau)\neq 0$.
 
Therefore $\{u\in W^v| \pi^H_u(F_{r_{\beta^\vee}}-\zeta_{r_{\beta^\vee}})\notin \C(Y)_\tau\}\subset \{1\}$. By Lemma~\ref{lemComparison_poles_H_B} we deduce that  $\{u\in W^v| \pi^B_u(F_{r_{\beta^\vee}}-\zeta_{r_{\beta^\vee}})\notin \C(Y)_\tau\}\subset \{1\}$. Using Lemma~\ref{lemConstant_term_intertwining_operators} we deduce that $\{u\in W^v| \pi^B_u(F_{r_{\beta^\vee}}-\zeta_{r_{\beta^\vee}})\notin \C(Y)_\tau\}=\emptyset$. By Lemma~\ref{lemComparison_poles_H_B}, $\{u\in W^v| \pi^H_u(F_{r_{\beta^\vee}}-\zeta_{r_{\beta^\vee}})\notin \C(Y)_\tau\}=\emptyset$, which proves the lemma. 
\end{proof}

\subsection{Description of generalized weight spaces }\label{subDescription_generalized_weight_spaces}

In this subsection, we describe $I_\tau(\tau,\mathrm{gen})$ for $\tau\in \UC_\C$ when $\Wta=W_\tau$, using the $K_{r_1}\ldots K_{r_k}(\tau)$, for $r_1,\ldots, r_k\in \SCC_\tau$  (see Theorem~\ref{thmGeneralized_weigth_space}).

For $r\in \RCC$, one sets $K_{r}=F_{r}-\zeta_{\alpha_r^\vee}\in \ATC$\index{$K_r$}. By Lemma~\ref{lem4.3 de reeder} we have: \begin{equation}\label{eqRelation_commutation_K}
 \theta* K_r=K_r*\theta^r+(\theta^r-\theta)\zeta_r \text{ for all }\theta\in \C(Y).
 \end{equation}

\begin{lemma}\label{lemProductFuFv}
Let $w_1,w_2\in W^v$. Then there exists $P\in \C(Y)^\times$ such that $F_{w_1}*F_{w_2}=F_{w_1w_2} *P$. If moreover $\tau\in \UC_\C$, then one can write $P=\frac{f}{g}$ with $f,g\in \C[Y]^\times$ and $f(w.\tau)\neq 0$ for all $w\in W^v$.
\end{lemma}

\begin{proof}
Let $u,v\in W^v$. Let us prove that if $\chi\in T_\C^{\mathrm{reg}}$, then $F_u*F_v\in \ATC_\chi$. Write $F_u=\sum_{u'\leq u} H_{u'}\theta_{u'}$, where $\theta_{u'}\in \C(Y)$ for all $u'\leq u$. Then by Lemma~\ref{lem4.3 de reeder}, \[F_u*F_v=\sum_{u'\leq u} H_{u'}\theta_{u'}*F_v=\sum_{u'\leq u}H_{u'}*F_v*(\theta_{u'})^{v^{-1}}.\] By Lemma~\ref{lem4.3 de reeder}, $\theta_{u'}\in \ATC_{\chi}$ for all $\chi\in T_\C^{\mathrm{reg}}$ and thus $(\theta_{u'})^{v^{-1}}\in \ATC_\chi$ for all $\chi\in T_\C^{\mathrm{reg}}$. Let $\chi\in T_\C^{\mathrm{reg}}$. As $\ATC_\chi$ is an $\HCW-\C(Y)_\chi$ bimodule, we deduce that $F_u*F_v\in \ATC_\chi$.

Let $u,v\in W^v$. Let us prove that there exists $Q\in \C(Y)$ such that $F_u*F_v=F_{uv}*Q$. Let $\lambda\in  Y$. Then by Lemma~\ref{lem4.3 de reeder}, one has $Z^\lambda F_u*F_v=F_u*F_v*Z^{(uv)^{-1}.\lambda}$. Therefore for all $\chi\in  T^{\mathrm{reg}}_\C$, there exists $a(\chi)\in \C$ such that  $F_u*F_v(\chi)=a(\chi)F_{uv}(\chi)$. Write $F_u*F_v=\sum_{w\in W^v}H_w*\theta_w$ and $F_{uv}=\sum_{w\in W^v} H_w*\tilde{\theta}_w$, where $(\theta_w),(\tilde{\theta}_w)\in \C(Y)^{(W^v)}$. Let $Q=\frac{\theta_{uv}}{\tilde{\theta}_{uv}}={\theta}_{uv}$. Let $w\in W^v$ be such that $\tilde{\theta}_w=0$. Then for all $\chi\in T^{\mathrm{reg}}_\C$, $\theta_w(\chi)=0$ and by Lemma~\ref{lemDensity_regular_elements}, $\theta_w=0=Q\tilde{\theta}_w$. Let $w\in W^v$ be such that $\theta_w\neq 0$. Then $U:=\{\chi\in T_\C|\theta_w\in \ATC_\chi\text{ and }\theta_w(\chi)\neq  0\}$ is open and dense in $T_\C$. By Remark~\ref{rkMeasure_irreducible_representations}, $T_\C^{\mathrm{reg}}$ has full measure in $T_\C$ and thus $U\cap T_\C^{\mathrm{reg}}$ is dense in $T_\C$. Moreover $\theta_w(\chi)=Q(\chi)\tilde{\theta}(\chi)$ for all $\chi\in U\cap T_\C^{\mathrm{reg}}$ and thus $\tilde{\theta}_w=Q\theta_w$.   Consequently, there exists $Q\in \C(Y)$ such that $F_u*F_v= F_{uv}*Q$.

Let $\tau\in \UC_\C$. Let $w_1\in W^v$. Let $u\in W^v$ be such that there exists $\theta=\frac{f}{g}\in \C(Y)^\times$ such that  $F_{w_1}*F_u=F_{w_1u}*\theta$, with $f(w.\tau)\neq 0$ for all $w\in W^v$. Let $s\in \SCC$ be such that $us>u$.  Then by Lemma~\ref{lemFs2},
 \[F_{w_1}*F_{us}=F_{w_1u}*\theta*F_s=F_{w_1u}*F_s*\theta^s.\]
  Suppose $w_1us>w_1u$. Then $F_{w_1u}*F_s=F_{w_1us}$ and thus $F_{w_1}*F_{us}=F_{w_1us}*\theta^s$ and $f^s(w.\tau)\neq 0$ for all $w\in W^v$. Suppose $w_1us<w_1u$. Then $F_{w_1u}*F_s=F_{w_1us}*(F_s)^2$ and thus by Lemma~\ref{lemFs2}, $F_{w_1}*F_{us}=F_{w_1us}*(\theta^s \zeta_s\zeta_s^s)$. By definition of $\UC_\C$,  one can write $F_{w_1}*F_{us}=F_{w_1us}*\frac{\tilde{f}}{\tilde{g}}$ with $\tilde{f},\tilde{g}\in \C[Y]^\times$ such that $\tilde{f}(w.\tau)\neq 0$ for all $w\in W^v$ and the lemma follows.
\end{proof}

\begin{remark}
In \cite[Lemma~4.3 (2)]{reeder1997nonstandard}, Reeder gives an explicit expression of $F_u*F_v$, for $u,v\in W^v$.
\end{remark}

Let $r\in \RCC$. Let $\Omega_r:\C(Y)\rightarrow \C(Y)$\index{$\Omega_r$} be defined by $\Omega_r(\theta)=\zeta_r(\theta^r-\theta)$ for all $\theta\in \C(Y)$.

\begin{lemma}\label{lemStability_C(Y)tau}
Let $r\in \SCC_\tau$. Then $\Omega_r(\C(Y)_\tau)\subset \C(Y)_\tau$.
\end{lemma}

\begin{proof}
 Write $r=r_{\beta^\vee}$, where $\beta^\vee\in \Phi^\vee$. Then  one has $r(\lambda)=\lambda-\beta(\lambda)\beta^\vee$ for all $\lambda\in Y$. Let $\lambda\in Y$. Then  with the same computation as in Remark~\ref{remIH algebre dans le cas KM deploye}~(\ref{itPolynomiality_Bernstein_Lusztig}), we have that $\Omega_r(Z^\lambda)\in \C(Y)_\tau$. Thus $\Omega_r(\theta)\in \C(Y)_\tau$ for all $\theta\in \C[Y]$.

Let $\theta\in \C(Y)_\tau$. Write $\theta=\frac{f}{g}$, where $f,g\in \C[Y]$ and $g(\tau)\neq 0$. Then $\zeta_r (\theta^r-\theta)=\zeta_r(\frac{f^rg-(f^rg)^r}{gg^r})$. Moreover, $g^r(\tau)=g(r.\tau)=g(\tau)\neq 0$ and as $f^rg\in \C[Y]$, we have that $\zeta_r(\theta^r-\theta)\in \C(Y)_\tau$. 
\end{proof}

We now assume that $\tau\in \UC_\C$.

For each $w\in \Wta$ we fix a reduced writing $w=r_1\ldots r_k$, with $k=\ell(w)$ and $r_1,\ldots,r_k\in \SCC_\tau$ and we set $\underline{w}=(r_1,\ldots,r_k)$. Let $K_{\underline{w}}=K_{r_1}\ldots K_{r_k}\in \ATC$\index{$K_{\underline{w}}$}.

\begin{lemma}\label{lemStability_HCYtau_K}
Let $r\in \SCC_\tau$. Then $\ATC_\tau*K_r\subset \ATC_\tau$. In particular, $K_{\underline{w}}\in \ATC_\tau$ for all $w\in \Wta$. 
\end{lemma}

\begin{proof}
Let $w\in W^v$ and $\theta\in \C(Y)_\tau$. Then $H_w\theta*K_r=H_wK_r \theta^r+H_w*\Omega_r(\theta)$. Using Lemma~\ref{lemReeder14.3}, Lemma~\ref{lemStability_C(Y)tau} and the fact that $\ATC_\tau$ is a $\HCW-\C(Y)_\tau$-bimodule, we deduce that $H_w\theta*K_r\in \ATC_\tau$. Hence $\ATC_\tau*K_r\subset \ATC_\tau$.
\end{proof}

\begin{lemma}\label{lemMaximal_support_K}
Let $w\in \Wta$. Then $\max\supp\big(K_{\underline{w}}(\tau)\big)=\{w\}$, where $\max$ is defined with respect to the order $\leq$ on $W^v$.
\end{lemma}

\begin{proof}
Write  $\underline{w}=(r_{1},\ldots ,r_{k})$ with $r_1,\ldots, r_k\in \SCC_\tau$. Then \[K_{\underline{w}}=(F_{r_{i_1}}-\zeta_{r_{i_1}})\ldots (F_{r_{i_k}}-\zeta_{r_{i_k}})=F_{r_{i_1}}*F_{r_{i_2}}*\ldots*F_{r_{i_k}}+\sum_{v<_\tau w} F_v P_v,\]  for some $P_v\in \C(Y)$. 
By Lemma~\ref{lemProductFuFv}, there exist $f,g\in \C[Y]^\times$ such that  $F_{r_{i_1}}*F_{r_{i_2}}*\ldots*F_{r_{i_k}}= F_{w}*\frac{f}{g}$ and $f(\tau)\neq 0$. One has $\pi^T_w(F_w)=1$ and by Lemma~\ref{lemComparison_Bruhat_orders}, $\pi_v^T(F_v)=0$ for all $v\in [1,w)_{\leq_\tau}$. Thus using Lemma~\ref{lemStability_HCYtau_K}, one can moreover assume $g(\tau)\neq 0$. Therefore $\pi^T_w(K_{\underline{w}})=\frac{f}{g}\in \C(Y)_\tau$ and $f(\tau)\neq 0$, which proves the lemma.
\end{proof}

Let $\KC(\Wta)=\bigoplus_{w\in \Wta} F_w \C(Y)$\index{$\KC(\Wta)$}. By Lemma~\ref{lemProductFuFv} and Lemma~\ref{lem4.3 de reeder}, $\KC(\Wta)$ is a subalgebra of $\ATC$. Let $\KC_\tau=\KC(\Wta)\cap \ATC_\tau$\index{$\KC_\tau$}. For $w\in \Wta$, set $\KC(\Wta)^{<_\tau w}=\bigoplus_{v\in \Wta, v<_\tau w}F_w\C(Y)$ and $\KC_\tau^{<_\tau w}=\bigoplus_{v<_\tau w} K_{\underline{v}}\C(Y)_\tau$.

\begin{lemma}\label{lemCommutation_relation_K_theta}
Let $\theta\in \C(Y)_\tau$ and $w\in \Wta$. Then there exists $k_{\underline{w}}(\theta)\in \KC_\tau^{<_\tau w}$  such that $\theta*K_{\underline{w}}=K_{\underline{w}}*\theta^{w^{-1}}+k_{\underline{w}}(\theta)$.
\end{lemma}

\begin{proof}
If $w=1$, this is clear. Suppose $w>_\tau 1$. Write $w=vr$ with $v\in \Wta$ and $r\in \SCC_\tau$ such that $v<_\tau w$. Suppose that $\theta *K_{\underline{v}}=K_{\underline{v}}*\theta^{v^{-1}}+k_{\underline{v}}(\theta)$ with $k_{\underline{v}}(\theta)\in \KC_\tau^{<_\tau v}$. One has 
\[ \theta* K_{\underline{w}}=\theta* K_{\underline{v}}*K_r=\big(K_{\underline{v}} \theta^{v^{-1}}+k_{\underline{v}}(\theta)\big)*K_r=K_{\underline{w}}*\theta^{w^{-1}}+K_{\underline{v}}*\Omega_r(\theta^{v^{-1}})+k_{\underline{v}}
(\theta)*K_r.\]

 The sets $\KC( \Wta)^{\leq_\tau v}=\bigoplus_{v'\leq_\tau v} F_{v'}\C(Y)$ and  $\ATC_\tau
 $ are right $\C(Y)_\tau$-submodules of $\ATC$ and thus  by Lemma~\ref{lemStability_HCYtau_K} and Lemma~\ref{lemStability_C(Y)tau}, $K_{\underline{v}}*\Omega_r(\theta^{v^{-1}})\in \KC_\tau^{\leq_\tau v}\subset \KC_\tau^{<_\tau w}$.
 
By Lemma~\ref{lemStability_HCYtau_K}, $k_{\underline{v}}(\theta)*K_r\in \ATC_\tau$. By Lemma~\ref{lem4.3 de reeder} and \cite[Corollary 1.3.19]{kumar2002kac}, $k_{\underline{v}}F_r\in \KC(\Wta)^{<_\tau\max (vr,v)}=\KC(\Wta)^{<_\tau w}$. Consequently $k_{\underline{v}}*K_r\in \KC_\tau^{<_\tau w}$ and  $K_{\underline{v}}\Omega_r(\theta^{v^{-1}})+k_{\underline{v}}(\theta)K_r\in \KC_\tau^{<_\tau w}$, which proves the lemma.
\end{proof}

\begin{lemma}\label{lemDescription_Ktau}
One has $\KC_\tau=\bigoplus_{w\in \Wta} K_{\underline{w}}\C(Y)_\tau$.
\end{lemma}

\begin{proof}
By Lemma~\ref{lemStability_HCYtau_K}, $\KC_\tau\supset \bigoplus_{w\in \Wta} K_{\underline{w}}\C(Y)_\tau$.

For $w\in \Wta$, set $\KC(\Wta)^{\leq_\tau w}=\bigoplus_{v\leq_\tau w}F_v \C(Y)\subset \KC(\Wta)$. Let $w\in \Wta$. Suppose that for all $v\in [1,w)_{\leq_\tau}$, one has $\KC_\tau^{\leq_\tau v}= \bigoplus_{v'\in [1,v]_{\leq_\tau}} K_{\underline{v'}}\C(Y)_\tau$. By Lemma~\ref{lemMaximal_support_K}, one can write $\pi^T_w(K_{\underline{w}})=\frac{f}{g}$, with $f,g\in \C[Y]$ such that $f(\tau)g(\tau)\neq 0$. Let $x\in \KC_\tau^{\leq_\tau w}$ and $\theta=\pi^T_w(x)\in \C(Y)_\tau$. By Lemma~\ref{lemStability_HCYtau_K}, $\theta \frac{g}{f}K_{\underline{w}}\in \ATC_\tau$. Moreover,  $x-\theta \frac{g}{f}K_{\underline{w}}\in \sum_{v\in [1,w)_{\leq_\tau}} \KC_\tau^{\leq_\tau v}$. Therefore, $x\in \bigoplus_{v\in [1,w]_{\leq_\tau}}K_{\underline{v}}\C(Y)_\tau$ and the lemma follows.
\end{proof}

\begin{theorem}\label{thmGeneralized_weigth_space}
Let $\tau\in \UC_\C$ be such that $\Wta=W_{\tau}$. Then $I_\tau(\tau,\mathrm{gen})=\ev_\tau(\KC_\tau)\otimes_\tau 1$. 
\end{theorem}

 \begin{proof}
 Let $w\in \Wta$ and $\theta\in \C(Y)_\tau$. As $w\in W_\tau$, $\theta^{w^{-1}}\in \C(Y)_{w.\tau}=\C(Y)_\tau$ and $\tau(\theta^{w^{-1}})=\tau(\theta)$.  Then by Lemma~\ref{lemCommutation_relation_K_theta}, $(\theta-\tau(\theta))K_{\underline{w}}(\tau)\otimes_\tau 1\in \KC^{<_\tau w}(\tau)\otimes_\tau 1$. By an induction using Lemma~\ref{lemDescription_Ktau} we deduce that $\KC_{\tau}(\tau)\otimes_\tau 1 \subset I_\tau(\tau,\mathrm{gen})$.
 
 Let $w\in W^v$ and  $E_w=\big(\ev_\tau(\KC_\tau)\otimes_\tau 1\big)\cap I_\tau^{\leq w}$. By Lemma~\ref{lemMaximal_support_K}, $\dim E_w= |\Wta\cap \{v\in W^v|v\leq w\}|$. By Proposition~\ref{propDecomposition_generalized_weight_spaces_Itau}, $\dim I_\tau(\tau,\text{gen})^{\leq w}=|\{v\in W_\tau|v\leq w\}|=\dim E_w$. As $(W^v,\leq)$ is a directed poset, $I_\tau=\bigcup_{w\in W^v} I_\tau^{\leq w}$, which proves the theorem.
 \end{proof}

\subsection{Irreducibility of $I_\tau$ when $W_\tau=W_{(\tau)}$ is the infinite dihedral group}\label{subIrreducibility_Itau_infinite_dihedral_group}

In this subsection, we prove that if $\tau\in \UC_\C$ is such that $W_\tau=\Wta$  and $\Wta$ is isomorphic to  the infinite dihedral group, then $I_\tau$ is irreducible (see Lemma~\ref{lemIrreducibility of Itau}). Let us sketch the proof of this lemma. We prove that $I_\tau(\tau)=\C 1\otimes_\tau 1$. For $w\in \Wta$, let $\pi^K_w:\I_\tau(\tau,\mathrm{gen})\rightarrow \C$ be defined as $\pi^K_{w}\big(\sum_{v\in W^v} K_{\underline{v}}(\tau)x_v\big)=x_w$, for all $(x_v)\in \C^{(\Wta)}$, which is well-defined by Lemma~\ref{lemMaximal_support_K} and Theorem~\ref{thmGeneralized_weigth_space}. We suppose that  $I_\tau(\tau)\setminus\C 1 \otimes_\tau 1$ is nonempty and we consider one of its elements $x$.  We reach a contradiction by computing $\pi^K_w(x)$, where $w\in \Wta$ is such that $\ell_\tau(w)=\max \{\ell_\tau(v)|v\in \supp(x)\cap \Wta\}-1$.

Let $\tau\in \UC_\C$. Assume that $(\Wta,\SCC_\tau)$ is isomorphic to the infinite dihedral group (in particular, $|\SCC_\tau|=2$ and every element of $\Wta$ admits a unique reduced writing). 

\medskip

The following lemma is easy to prove.

\begin{lemma}\label{lemBruhat_order_K}
Let $w\in \Wta$ and $r\in \SCC_\tau$ be such that $\ell_\tau(wr)=\ell_\tau(w)+1$. Let $u\in [1,w)_{\leq \tau}$. Then $ur\neq w$.
\end{lemma}

\begin{lemma}\label{lemComputation_coefficient}
Let $\tau\in \UC_\C$. Let $r=r_{\beta^\vee}\in \SCC_\tau$, where $\beta^\vee\in \Phi^\vee$. Then there exists $a\in \C^*$ such that for all $\lambda\in Y$, \[\tau\big((Z^{r.\lambda}-Z^{\lambda})\zeta_r\big)=a\tau(\lambda)\beta(\lambda).\]
\end{lemma}

\begin{proof}
One has \[\zeta_r=\frac{1}{\zeta_{\beta^\vee}^{\mathrm{den}}}.\prod_{\alpha^\vee\in N_{\Phi^\vee}(r)}\zeta_{\alpha^\vee}^{\mathrm{num}}. \prod_{\alpha^\vee\in N_{\Phi^\vee}(r)\setminus\{\beta^\vee\}}\frac{1}{\zeta_{\alpha^\vee}^{\mathrm{den}}}.\] By Lemma~\ref{lemComparison_definition_tau_simples} and by definition of $\UC_\C$, \[\tau(\prod_{\alpha^\vee\in N_{\Phi^\vee}(r)\setminus\{\beta^\vee\}} \zeta_{\alpha^\vee}^{\mathrm{den}})\neq 0\text{ and }\tau(\prod_{\alpha^\vee\in N_{\Phi^\vee}(r)}\zeta_{\alpha^\vee}^{\mathrm{num}})\neq~0.\] 

If $\sigma_{\beta^\vee}=\sigma'_{\beta^\vee}$, one  has $\frac{Z^{r.\lambda}-Z^\lambda}{\zeta_{\beta^\vee}^{\mathrm{den}}}=\frac{Z^{r.\lambda}-Z^\lambda}{1-Z^{\beta^\vee}}=Z^\lambda \frac{Z^{-\beta(\lambda)\beta^\vee}-1}{1-Z^{\beta^\vee}}$. By Lemma~\ref{lemComparison_definition_tau_simples}, $r\in \RCC_{(\tau)}$ and thus $\tau(\beta^\vee)=1$.  Thus by the same computation as in Remark~\ref{remIH algebre dans le cas KM deploye}, $\tau(\frac{Z^{r.\lambda}-Z^\lambda}{1-Z^{\beta^\vee}})=\beta(\lambda)\tau(\lambda)$. Using a similar computation when $\sigma_{\beta^\vee}\neq \sigma'_{\beta^\vee}$, we deduce the lemma.
\end{proof}

\begin{lemma}\label{lemCoefficient in w, 1 K}
Let $w\in \Wta$ and  $r\in \SCC_\tau$ be such that $\ell_\tau(wr)=\ell_\tau(w)+1$. Then there exists $a\in \C^*$ such that for all  $\lambda\in Y$, one has: \[\pi^K_w\big(Z^\lambda*K_{\underline{wr}}(\tau)\otimes_{\tau} 1\big)=a\tau(\lambda)\alpha_r(w^{-1}.\lambda).\]
\end{lemma}

\begin{proof}
Let $\lambda\in Y$. Write $Z^\lambda* K_{\underline{w}}=K_{\underline{w}}*Z^{w^{-1}.\lambda}+k$, where $k\in \KC_\tau^{<_\tau w}$, which is possible by Lemma~\ref{lemCommutation_relation_K_theta}. One has \[Z^\lambda *K_{\underline{wr}}=(K_{\underline{w}}*Z^{w^{-1}.\lambda}+k)*K_r=K_{\underline{wr}}*Z^{rw^{-1}.\lambda}+K_{\underline{w}} \big((Z^{rw^{-1}.\lambda} -Z^{w^{-1}.\lambda})\zeta_r\big)+k*K_r.\]  
Therefore, using Lemma~\ref{lemBruhat_order_K} and Lemma~\ref{lemComputation_coefficient} we deduce \[\pi^K_w\big(Z^\lambda K_{\underline{wr}}(\tau)\otimes_{\tau}1\big)=\tau\big((Z^{rw^{-1}.\lambda}-Z^{w^{-1}.\lambda})\zeta_r\big)=a\tau(\lambda)\beta(w^{-1}.\lambda) ,\] for some $a\in \C^*$.
\end{proof}

\begin{lemma}\label{lemBruhat_order_K_2}
Let $w\in \Wta$ and $r\in \SCC_\tau$ be such that $\ell_\tau(rw)=\ell_\tau(w)+1$.

 One has $\pi^K_w(K_r*\KC(\Wta)^{\leq_\tau w})=\{0\}$.

\end{lemma}

\begin{proof}
Let $u\in \Wta$ and $r\in \SCC_\tau$ be such that $ru>_\tau u$. Then by Lemma~\ref{lemProductFuFv} and \cite[Corollary~1.3.19]{kumar2002kac}, $F_r*\KC(\Wta)^{\leq_\tau u}\subset \KC(\Wta)^{\leq_\tau \max(u,ru)}$ and thus $K_r*\KC(\Wta)^{\leq_\tau u}\subset \KC(\Wta)^{\leq_\tau \max(u,ru)}$.

 Let $v\in [1,w)_{\leq_\tau}$. If $rv>_\tau v$, then by Lemma~\ref{lemProductFuFv}, there exists $Q\in \C(Y)$ such that $F_r*F_v=F_{rv}*Q$ and thus $K_r*F_v\in F_{rv}*Q+F_v\C(Y)$. By Lemma~\ref{lemBruhat_order_K}, $rv\neq w$. Using Lemma~\ref{lemMaximal_support_K} and the fact $w$ and $rv$  have the same length, we deduce that $\pi^K_w(K_r*F_v)=0$. 
 
 If  $rv<_\tau v$, then $K_r*F_v\in \KC(\Wta)^{\leq_\tau v}$  and thus $\pi^K_w(K_r*F_v)=0$ which finishes the proof of the lemma.
\end{proof}

\begin{lemma}\label{lemCoefficient in w,2_K}
Let $w\in W_\tau$, $r\in \SCC_\tau$ be such that $\ell_\tau(rw)=\ell_\tau(w)+1$. Then there exists $b\in \C^*$ such that for all $\lambda\in Y$: \[\pi^K_w(Z^\lambda.K_{\underline{rw}}(\tau)\otimes_{\tau} 1)=b\tau(\lambda)\alpha_r(\lambda).\]
\end{lemma}

\begin{proof}
One has \[Z^\lambda K_{\underline{rw}}=(Z^\lambda*K_r)*K_{\underline{w}}=\big(K_r. Z^{r.\lambda}+(Z^{r.\lambda}-Z^\lambda)\zeta_r\big)*K_{\underline{w}}(\tau).\] 
One has   $Z^{r.\lambda}*K_{\underline{w}}\in \KC(\Wta)^{\leq_\tau w}$. Thus by Lemma~\ref{lemBruhat_order_K_2}, $\pi^K_w(K_r.Z^{r.\lambda}*K_{\underline{w}})=0$. Moreover, by Lemma~\ref{lemComputation_coefficient}, there exists $b\in \C^*$ such that  \[\pi^{K}_w\big((Z^{r.\lambda}-Z^\lambda)\zeta_r K_{\underline{w}}(\tau)\otimes_\tau 1\big)=w.\tau\big((Z^{r.\lambda}-Z^\lambda)\zeta_r\big) =b\tau(\lambda)\alpha_r(\lambda),\] which proves the lemma.
\end{proof}

\begin{lemma}\label{lemIrreducibility of Itau}
Let $\tau\in \UC_\C$ be such that $W_\tau=\Wta$ and such  that there exists $r_1,r_2\in \SCC_\tau$ such that  $(\Wta,\{r_1,r_2\})$ is isomorphic to the infinite dihedral group. Then $I_\tau$ is irreducible.
\end{lemma}

\begin{proof}
Let us prove that $I_\tau(\tau)=\C.1\otimes_{\tau}1$. Let $x\in I_\tau\setminus \C. 1\otimes_{\tau} 1$ and assume that $x\in  I_\tau(\tau)$. Let $n=\max \{\ell_\tau(w)| w\in\supp(x)\}$. Let $w\in \Wta$ be such that $\ell_\tau(w)=n-1$. Then there exist $r,r'\in \SCC_\tau$ such that $\{v\in \Wta|\ell_\tau(v)=n\}=\{rw,wr'\}$. By Theorem~\ref{thmGeneralized_weigth_space}, $x\in \sum_{v\in \Wta} \C K_{\underline{v}}(\tau)\otimes_\tau 1$. Let $\gamma=\pi^K_{rw}(x)$ and $\gamma'=\pi^K_{wr'}(x)$.

Set $\gamma_w=\pi^K_w(x)$. Then by Lemma~\ref{lemCoefficient in w, 1 K} and Lemma~\ref{lemCoefficient in w,2_K}, there exist $a,a'\in \C^*$ such that for all $\lambda\in Y$, \[\pi^K_w(Z^\lambda.x)=\tau(\lambda)\big(a\gamma \alpha_{r}(\lambda)+a'\gamma'w.\alpha_{r'}(\lambda)+\gamma_w\big)=\tau(\lambda)\gamma_w.\] Therefore $\{\alpha_{r}, w.\alpha_{r'}\}$ is linearly dependent and hence $w.\alpha_{r'}\in \{\pm \alpha_r\}=\{\alpha_r,r.\alpha_r\}$. By Lemma~\ref{lemKumar1.3.11} we deduce $rw=wr'$: a contradiction because $|\{rw,wr'\}|=|\{v\in \Wta|\ell_\tau(v)=n\}|=2$.

Therefore $I_\tau=\C 1\otimes_\tau 1$ and by Theorem~\ref{thmIrreducibility criterion}, $I_\tau$ is irreducible.
\end{proof}

\subsection{Kato's criterion when the Kac-Moody matrix has size $2$}\label{subKatos_irreducibility_criterion}

In this subsection, we prove Kato's irreducibility criterion  when  $|\SCC|=2$ (see Theorem~\ref{thmKato's_theorem_dim2}). As the case where $W^v$ is finite is a particular case of Kato's theorem \cite[Theorem~2.2]{kato1982irreducibility} we assume that $W^v$ is infinite.

This is equivalent to assuming that the Kac-Moody matrix of the root generating system $\mathcal{S}$ is of the form $\begin{pmatrix}
2 & a\\ b& 2
\end{pmatrix}$, with $a,b\in \Z_{<0}$ and $ab\geq 4$ (\cite[Proposition 1.3.21]{kumar2002kac}). The system  $(W^v,\SCC)$ is then the infinite dihedral group. Write $\SCC=\{s_1,s_2\}$. Then every element of $W^v$ admits a unique reduced writing involving $s_1$ and $s_2$.

Let $G$ be a group and $a,b\in G$. For $k\in \N$, we define $P_k(a,b)=aba\ldots$ where the products has $k$ terms. 

\begin{lemma}\label{lemClassification_subgroups_W_dim2}
The subgroups of $W^v$ are exactly the ones of the following list: \begin{enumerate}
\item $\{1\}$

\item $\langle r\rangle=\{1,r\}$, for some $r\in \RCC$

\item $Z_k=\langle P_{2k}(s_1,s_2)\rangle=\langle P_{2k}(s_2,s_1)\rangle\simeq  \Z$ for $k\in \Ne$

\item $R_{k,m}=\langle P_{2k+1}(s_1,s_2),P_{2m+1}(s_2,s_1)\rangle\simeq W^v$ for $k,m\in \N$.
\end{enumerate}
\end{lemma}

\begin{proof}
Let $\{1\}\neq H\subset W^v$ be a subgroup. Let $n=\min \{\ell(w)|w\in H\setminus \{1\}\}.$ 

First assume that $n$ is even and set $k=\frac{n}{2}$. Then $P(s_1,s_2,n)=P(s_2,s_1,n)^{-1}$ and as these are the only elements having length $n$ in $W^v$, $H\supset Z_k$. Let $w=P_n(s_1,s_2)$. Let $h\in H\setminus \{1\}$. Write $\ell(h)=an+r$ with $a\in \Ne$ and $r\in \llbracket 0,r-1\rrbracket$. Then there exists $\epsilon\in \{-1,1\}$ such that $h=w^{\epsilon a}.h'$, with $\ell(h')=r$. Moreover, $h'\in H$ and thus $h'=1$. Therefore $H=Z_k$.

We now assume that $n$ is odd. Maybe considering $vHv^{-1}$ for some $v\in W^v$ and exchanging the roles of $s_1$ and $s_2$, we may assume that $s_1\in H$. Assume $H\neq \langle s_1\rangle$. Let $n'=\min \{\ell(w)|w\in H\setminus\langle s_1\rangle\}$. Let $w\in H\setminus \langle s_1\rangle$ be such that $\ell(w)=n'$. Then the reduced writing of $w$ begins and ends with $s_2$. Thus $n'=2n''+1$ for some $n''\in \N$. Then it is easy to see that $H=R_{1,n''}$, which finishes the proof.
\end{proof}

We prove in Appendix~\ref{secExample_possibilities_Wtau} that there exists size $2$ Kac-Moody matrices such that for each subgroup of $W^v$, there exists $\tau\in T_\C$ such that $\Wta$ is isomorphic to this subgroup.

\begin{theorem}\label{thmKato's_theorem_dim2}
Assume that the matrix of the root generating system $\mathcal{S}$ is of size $2$. Let $\tau\in T_\C$. Then $I_\tau$ is irreducible if and only if $\tau\in \UC_\C$ and $W_\tau=\Wta$.
\end{theorem}

\begin{proof}
If $W^v$ is finite, this is a particular case of Kato's theorem (\cite[Theorem~2.2]{kato1982irreducibility}). Suppose that $W^v$ is infinite. By Lemma~\ref{lemIrreducibility implies isomorphisms} and Proposition~\ref{propKato's weak theorem}, if $I_\tau$ is irreducible, then $\tau\in \UC_\C$ and $W_\tau=\Wta$. Reciprocally, suppose $\tau\in \UC_\C$ and $W_\tau=\Wta$. Then by Lemma~\ref{lemClassification_subgroups_W_dim2}, either $\Wta=\{1\}$, or $\Wta=\langle r\rangle$ for some $r\in \RCC$ or $\Wta=\langle r_1,r_2\rangle$ for some $r_1,r_2\in \RCC$ and $(\Wta,\{r_1,r_2\})$ is  isomorphic to the infinite dihedral group. In the first two cases, $I_\tau$ is irreducible by Corollary~\ref{corMatsumoto theorem} or Corollary~\ref{corIrreducibility when Wchi=s}. Suppose $\Wta=\langle r_1,r_2\rangle$. Then by Remark~\ref{rkCoxeter_generators}~(\ref{itCoxeter_generators_infinite_dihedral_group}), $(\Wta,\SCC_\tau)$ is isomorphic to the infinite dihedral group and $I_\tau$ is irreducible by Lemma~\ref{lemIrreducibility of Itau}.
\end{proof}

\paragraph{Comments on the proofs of Kato's criterion}

There are several proofs of Kato's criterion in the literature. In \cite{reeder1991certain}, Reeder proves this criterion (see Corollary 8.7). In his proof, he uses the $R$-group $R_\tau=\{w\in W_\tau| w(\Phi^\vee_{(\tau)}\cap \Phi^\vee_+)=\Phi^\vee_{(\tau)}\cap \Phi^\vee_+\}$.  This group is reduced to $\{1\}$ when $W_\tau=\Wta$. His proof uses Harish-Chandra completeness theorem, which - under certain hypothesis on $\tau$ - majorizes the dimension of the space of  intertwining operators of $I_\tau$. Unfortunately, it seems that there exists up to now no equivalent of Harish-Chandra completeness theorem available in the Kac-Moody framework.

In \cite{rogawski1985modules}, Rogawski gives a proof of a particular case of Kato's criterion (see Corollary 3.2). However, it seems that its proof uses the fact that every element $x$ of $I_\tau(\tau)$ can be written as a sum $x=\sum_{j\in J} x_j$ where $J$ is a finite set and for all $j\in J$, $|\max \supp(x_j)|=1$ and $x_j\in I_\tau(\tau)$. I do not know how to prove such a property.

In \cite{reeder1997nonstandard}, Reeder gives two proofs of Kato's criterion or of weak versions of it (see Corollary~4.6 and Theorem~14.7).  Our proof of Theorem~\ref{thmKato's_theorem_dim2} is strongly inspired by the proof of \cite[Theorem~14.7]{reeder1997nonstandard}.

\section{Towards principal series representations of $G$}\label{secTowards_principal_series_representations}

Suppose that $\HC_\C$ is associated with a reductive group $G$. Then  for every open compact subgroup $K'$ of $G$ and  every smooth representation $V$, $V^{K'}$ is naturally equipped with the structure of an $\HC_{K',\C}$ module, where $\HC_{K',\C}$ is the Hecke algebra associated with $K'$ with coefficients in $\C$. Moreover, the assignment $V\mapsto V^{K'}$  induces a bijection between the following sets: \begin{itemize}
\item equivalence classes of irreducible smooth representations $V$ of $G$
such that $V^{K'} \neq \{0\}$,

\item  isomorphism classes of simple $\HC_{K',\C}$-modules (see \cite[4.3]{bushnell2006local} for example).
\end{itemize}

In the Kac-Moody case, we do not know how to define ``smooth'' for a representation of $G$. We know that for any topological group structure on $G$,   $K_I$   is not compact open (see \cite[Theorem 3.1]{abdellatif2019completed}). The hope is that there should be a link between  representations of $G$ satisfying some regularity conditions and representations of $\HC_\C$ or $\AC_\C$.

Let $\epsilon\in \{+,\emptyset\}$. In this section, we associate to every $\tau\in T_\FC^\epsilon$ a representation $\widehat{I(\tau^\epsilon)^\epsilon}$ of $G^\epsilon$. The principal series representation associated with $\tau$ should  correspond to the space of elements of $\widehat{I(\tau^\epsilon)^\epsilon}$ which satisfy some regularity condition.  We define an action of $\HC_\FC$ on some subspace   $I_{\tau^\epsilon,G^\epsilon}$ of $\big(\widehat{I(\tau^\epsilon)^\epsilon}\big)^{K_I}$. We then prove that $I_{\tau^\epsilon,G^\epsilon}$ is isomorphic (as an $\HC_\FC$-module) to the representation $I_{\tau^\epsilon|G^+}^+$ introduced in section~\ref{secIH algebras}. We then study the extendability of $\widehat{I(\tau^\epsilon)^\epsilon}$ and $I_{\tau^\epsilon,G^\epsilon}$ to representations of $G$ and $\AC_\FC$.

For simplicity, we only introduce split Kac-Moody groups, although our results also apply to almost-split Kac-Moody groups over local fields, see \cite{rousseau2017almost}.

\medskip

In subsection~\ref{subKM_groups}, we introduce split Kac-Moody groups over local fields, masures, their Iwahori-Hecke algebras and principal series representations.

In subsection~\ref{subAction_IH_algebra} we prove that the actions of $\HC_\FC$ on $I_{\tau,G^+}$ and $I_{\tau,G}$  are well-defined and prove that $I_{\tau,G^+}$ is isomorphic to $I_\tau$.

In subsection~\ref{subExtendability} we study under which condition $I_{\tau,G^+}$ and $I_\tau^+$ extend to representations of $G$ and  of $\AC_\FC$, for $\tau\in T_{\FC}^+$. We give examples of $\tau\in T_\FC$ (for particular choices of $G$)  such that $I_{\tau,G^+}$ and $I_\tau^+$ do not extend to representations of $G$ and of $\AC_\FC$.

\subsection{Kac-Moody groups over local fields and masures}\label{subKM_groups}

\subsubsection{Split Kac-Moody groups over local fields and masure}

 Let  $\mathbf{G}_{\mathcal{S}}$\index{$G$} be the group functor associated in \cite{tits1987uniqueness} with  the  generating root datum $\mathcal{S}$, see also \cite[8]{remy2002groupes}. Let $(\KC,\omega)$\index{$\KC$}\index{$\omega$} be a non-Archimedean local  field where $\omega:\KC\twoheadrightarrow \Z\cup\{+\infty\}$ is a valuation. Let $G=\mathbf{G}_{\mathcal{S}}(\KC)$ be the \textbf{split Kac-Moody group over $\KC$ associated with $\mathcal{S}$}.  The group $G$ is generated by the following subgroups:\begin{itemize}
\item the fundamental torus $T=\mathbf{T}(\KC)$\index{$T$}, where $\textbf{T}=\mathrm{Spec}(\Z[X])$,

\item the root subgroups $U_\alpha=\textbf{U}_\alpha(\KC)$\index{$U_\alpha$}, each isomorphic to $(\KC,+)$ by an isomorphism $x_\alpha$\index{$x_\alpha$}.
\end{itemize}

In \cite{gaussent2008kac} and \cite{rousseau2016groupes} (see also \cite{rousseau2017almost}) the authors associate a masure $\I$\index{$\I$} on which the group $G$ acts. We recall briefly the construction of this masure. Let $N$\index{$N$} be the normalizer of $T$ in $G$. Then they define an action of $N$ on $\A$, see \cite[3.1]{gaussent2008kac}. For $n\in N$ denote by $\nu(n):\A\rightarrow \A$\index{$\nu$} the affine automorphism of $\A$ induced by the action of $N$ on $\A$. Then $\nu(t)$ is a translation, for every $t\in T$ and $\nu(N)=W^v\ltimes Y$. For every $\wb\in W^v\ltimes Y$, we choose $n_\wb\in N$\index{$n_\wb$, $n_\lambda$} such that $\nu(n_\wb)=\wb$.

The masure  $\I$ is defined to be the set $G\times \A/\sim$, for some equivalence relation $\sim$ (see \cite[Definition 3.15]{gaussent2008kac}). Then $G$ acts on $\I$ by $g.[h,x]=[gh,x]$ for $g,h\in G$ and $x\in \A$, where $[h,x]$ denotes the class of $(h,x)$ for $\sim$. The map $x\mapsto [1,x]$ is an embedding of $\A$ in $\I$ and we identify $\A$ with its image. Then $N$ is the stabilizer of $\A$ in $G$  and it acts on $\A$ by $\nu$. If $\alpha\in \Phi$ and $a\in \KC$, then $x_\alpha(a)\in U_\alpha$ fixes the half-apartment $D_{\alpha,\omega(a)}=\{y\in \A|\ \alpha(y)+\omega(a)\geq 0\}$ and for all $y\in \A\setminus D_{\alpha,\omega(a)}$, $x_\alpha(a).y\notin \A$.

 An \textbf{apartment} is a set of the form $g.\A$, for $g\in G$. We have $\I=\bigcup_{g\in G} g.\A$.  Then $\I$ satisfies axioms (MA i), (MA ii) and (MA iii) of \cite[Appendix A]{hebert2018study} or \cite{hebert2020new}.These axioms describe the following properties. \begin{itemize}
 \item[(MA i)] Let $A$ be an apartment of $\I$. Then $A=g.\A$, for some $g\in G$. We can then transport every notion which is preserved by $\nu(N)=W^v\ltimes Y$ to $A$ (in particular, we can define a segment, a hyperplane, ... in  $A$).

 \item[(MA ii)] This axiom asserts that if $A$ and $A'$ are two apartments such that $A\cap A'$ is ``large enough'', then $A\cap A'$ is a finite intersection of half-apartments (i.e of sets of the form $h.D_{\alpha,k}$, for $\alpha\in \Phi$, $k\in \Z$, if $A=h.\A$) and there exists $g\in G$ such that $A'=g.A$ and $g$ fixes $A\cap A'$. When $G$ is an affine Kac-Moody group, this is true for every pair of apartments $A,A'$, without any assumption on $A\cap A'$.

  \item[(MA iii)] This axiom asserts that for some pairs of filters on $\I$, there exists an apartment containing them. This axiom is the building theoretic translation of some decompositions of $G$ (e.g Iwasawa decomposition).
 \end{itemize}

A \textbf{filter} on a set $E$ is a nonempty set $\VC$ of nonempty subsets of $E$ such that, for all subsets $S$, $S'$ of $E$,  if $S$, $S'\in \VC$ then $S\cap S'\in \VC$ and, if $S'\subset S$, with $S'\in \VC$ then $S\in \VC$.

Let $E,E'$ be  sets, $E'\subset E$ and $\VC$ be a filter on $E'$. One says that a set $\Omega\subset E$ contains $\VC$ if there exists $\Omega'\in \VC$ such that $\Omega'\subset \Omega$ (or equivalently if $\Omega\in \VC$ if $E=E'$). Let $f:E\rightarrow E$. One says that $f$ fixes $\VC$ if there exists $\Omega'\in \VC$ such that $f$ fixes $\Omega'$.

\subsubsection{Cartan decomposition, Tits preorder on $\I$ and sub-semi-group $G^+$}
Let $K=\mathbf{G}_\SC(\mathcal{O})$\index{$K$}, where $\OC$\index{$\OC$} is the ring of integers of $\KC$. Then $K$ is the fixator of $0\in \A\subset \I$ in $G$. For $\lambda\in Y$, choose $n_\lambda\in T$ such that $n_\lambda$ induces the translation on $\A$ by the vector $\lambda$.
Unless $G$ is reductive, the Cartan decomposition of $G$ does not hold: $\bigsqcup_{\lambda\in Y^{++}} Kn_\lambda K\subsetneq G$, where $Y^{++}=\overline{C^v_f}\cap Y$. For $x,y\in \A$, one writes $x\leq y$\index{$\leq$} if $y-x\in \T$ (where $\T$ is the Tits cone). If $x,y\in \I$, one writes $x\leq y$ if there exists $g\in G$ such that $g.x,g.y\in \A$ and $g.x\leq g.y$. This defines a $G$-invariant  preorder on $\I$ by \cite[Théorème 5.9]{rousseau2011masures}. We call it the \textbf{Tits preorder on $\I$}. Let $G^+=\{g\in G| g.0\geq 0\}$\index{$G^+$} (see \cite[1.2.2]{braverman2016iwahori} for a more explicit description of $G^+$, when $G$ is affine).  Then $G^+$ is a sub-semi-group of $G$ (as $\leq$ is transitive) and we have $G^+=\bigsqcup_{\lambda\in Y^{++}} K n_\lambda K$: the Cartan decomposition holds on $G^+$. Note that when $G$ is reductive, $G=G^+$ since $\T=\A$. A \textbf{type $0$ vertex} is a point of the form $g.0$ for some $g\in G$. We set $\I_0=G.0$\index{$\I_0$}. Then the map $g\mapsto g.0$ induces a bijection between $G/K$ and $\I_0$.

Let $x,y\in \I$ be such that $x\leq y$. Let $A_1,A_2$ be  apartments containing $x$ and $y$. Let $[x,y]_{A_1}$ (resp. $[x,y]_{A_2}$) be the segment in $A_1$ (resp. $A_2$) joining $x$ to $y$.  Then by \cite[Proposition 5.4]{rousseau2011masures}, $[x,y]_{A_1}=[x,y]_{A_2}$ and  there exists $g\in G$ such that $g.A_1=A_2$ and $g$ fixes $[x,y]_{A_1}$. We thus simply write $[x,y]$. Let $h \in G$ be such that $h.A_1=\A$. Then as $\leq$ is $G$-invariant, $h.x\leq h.y$ and thus $h.y-h.x\in \T$. Replacing $h$ by $nh$ for some $n \in N$, we may assume that $h.y-h.x\in \overline{C^v_f}$. One sets $\dv(x,y)=h.y-h.x\in \overline{C^v_f}$. We thus get a $G$-invariant \textbf{vectorial distance} $\dv:\I\times_{\leq}\I\rightarrow \overline{C^v_f}$, where $\I\times_{\leq}\I$ is the set of pairs $x,y\in \I$ such that $x\leq y$.  It is denoted $d^v$ in \cite{gaussent2014spherical}. When moreover $x,y\in \I_0$, then $\dv(x,y)\in Y^{++}$\index{$\dv$}. This distance parametrizes the $K$ double cosets: if $g\in G^+$ and $\lambda\in Y^+$, then $g\in Kn_\lambda K$ if and only if $\dv(0,g.0)=\lambda$.

\subsubsection{Local faces and chambers}

Recall the definition of vectorial faces from subsection~\ref{subRootGenSyst}.  A \textbf{local face of $\A$} (we omit the adjective ``local'' in the sequel) is a filter on $\A$ associated with a point $x$ and with a vectorial face $F^v$.  The point $x$ is the \textbf{vertex} of $F$ and $F^v$ is its direction. More precisely the chamber $F=F_{x,F^v}$\index{$F_{x,F^v}$} associated to $x$ and $F^v$ is the filter on $\A$ consisting of the sets   $\Omega\cap (x+F^v)$, where $\Omega$ is a neighborhood of $x$ in $\A$. We call $F$ \textbf{positive} (resp. \textbf{negative}) if $F^v$ is. When $F^v$ is a vectorial chamber (resp. a vectorial panel, that is when  $F^v$ is a codimension one face of a vectorial chamber), we call $F$ a \textbf{chamber} (resp. \textbf{panel}).  As the sets of local faces, of positive faces, of local chambers, ... are stable under the action of $W^v\ltimes Y$, we extend these notions to $\I$: a local face $F$ (resp. positive, negative) is a filter on $\I$ generated by $g.F$ for some local face (resp. positive, negative) $F_0$ and some $g\in G$. Its \textbf{vertex} is $\ve(F)=g.\lambda$\index{$\ve$}, where $\lambda$ is the vertex of $F_0$. This does not depend on the choices of $g$ and $F_0$ such that $F=g.F_0$. 

We denote by $C_0^+$\index{$C_0^+$} the local positive chamber associated with $0$ and $C^v_f$. A type $0$ positive local chamber is a filter of the form $g.C_0^+$ for some $g\in G$. Equivalently, this is a positive  chamber based at a type $0$ vertex. We denote by $\mathscr{C}_0^+$\index{$\CCC_0^+$} the set of positive type $0$   chambers of $\I$.

 We say that a chamber $C$ of $\A$  \textbf{dominates} a panel $P$ of $\A$ if $C$ and $P$ are based at the same vertex and if $P^v\subset \overline{C^v}$, where $C^v$ and $P^v$ are the vectorial faces defining $C$ and $P$.

  We say that a chamber $C$ of $\I$ \textbf{dominates} a panel $P$ of $\I$ if there exists $g\in G$ such that $g. C, g.P\subset \A$ and such that  $g.C$ dominates $g.P$. Then every type $0$ local panel is dominated by exactly $q+1$ chambers, where $q$ is the cardinal of the residue cardinal of $\KC$. In particular, $\I$ has \textbf{finite thickness}: every panel is dominated by finitely many chambers. This property is crucial in order to apply the finiteness results of \cite{gaussent2014spherical} and \cite{bardy2016iwahori}.

Let $W^+=W^v\ltimes Y^+$\index{$W^+$}. Then $W^+$ is a sub-semi-group of $W^v\ltimes Y$.If $C,C'\in \CCC_0^+$, we write $C\leq C'$ if $\ve(C)\leq \ve(C')$. Let $\CCC_0^+\times_{\leq}\CCC_0^+=\{(C,C')\in \CCC_0^+|C\leq C'\}$. Let $(C,C')\in \CCC_0^+\times_{\leq} \CCC_0^+$. Then by \cite[Proposition 5.5]{rousseau2011masures} or \cite[Proposition 5.17]{hebert2020new}, there exists an apartment $A=g.\A$ containing $C$ and $C'$. Then $g.C\subset \A$ and thus there exists $\wb\in W^v\ltimes Y$ such that $g.C=\wb.C_0^+$. Maybe replacing $g$ by $n_\wb^{-1} g$, we may assume that $g.C=C_0^+$. Then $g.C'\geq C$ and thus there exists $\vb\in W^+$ such that $g.C'=\vb.C_0^+$. One sets $\dw(C,C')=\vb$.  By \cite[Proposition 5.5]{rousseau2011masures} or \cite[Theorem 4.4.17]{hebert2018study}, $\vb$ does not depend on the choice of $A$. This defines a $G$-invariant ``$W$-distance'' $\dw:\CCC_0^+\times_\leq \CCC_0^+\rightarrow W^+$\index{$\dw$}.

 Let $C,C'$ be two chambers of the same sign and based at the same vertex. We say that $C$ and $C'$ are \textbf{adjacent} if they dominate a common panel. A gallery $\Gamma$ between $C$ and $C'$ is a finite sequence $\Gamma=(C_1,\ldots,C_n)$ such that $n\in \N$, $C_1=C$, $C_n=C'$ and $C_i,C_{i+1}$ are adjacent for every $i\in \llbracket 1,n-1\rrbracket$. The gallery $\Gamma$ is called \textbf{minimal} if $n$ is the minimum length among all the lengths of  the galleries joining $C$ to $C'$. If the vertex of $C$ and $C'$ is in $\I_0$, then the length of a minimal gallery between $C$ and $C'$ is $\ell(w)$, where $w=\dw(C,C')\in W^v$. 

\subsubsection{Iwahori subgroup and Iwahori-Hecke algebras associated with $G$}
Let $K_I$\index{$K_I$} be the fixator of $C_0^+$ in $G$. This is the \textbf{Iwahori subgroup} of $G$ (see also \cite[(3.8)]{braverman2016iwahori} for a more explicit description in the affine case). The map $g\mapsto g.C_0^+$ induces a bijection between $G/K_I$ and $\CCC_0^+$.  For $\wb\in W^v\ltimes Y$, we choose $n_\wb\in N$ such that $n_\wb$ induces $\wb$ on $\A$.   Then we have the Bruhat decomposition (see \cite[1.11]{bardy2016iwahori}): \[G^+= \bigsqcup_{\wb\in W^+} K_I n_\wb K_I.\]

In terms of masures, this decomposition has the following interpretation: for every $C,C'\in \CCC_0^+$ such that $\ve(C)\leq \ve(C')$, there exists an apartment containing $C$ and $C'$. Note that $\dw$ parametrizes the $K_I$ double cosets: if $g\in G^+$, then $g\in K_In_{\wb} K_I$ if and only if $\wb=\dw(C_0^+,g.C_0^+)$.

Let $\RCC$ be a ring. For $\wb\in W^+$, we denote by $T_{\wb}$\index{$T_\wb$} the indicator function of $K_I n_\wb K_I$.  Then the \textbf{Iwahori-Hecke algebra of }$G$ with coefficients in $\RCC$ is the free $\RCC$-module $\HC_{G,\RCC}$ with basis $(T_{\wb})_{\wb\in W^+}$ equipped with the product $*$ such that $T_{\mathbf{v}}*T_{\wb}=\sum_{\mathbf{u}\in W^+} a^{\mathbf{u}}_{\mathbf{v},\mathbf{w}}$, with $a^{\mathbf{u}}_{\mathbf{v},\mathbf{w}}=|(K_I n_\mathbf{v} K_I \cap n_\mathbf{u} K_I n_\mathbf{w}^{-1} K_I)/K_I|$ for $\mathbf{u},\mathbf{v},\mathbf{w}\in W^+$. The fact that such an algebra is well-defined is \cite[Theorem 2.4]{bardy2016iwahori} (the definition of the $T_{\wb}$ in \cite[2]{bardy2016iwahori} is slightly different but we obtain the same algebra).

Let $\FC$ be a field as in Definition~\ref{defBernstein-Lusztig algebra}. Let $q$ be the residue cardinal of $\KC$. As in \cite[5.7]{bardy2016iwahori}, we assume that there exists $\delta^{1/2}\in T_\FC$\index{$\delta^{1/2}$} such that $\delta^{1/2}(\alpha_s^\vee)=\sqrt{q}$ for every $s\in \SCC$. If $\FC=\C$, such a map exists by Lemma~\ref{lemExistence_character_values_alphavee}.  For $w\in W^v\subset W^+$, set  $H_w=q^{-\frac{1}{2}\ell(w)} T_w\in \HC_{G,\FC}$. For $\lambda\in Y^{++}$, set $Z^\lambda=\delta^{-\frac{1}{2}}(\lambda) T_{\lambda}\in \HC_{G,\FC}$. By \cite[5]{bardy2016iwahori}, we have the following proposition.

\begin{proposition}
 Let $\iota:\{Z^\lambda|\lambda\in Y^{++}\}\cup \{T_w|w\in W^v\}\subset \HC_{G,\FC}\rightarrow \AC_{\FC}$ be defined by $\iota(Z^\lambda)=Z^\lambda$ and $\iota(T_w)=T_w$ for $\lambda\in Y^{++}$ and $w\in W^v$. Then $\iota$ extends uniquely to an algebra morphism $\iota:\HC_{G,\FC}\rightarrow \AC_{\FC}$. Moreover, $\iota(\HC_{G,\FC})=\HC_\FC$ and $\iota$ is injective.
\end{proposition}

\subsubsection{Iwasawa decomposition and retractions centered at $\epsilon\infty$}

Let $\epsilon\in \{-,+\}$ and $U_\epsilon=\langle U_\alpha|\ \alpha\in \Phi_\epsilon\rangle$\index{$U_+,U_-$}. 
 We denote by $\epsilon\infty$\index{$-\infty$, $+\infty$} the germ of $\epsilon C^v_f$ at infinity: this is the filter on $\I$ composed with the sets containing a translate of $\epsilon C^v_f$. Then $U_\epsilon$ fixes $\epsilon\infty$, which means that for every $u\in U_\epsilon$, there exists $x\in \A$ such that $u$ fixes $x+\epsilon C^v_f$.
 
 Let $C$ be a chamber of $\I$. Then there exists an apartment containing $C$ and $\epsilon\infty$. This means that there exists $\Omega\in C$, $y\in \A$ and  an apartment containing  $\Omega\cup y+\epsilon C^v_f$. In particular for every $x\in \I$, there exists an apartment containing $x$ and $\epsilon\infty$. When $C\in \CCC_0^+$ and $x\in \I_0$, these results correspond to  the following decompositions: \[G=\bigsqcup_{\wb\in W^v\ltimes Y} U_\epsilon n_\wb K_I\text{ and } G=\bigsqcup_{\lambda\in Y} U_\epsilon n_\lambda K.\]

Let $x\in \I$. Let $A$ be an apartment containing $x$ and $\epsilon\infty$. Then by (MA ii), there exists $h\in G$ such that  $h.A=\A$ and $h$ fixes $A\cap \A$. We set $\rho_{\epsilon\infty}(x)=h.x$. This is well-defined, independently  of the choices of $A$ and $h$. Then $\rho_{\epsilon\infty}(x)$ is the unique element of $U_\epsilon.x\cap \A$. Then $\rho_{\epsilon\infty}:\I\rightarrow \A$\index{$\rho_{\epsilon\infty}$} is a retraction called the \textbf{retraction onto $\A$ centered at $\epsilon\infty$}.

\subsubsection{Towards principal series representations of $G^+$ and $G$}

Let $B=TU_+$\index{$B$} be the \textbf{positive standard Borel subgroup of $G$.} In term of masures, $B$ is stabilizer of $+\infty$ in $G$ (by \cite[Lemma 3.4.1]{hebert2018study}), which means that $B$ is the set of $g\in G$ such that there exists $a,a'\in \A$ such that $g.(a+C^v_f)=(a'+C^v_f)$ and such that there exists a translation $f$ of $\A$ such that $g.x=f(x)$ for every $x\in a+C^v_f$. Let $B^+=G^+\cap B$ and $T^+=T\cap G^+$.

\begin{lemma}
We have $T^+\subset B^+\subset T^+U_+$.
\end{lemma}

\begin{proof}
Let $g\in B^+$. Write $g=tu$ with $t\in T$ and $u\in U_+$. Then as $t$ normalizes $U_+$ (by \cite[8.3.3]{remy2002groupes}), there exists $u'\in U_+$ such that $g=u't$. Then $\rho_{+\infty}(g.0)=t.0$. Moreover by \cite[Corollaire 2.8]{rousseau2011masures}, $\rho_{+\infty}(g.0)\geq 0$ and thus $t.0\geq 0$, which proves the lemma.
\end{proof}

\begin{remark}
Unless $G$ is reductive, $T^+U_+\supsetneq B^+$. Indeed, let us prove that $U_+$ is not contained in $G^+$. Let $s\in \SCC$.  Take $a\in \KC$ such that $\omega(a)=-2$. Set $u=x_{\alpha_s}(a)\in U_+$. Let $A'=u.\A$. Then $A'\cap \A$ is the half-apartment $D_{\alpha_s,-2}=\{x\in \A|\alpha_s(x)-2\geq 0\}$. Let $D_{A'}$ be the half-apartment of ${A'}$ opposite to $D_{\alpha_s,-2}$. By \cite[Proposition 2.9 2)]{rousseau2011masures}, $\tilde{A}:=D_{-\alpha_s,2}\cup D_{A'}$ is an apartment of $\I$. As $0\notin D_{\alpha_s,-2}$, $u.0\in D_{A'}$. Then $\tilde{A}\ni 0,u.0$. Let $g\in G$ be such that $g.\tilde{A}=\A$ and such that $g$ fixes $D_{-\alpha_s,2}$. Let $r:\A\rightarrow \A$ be defined by $r(x)=s.x+2\alpha_s^\vee$ for $x\in \A$. Then by \cite[Lemma 3.4]{hebert2016distances}, $g.u.0=r.0=2\alpha_s^\vee$. By the lemma below, $g.u.0$ and $0=g.0$ are not comparable for $\leq$. We deduce that $u.0$ and $0$ are not comparable for $\leq$, which proves that $u\notin G^+$. 
\end{remark}

Recall the definition of indecomposable Kac-Moody matrices from \cite[§1.1]{kac1994infinite}.

\begin{lemma}
Assume that $G$ is associated with an indecomposable Kac-Moody matrix $A$ which is not a Cartan matrix. Then for all $s\in \SCC$, $\alpha_s^\vee\in \A\setminus(\T\cup -\T)$.
\end{lemma}

\begin{proof}
We first assume that $A$ is of affine type (see \cite[Theorem 4.3]{kac1994infinite} for the definition). Then there exists $\delta\in \bigoplus_{s\in \SCC} \R_+\alpha_s$ such that $\T=\delta^{-1}(\R^*_+)\sqcup \bigcap_{s\in \SCC} \alpha_s^{-1}(\{0\})$ (see \cite[Corollary 2.3.8]{hebert2018study}). By \cite[Proposition 5.2 a) and Theorem 5.6b)]{kac1994infinite}, $w.\delta=\delta$ for every $w\in W^v$.  Let  $x\in \A$ be such that $\delta(x)=0$ and $x\geq 0$. Then there exists $w\in W^v$ such that $w.x\in \overline{C^v_f}$. Then $\delta(x)=\delta(w.x)=0$. Thus $w.x\in \bigcap_{s'\in \SCC} \alpha_{s'}^{-1}(\{0\})$. As $\alpha_s(\alpha_s^\vee)=2$, $\alpha_s^\vee\notin \T$. As $s.\alpha_s^\vee=-\alpha_s^\vee$ we have $\alpha_s^\vee\in \A\setminus(\T\cup -\T)$.

We now assume that $A$ is of indefinite type. Then by \cite[Proposition 5.8 c)]{kac1994infinite} and \cite[2.9 Lemma]{gaussent2014spherical}, $\alpha_s^\vee\in \A\setminus \overline{\T}$. As $s.\alpha_s^\vee=-\alpha_s^\vee$ we deduce that $\alpha_s^\vee\in \A\setminus (\overline{\T}\cup -\overline{\T})$.
\end{proof}

Let $T_\FC^+= \mathrm{Hom}_{\mathrm{Mon}}(Y,\FC^*)$. Let $\tau\in T_\FC$ (resp. $\tau\in T_\FC^+$). We regard $\tau$ as a homomomorphism $T\rightarrow \FC^*$ (resp. as a monoid morphism $T^+\rightarrow \FC$) by setting $\tau(t)=\tau(t.0)$ for every $t\in T$ (resp. $t\in T^+$). We extend $\tau$ to a homomorphism $B\rightarrow \FC^*$ (resp. to a monoid morphism $B^+\rightarrow \FC$) by setting $\tau(tu)=\tau(t)$, for every $t\in T$ and $u\in U_+$ (resp $\tau(tu)=\tau(t)$ for every $t\in T^+$ and $u\in U_+$ such that $tu\in B^+$). By  \cite[Proposition 1.5 (DR5)]{rousseau2006immeubles} (note that there is a misprint in this proposition, $Z$ is in fact $T$), $T\cap U_+=\{1\}$. This implies that $\tau:B\rightarrow \FC^*$ is well-defined. The fact that  $\tau$ is a homomorphism follows from the fact that $t$ normalizes $U$ for every $t\in T$ (by \cite[8.3.3]{remy2002groupes}).

\begin{lemma}\label{lemDecomposition_borel_G}
\begin{enumerate}
\item Let $g\in G$ and $v\in W^v$. Then $g\in Bn_v K_I$ if and only if $\rho_{+\infty}(g.C_0^+)\in v.C_0^++Y$. In particular $G=\bigsqcup_{v\in W^v} Bn_vK_I$.

\item We have $G^+=\bigsqcup_{v\in W^v} B^+ n_v K_I$.
 \end{enumerate}
\end{lemma}

\begin{proof}
There exists $v\in W^v$ and $\lambda\in Y$ such that $\rho_{+\infty}(g.C_0^+)=v.C_0^++\lambda$. Thus there exists $t\in T$ and $v\in W^v$ such that $\rho_{+\infty}(g.C_0^+)=tn_v.C_0^+$. Hence $g.C_0^+=utn_v.C_0^+$ and $g\in utn_v K_I\subset Bn_vK_I$, for some $u\in U_+$. Conversely if $g\in Bn_vK_I$, then $\rho_{+\infty}(g.C_0^+)\in v.C_0^++Y$, which proves (1).

As $G^+$ is a sub-semi-group of $G$, $\bigsqcup_{v\in W^v} B^+ n_v K_I\subset G^+$. Let $g\in G^+$. By (1), we can write $g=bn_vk$, with $b\in B$, $v\in W^v$ and $k\in K_I$. Then $b.0=g.0\geq 0$ and hence $b\in B^+$, which proves (2).
\end{proof}

  \subsection{Action of $\HC_\FC$ on $I_{\tau,G^+}$ and $I_{\tau,G}$}\label{subAction_IH_algebra}

\subsubsection{Well-definedness of the action}\label{subsubWell_definedness_action}

Let $\epsilon\in \{+,\emptyset\}$. For $\tau\in T_\FC^\epsilon$ , we define $\widehat{I(\tau)^\epsilon}$  to be the set of functions $f$ from $G^\epsilon$  to $\FC$ such that for all $b\in B^\epsilon$  and $g\in G^\epsilon$, one has  $f(bg)=(\delta^{1/2}\tau)(b)f(g)$. The group $G$ (resp. semi-group $G^+$) acts on $\widehat{I(\tau)}$ (resp. $\widehat{I(\tau)^+}$) by right translation. When $G$ is reductive, the principal series representation associated with  $\tau$ is the subset $I(\tau)$ of functions of $\widehat{I(\tau)}$ which are locally constant.  Then $I_\tau=I(\tau)^{K_I}$. When $G$ is not reductive, we do not know which condition could replace ``locally constant''.   The hope is that the principal series representation of $G$ associated with $\tau$ should be the set of functions of $\widehat{I(\tau)}$ satisfying some ``regularity condition''.

Let $\tau\in T_\FC^\epsilon$. Let $\widehat{I(\tau)_{\mathrm{fin}}^\epsilon}$ be the set of $f\in \widehat{I(\tau)^\epsilon}$  such that there exists a finite set $F\subset W^v$ such that $\supp(f)\subset\bigcup_{v\in F} B n_v K_I$. Let $I_{\tau,G^\epsilon}=(\widehat{I(\tau)_{\mathrm{fin}}^\epsilon})^{K_I}$\index{$I_{\tau,G}$, $I_{\tau,G^+}$}  be the set of elements of $\widehat{I(\tau)_{\mathrm{fin}}^\epsilon}$  which are invariant under the action of $K_I$. For $v,w\in W^v$, define $f_w\in I_{\tau,G^\epsilon}$\index{$f_v$} by $f_w(n_v)=1$ if and only $v=w$. Then by Lemma~\ref{lemDecomposition_borel_G}, $(f_w)_{w\in W^v}$ is a basis of $I_{\tau,G^\epsilon}$. 

\medskip 

Fix $\tau\in T_{\FC}^\epsilon$. Following \cite[4.2.2]{bushnell2006local}, we would like to define an action of $\HC_\FC$ on $I_{\tau,G^\epsilon}$ by \[\phi.f= \sum_{g\in G^+/K_I} \phi(g) g.f,\forall (\phi,f)\in \HC_\FC\times I_{\tau,G^\epsilon}.\] However, we need to prove that such an action is well-defined. The main  difficulties are to prove that if $\phi\in \HC_\FC$, $f\in I_{\tau,G^\epsilon}$ and $h\in G$, then: \[\sum_{g\in G^+/K_I} \phi(g) f(hg)\] only involves finitely many terms and that $\phi.f$ also has finite support.  The aim of this section is to prove these results. For this, we use the masure $\I$, finiteness results of \cite{gaussent2008kac} and \cite{gaussent2014spherical} and the theory of Hecke paths introduced by Kapovich and Millson in \cite{kapovich2008path}.  In \cite{gaussent2008kac} and \cite{gaussent2014spherical}, the authors mainly use $\rho_{-\infty}$. As we use $\rho_{+\infty}$, we adapt their results to our framework.

Let $\lambda\in Y^{++}$.  A \textbf{$\lambda$-path} of $\A$ is a continuous piecewise linear map $\pi:[0,1]\rightarrow \A$ such that for every $t\in ]0,1[$, $\pi'_-(t),\pi'_+(t)\in W^v.\lambda$ (where $\pi'_-(t)$ and $\pi'_+(t)$ denote the left-hand and right-hand derivatives of $\pi$ at $t$) and  $\pi'_+(0),\pi'_-(1)\in W^v.\lambda$. A Hecke path of $\A$ of shape $\lambda$ with respect to $C^v_f$ is a $\lambda$-path satisfying \cite[1.8 Definition]{gaussent2014spherical}, with $\beta_i$ satisfying $\beta_i(C^v_f)<0$. Hecke paths are the images by retractions of preordered segments in $\I$. More precisely:

\begin{theorem}\label{thmHecke_path_retraction}(see \cite[Theorem 6.2]{gaussent2008kac})

Let  $x,y\in \I$ be such that $x\leq y$ and $\lambda=\dv(x,y)\in \overline{C^v_f}$. Let $\gamma:[0,1]\rightarrow  A$ be an affine parametrization of the segment $x,y$. Then $\rho_{+\infty}\circ \gamma$ is a Hecke path of shape $\lambda$ with respect to $C^v_f$ from $\rho_{+\infty}(x)$ to $\rho_{+\infty}(y)$. 
\end{theorem}

By definition of Hecke paths and by \cite[Lemma 1.3.13]{kumar2002kac}, we have the following lemma.

\begin{lemma}\label{lemBruhat_order_Hecke_paths}
Let $\lambda\in\overline{C^v_f}$ and $\pi:[0,1]\rightarrow \A$ be a Hecke path of shape $\lambda$ with respect to $C^v_f$. For $t\in [0,1]$ where it makes sense, we write $\pi'_+(t)=w'_+(t).\lambda$ and $\pi'_-(t)=w'_-(t).\lambda$, where $w'_-(t), w'_+(t)\in W^v$ have minimal lengths for these properties. Then for all $t,t'\in [0,1]$ such that $0\leq t< t'\leq 1$, we have $w'_-(t)\leq w'_+(t) \leq w'_-(t')\leq w'_+(t')$, where we delete the derivatives that do not make sense (for $t=0$ or $t'=1$).
\end{lemma}

\begin{theorem}\label{thmGR14}(see \cite[5.2]{gaussent2014spherical})
Let $x\in \I_0$, $\lambda\in Y^{++}$ and $\mu\in Y$. Then \[\{y\in \I_0|\ y\geq x,\  \dv(x,y)=\lambda\text{ and }\rho_{+\infty}(y)=\mu\}\] is finite. 
\end{theorem}

\begin{lemma}\label{lemFiniteness_chambers_retracting_vertex}
Let $y\in \I_0$ and $C$ be a type $0$ positive local chamber of $\A$.   Then \[\{C'\in \CO^+|\ve(C')=y\text{ and }\rho_{+\infty}(C')=C\}\] is finite.
\end{lemma}

\begin{proof}
Let $A$ be an apartment containing $y$ and $+\infty$. Then by (MA ii), there exists $g\in G$ such that $g.A=\A$ and $g$ fixes $A\cap \A$. Maybe working with $\rho_{+\infty,A}=g^{-1}.\rho_{+\infty}$ instead of $\rho_{+\infty}$, we can thus assume that $y$ is in $\A$. Let $C'\in \CO^+$ be such that $\ve(C')=y$ and $\rho_{+\infty}(C')=C$. Let $A'$ be an apartment containing $C'$ and $+\infty$. Then $A'$ contains $y$ and  by (MA ii), $A'$ contains $y+C_0^+$. Let $h\in G$ be such that $h$ fixes $A'\cap \A$ and $h.A'=\A$. Then $\rho_{+\infty}(C')=h.C'$. Therefore \begin{equation}\label{eqDistance_Retraction}
\dw(C',y+C_0^+)=\dw\big(h.C',h.(y+C_0^+)\big)=\dw(C,y+C_0^+)\in W^v.
\end{equation} Using \cite[Lemma 5.5]{abdellatif2019completed} we deduce that $\{C'\in \CO^+|\ve(C')=y\text{ and }\rho_{+\infty}(C')=C\}$ is finite.
\end{proof}

Let $x\in \I_0$ and $C\in \CCC_0^+$ be such that $C\geq x$ (i.e $\ve(C)\geq x$). By \cite[Proposition 5.17]{hebert2020new}, there exists an apartment $A$ containing $x$ and $C$. Then there exists $g\in G$ such that $g.A=\A$, $g.x=0$ and $g.C_0^+\in Y+C_0^+$. Then $g.\ve(C)\geq g.0$ and thus $g.\ve(C)\in Y^+$. One sets $\dy(0,C)=g.\ve(C)$. This does not depend on the choices we made by \cite[Theorem 4.4.17]{hebert2018study}. This defines a $G$-invariant ``distance'' $\dy:\I_0\times_{\leq} \CCC_0^+\rightarrow Y^+$\index{$\dy$}.

\begin{lemma}\label{lemFiniteness_chamber_support_well_definedness}
Let $v\in W^v$, $\lambda\in Y^{+}$. Then \[E:=\{C\in \CO^+|C\geq 0, \rho_{+\infty}(C)\in v.C_0^++Y\text{ and }\dy(0,C)=\lambda\}\] is finite.

Suppose moreover that $\lambda\in Y^{++}$ and that $v=1$. Then $E=\{\lambda+C_0^+\}$.
\end{lemma}

\begin{proof}
In order to prove that $E$ is finite, we begin by proving that $\ve(E):=\{\ve(C)|C\in E\}$ is finite.   To that end, our idea is to study, for each $C\in E$, the path $\tilde{\pi}=\rho_{+\infty}\circ \tilde{\gamma}:[0,1]\rightarrow \A$, where $\tilde{\gamma}$ is the segment joining $0$ to $\ve(C)$. We want to prove that $\tilde{\pi}'_{-}(1)$ lies in a finite set depending only on $v$ and $\lambda$.   In order to use the assumption that $\rho_{+\infty}(C)\in Y+v.C_0^+$, it is convenient to extend slightly the segment $\tilde{\gamma
}$ and this is why we consider a segment $\gamma:[0,1]\rightarrow \I$ such that $\gamma(0)=0$ and $\gamma(\frac{1}{2})=\ve(C)$.

Let $C\in E$. Let $A$ be an apartment containing $0$ and $C$. Let $g\in G$ be such that $g.\A=A$,  $g.0=0$ and $g.(\lambda+C_0^+)=C$.     Let $\gamma:[0,1]\rightarrow A$ be defined by $\gamma(t)=g.2t\lambda$.   Then $\pi=\rho_{+\infty}\circ \gamma$ is a Hecke path with respect to $+\infty$ of shape $2\lambda$.   Let $w_\lambda\in W^v$ be such that $(w_\lambda)^{-1}.\lambda\in Y^{++}$ and such that $w_\lambda$ has minimum length for this property.
  Set $C_\lambda=g.(\lambda+w_\lambda.C_0^+)$. Then: 
  \[\dw(C,C_\lambda)=\dw\big(g.(\lambda+C_0^+\big),g.(\lambda+w_\lambda.C_0^+)\big)=\dw(\lambda+C_0^+,\lambda+w_\lambda.C_0^+)=w_\lambda.\]
   Take a minimal gallery $\Gamma$ from $C$ to $C_\lambda$. Then $\Gamma$ has length $\ell(w_\lambda)$ and $\rho_{+\infty}(\Gamma)$ is a gallery from $\rho_{+\infty}(C)$ to $\rho_{+\infty}(C_\lambda)$. Therefore
    \[w:=\dw(\rho_{+\infty}(C),\rho_{+\infty}(C_\lambda)\big)\in W^v\text{ and }\ell(w)\leq \ell(w_\lambda).\]   Moreover, by definition of $E$, $\rho_{+\infty}(C)=\nu+v.C_0^+$, for some $\nu\in Y$. Consequently, $\rho_{+\infty}(C_\lambda)=\nu+vw.C_0^+$. Therefore for $\epsilon\in ]0,\frac{1}{2}]$ small enough, $\pi\big([\frac{1}{2},\frac{1}{2}+\epsilon]\big)\subset  \nu +vw.\overline{C^v_f}$  and thus $\pi'_+(\frac{1}{2})=2vw.\lambda$.    By  Lemma~\ref{lemBruhat_order_Hecke_paths}, $\pi'_-(\frac{1}{2})=u.\lambda$ for some $u\in W^v$ such that $\ell(u)\leq \ell(v)+\ell(w_\lambda)$.  
    
     Let now $\tilde{\gamma}:[0,1]\rightarrow A$ be defined by $\tilde{\gamma}(t)=g.t\lambda$ for $t\in[0,1]$ and $\tilde{\pi}=\rho_{+\infty}\circ \tilde{\gamma}$. Then by what we proved above, $\tilde{\pi}'_-(0)=u.\lambda$. By \cite[Lemma 1.8]{bardy2016iwahori} we have \[u.\lambda=\tilde{\pi}_-(1)\leq_{Q^\vee} \tilde{\pi}(1)-\tilde{\pi}(0)=\rho_{+\infty}\big(\ve(C)\big)\leq_{Q^\vee} \lambda^{++},\ \ell(u)\leq \ell(v)+\ell(w_\lambda),\] where $\lambda^{++}$ is the unique element of $Y^{++}\cap W^v.\lambda$. We deduce that  \[F:=\rho_{+\infty}\big(\ve(E)\big)=\{\rho_{+\infty}\big(\ve(C)\big)|C\in E\}\] is finite.
 
 Let $\nu\in  F$. Let $E_\nu=\{C\in E|\rho_{+\infty}(C)=\nu+v.C_0^+\}$. If $C\in E_\nu$, then $\dv\big(0,\ve(C)\big)=\lambda^{++}$  and $\rho_{+\infty}(\ve(C))=\nu$. Using Theorem~\ref{thmGR14} we deduce that $\{\ve(C)|C\in E_\nu\}$ is finite. By Lemma~\ref{lemFiniteness_chambers_retracting_vertex}, $E_\nu$ is finite and thus $E=\bigcup_{\nu\in F} E_\nu$ is finite.
 
 Suppose now that $v=1$ and that $\lambda\in Y^{++}$. Take $C\in E$. We use the same notation as in the beginning of the proof. Then we have $\pi'_{-}(\frac{1}{2})=\lambda=1.\lambda$ and by Lemma~\ref{lemBruhat_order_Hecke_paths} we deduce that there exists $\epsilon>0$ such that  $\pi(t)=2t\lambda$ for every $t\in [0,\frac{1}{2}+\epsilon]$. Moreover $\gamma(0)\in \A$ and thus by \cite[Lemma 3.4]{hebertGK} we deduce that $\gamma([0,\frac{1}{2}+\epsilon])\subset \A$. Therefore $C\subset \A$. Thus $\rho_{+\infty}(C)=C=\nu+C_0^+$ for some $\nu\in Y$. Moreover $d^{Y^+}(0,C)=\lambda+C_0^+$ and thus $\nu=\lambda$, which proves that $E=\{\lambda+C_0^+\}$ and completes the proof of the lemma.
\end{proof}

In the next lemma, we use the projection of a chamber on a vertex introduced in \cite[1.9]{bardy2016iwahori}. Let $x\in \A$ and $C$ be a positive chamber of $\A$ such that $y:=\ve(C)\geq x$.  Let $C^v$ be the positive vectorial chamber  of $\A$ such that $C=F_{y,C^v}$. Take $\xi\in C^v$. Then  there exists a positive vectorial chamber  $\tilde{C}^v\subset \A$ such that $x+\tilde{C}^v\supset \conv(x,]y,y+\epsilon\xi])$, for $\epsilon>0$ small enough, where $\conv$ denotes the convex hull. Then the chamber $\mathrm{pr}_x(C)=F_{x,\tilde{C}^v}$ is the projection of $C$ on $x$. Let now $x\in \I$ and $C$ be a positive chamber of $\I$ such that $\ve(C)\geq x$. Then there exists $g\in G$ such that $g.x,g.C\subset \A$. We set $\mathrm{pr}_x(C)=g^{-1}.\big(\mathrm{pr}_{g.x}(g.C)\big)$. This is the \textbf{projection of $C$ on $x$}. Then by \cite[Theorem 4.4.17]{hebert2018study}, $\mathrm{pr}_x(C)$ does not depend on the choice of $g$, every apartment containing $x$ and $C$ contains $\mathrm{pr}_x(C)$ and every $h\in G$ fixing $x$ and $C$ fixes $\mathrm{pr}_x(C)$.

\begin{lemma}\label{lemFinitess_well_definedness_action}
Let $\wb\in W^+$ and $v\in W^v$. Then: \begin{enumerate}
\item $\bigcup_{u\in W^v} (n_u K_I n_\wb K_I\cap B n_v K_I)/K_I$ is finite,

\item $\{u\in W^v| n_u K_I n_\wb K_I\cap B n_v K_I\neq \emptyset\}$ is finite.
\end{enumerate} 
\end{lemma}

\begin{proof}
Set $F=\bigcup_{u\in W^v} (n_u K_I n_\wb K_I\cap B n_v K_I)/K_I$. Let $u \in W^v$ and $g\in n_uK_In_\wb K_I$.  Set $C=g.C_0^+$. Then $\dw(u.C_0^+,C)=\wb$. Thus there exists $h\in G$ such that $h^{-1}.\A$ contains $u.C_0^+, C$ and such that $h.u.C_0^+=C_0^+, h.C=\wb.C_0^+$. Write $\wb=\lambda w$ (i.e $\wb.x=\lambda+w.x$ for every $x\in \A$). Set $h'=n_{w^{-1}} h$. Then $h'^{-1}.\A=h^{-1}.\A$ contains $0, C$, $h'.0=0$ and $h'.C=w^{-1}.\lambda+C_0^+$. Thus $\dy(0,C)=w^{-1}.\lambda$. Therefore  \[F.C_0^+ \subset \{C\in \CO^+|C\geq 0, \rho_{+\infty}(C)\in v.C_0^++Y\text{ and }\dy(0,C)=w^{-1}.\lambda\}.\] By Lemma~\ref{lemFiniteness_chamber_support_well_definedness}, $F.C_0^+$ is finite, which proves that $F$ is finite.

Let $u\in W^v$ be such that there exists $g\in n_uK_In_\wb K_I\cap Bn_v K_I$. Let $P=\{\mathrm{pr}_0(C')|C'\in F.C_0^+\}$. Let $C=g.C_0^+$. Then as $\dw(u.C_0^+,C)=\wb$, there exists $h\in G$ such that $h^{-1}.\A$ contains $u.C_0^+,C$, $h.u.C_0^+=C_0^+$ and $h.C=\wb.C$. Then $h.\mathrm{pr}_0(C)=\mathrm{pr}_0(\wb.C_0^+)$. Therefore \[w':=\dw\big(h.u.C_0^+,h.\mathrm{pr}_0(C)\big)=\dw\big(u.C_0^+,\mathrm{pr}_0(C)\big)=\dw\big(C_0^+,\mathrm{pr}_0(\wb.C_0^+)\big)\in W^v.\] Consequently there exists $C'\in P$ such that $\dw(u.C_0^+,C')=w'$. Consequently, \[\ell(u)= \ell\big(\dw(u.C_0^+,C_0^+)\big)\leq \ell(w')+\max_{C'\in P}\ell\big(\dw(C',C_0^+)\big).\]  This proves (2).
\end{proof}

\begin{definition/proposition}\label{def_prop_action}
Let $\epsilon\in \{+,\emptyset\}$ and $\tau\in T_{\FC}^\epsilon$.  Let $\phi\in \HC_\FC$ and $f\in I_{\tau,G^\epsilon}$. Define $\phi.f\in I_{\tau,G}$ by \[\phi.f=\sum_{g\in G^+/K_I} \phi(g) g.f.\] Then $.$ is well-defined and induces an action of $\HC_\FC$ on $I_{\tau,G^\epsilon}$.
\end{definition/proposition}

\begin{proof}
To prove that $\phi.f$ is  a well-defined element of $I_{\tau,G^\epsilon}$, it suffices to prove it for $\phi=T_{\wb}$ and $f=f_v$, for $v\in W^v$ and $\wb\in W^+$. Let $g\in G^+$ and $h\in G^\epsilon$. Suppose that  $T_{\wb}(g)f_v(hg)\neq 0$. Then $g\in K_In_{\wb} K_I\cap h^{-1} B n_v K_I$. Write $h=bn_u k$, with $b\in B^\epsilon$ and $k\in K_I$. Then $K_I n_\wb K_I\cap k^{-1}n_u^{-1} Bn_v K_I\neq \emptyset$. Therefore \begin{equation}\label{eqCondition_nonzero}g\in  K_I n_\wb K_I\cap k^{-1} n_u^{-1} B n_v K_I=k^{-1}(K_I n_\wb K_I\cap k^{-1} n_u^{-1} B n_v K_I).\end{equation} By Lemma~\ref{lemFinitess_well_definedness_action}, \[\sum_{g\in G^+/K_I} T_\wb(g) f_v(hg)=\sum_{g\in K_I n_\wb K_I\cap k^{-1}n_u^{-1} B n_v K_I/K_I} T_\wb(g)f_v(hg)\] is well-defined. Thus $T_\wb.f_v$ is a well-defined map $G^\epsilon\rightarrow \FC$. The fact that it is right $K_I$-invariant  and that $T_\wb.f(bh)=\delta^{1/2}\tau(b)T_\wb.f(h)$, for $B\in B^\epsilon$ are clear. 

Let $u\in W^v$. Suppose that  $T_\wb.f_v(n_u)\neq 0$. Then by \eqref{eqCondition_nonzero}, $K_I n_\wb K_I\cap n_u^{-1} B n_v K_I\neq  \emptyset$. By Lemma~\ref{lemFinitess_well_definedness_action} we deduce that $\{u\in W^v|\ T_{\wb}.f_v(n_u)\neq 0\}$ is finite, which proves that $T_\wb.f_v$ is an element of $I_{\tau,G^\epsilon}$. 

The fact that $(\phi*\phi').f=\phi.(\phi'.f)$ for every $f\in I_{\tau,G^\epsilon}$, $\phi,\phi'\in \HC_{\FC}$ is an easy consequence of the fact that $\phi*\phi'(h)=\sum_{g\in G^+/K_I}\phi(g)\phi'(g^{-1}h)$ for every $h\in G^+/K_I$.
\end{proof}

\subsubsection{Isomorphism between $I_\tau^\epsilon$ and $I_{\tau,G^\epsilon}$}

 Let $\tau:Y^+\rightarrow \FC$ be a monoid morphism. Then $\tau$ induces an algebra morphism $\tau:\FC[Y^+]\rightarrow \FC$ and thus this defines a representation $I_\tau^+=\mathrm{Ind}^{\HC_\FC}_{\FC[Y^+]}(\tau)=\HC_\FC\otimes_{\FC[Y^+]} \FC$.
Let $\epsilon\in \{+,\emptyset\}$. The aim of this section is to prove that  if $\tau\in (T_\FC)^\epsilon$ then  the map $I_{\tau}^\epsilon\rightarrow I_{\tau,G^\epsilon}$ defined by $h.1\otimes_\tau 1\mapsto h.f_1$, for $h\in \HC_\FC$ is well-defined and is an isomorphism of $\HC_{\FC}$-modules (see Proposition~\ref{propIso_group_theoretic_algebraic}). To that end, we prove that $Z^\lambda.f_1=\tau(\lambda) f_1$ for $\lambda\in Y^+$. For this we begin by proving that if $\lambda\in Y^{++}$, then $Z^\lambda.f_1=\tau(\lambda) f_1$. In the reductive case, this is sufficient to deduce the result for any $\lambda\in Y=Y^+$, since $Z^\lambda$ is invertible for $\lambda\in Y^{++}$. In the Kac-Moody case however, $Z^\lambda$ is not necessarily invertible for $\lambda\in Y^{++}$. We thus prove that if $f\in I_{\tau,G^\epsilon}$ is such that $Z^\lambda.f=0$ for $\lambda\in Y^{++}$ sufficiently dominant, then $f=0$.

\begin{lemma}\label{lemSupport_T_w_f_1}
Let $w\in W^v$. Then  $T_w.f_1=f_{w^{-1}}$.

\end{lemma}

\begin{proof}
Let $v\in W^v$.Then $T_w.f_1(n_v)=\sum_{g\in G^+/K_I} T_w(g)f_1(n_vg)$. Suppose that $T_w.f_1(n_v)\neq 0$. Then there exists $g\in K_In_w K_I\cap n_v^{-1} BK_I$ and thus $n_v K_I n_wK_I\cap B K_I\neq \emptyset$. 

Let $h\in n_vK_I n_wK_I\cap B K_I$ and $C=h.C_0^+$. Then $\dw(v.C_0^+,C)=w$ and $\rho_{+\infty}(C)\in Y+C_0^+$. Therefore $\ve(C)=0$ and hence $\rho_{+\infty}(C)=C_0^+$. By formula \eqref{eqDistance_Retraction} of the proof of Lemma~\ref{lemFiniteness_chambers_retracting_vertex}, we have $C=C_0^+$.  Consequently $C=C_0^+$, $v=w^{-1}$, $\supp(T_w.f_1)\subset B n_{w^{-1}} K_I$  and  $T_w.f(n_{w^{-1}})=1$. Therefore $T_w.f_1=f_{w^{-1}}$.
\end{proof}

\begin{lemma}\label{lemTriangular_support}
Let $w\in W^v$ and $\lambda\in Y\cap C^v_f$. Then: \begin{enumerate}
\item $\supp(T_\lambda.f_w)\subset  \bigcup_{v\leq w}B n_v K_I.$

\item  $T_\lambda.f_w(n_{w})\neq 0$.
\end{enumerate}
\end{lemma}

\begin{proof}
Let $v\in W^v$. Suppose that $T_\lambda.f_w(n_v)\neq 0$. Then $X:=n_v K_In_\lambda K_I\cap B n_w K_I$ is non-empty.  Let $g\in X$. Let $\gamma:[0,1]\rightarrow \I$ be defined by $\gamma(t)=g.t.\lambda$ for $t\in [0,1]$. Let $\pi=\rho_{+\infty}\circ \gamma$. Then $\pi$ is a Hecke path of shape $\lambda$ from $0$ to $\rho_{+\infty}\big(\ve(C)\big)$. For $t\in [0,1]$ where it makes sense, write $\pi'_-(t)=w_-(t).\lambda$, $\pi'_+(t)=w'_+(t).\lambda$, where $w_-'(t)$ and $w'_+(t)$ have minimum lengths for these properties.  By  the proof of Lemma~\ref{lemFiniteness_chamber_support_well_definedness}, $w'_-(1)\leq w$ (we have $w_\lambda=1$ in this case). Using Lemma~\ref{lemBruhat_order_Hecke_paths} we deduce that $w'_+(0)\leq w$. Let $C_{\pi(0^+)}$ (resp. $C_{\gamma(0^+)}$) be the local chamber based at $0$ and  containing $\pi(t)$ (resp. $\gamma(t)$)  for $t\in [0,1]$ near $0$. Then 
\[\dw(C_0^{+},C_{\gamma(0^+)})= \dw\big(\rho_{+\infty}(C_0^+),\rho_{+\infty}(C_{\gamma(0^+)})\big)=\dw(C_0^+,C_{\pi(0^+)})=w'_+(0).\]

Let us prove that $C_{\gamma(0^+)}=v.C_0^+$. Let $A$ be an apartment containing $v.C_0^+$ and $C$. Let $h\in G$ be such that $h.A=\A$ and such that $h$ fixes $v.C_0^+$. Then \[\begin{aligned} \dw(C_0^+,\lambda+C_0^+)&=\dw\big(h^{-1}.C_0^+,h^{-1}.(\lambda+C_0^+)\big)\\\ &=\lambda \\ &=\dw(v.C_0^+,h^{-1}.(\lambda+C_0^+)\big) \\ &=\dw(v.C_0^+,C).\end{aligned}\] As $A$ contains $v.C_0^+,C$ and $h^{-1}.(\lambda+C_0^+)$, we deduce that $h^{-1}.(\lambda+C_0^+)=C$. In particular, $h^{-1}.\lambda=g.\lambda$ and thus by \cite[Proposition 5.4]{rousseau2011masures}, $\gamma(t)=h^{-1}.t.\lambda$ for all $t\in [0,1]$. Let $\Omega'$ be a neighborhood of $0$ in $\A$ such that $h$ pointwise fixes  $\Omega=\Omega'\cap v.C^v_f$. Then for $t\in [0,1]$ small enough, $\gamma(t)\in \Omega$ and thus $C_{\gamma(0^+)}=v.C_0^+$. Consequently, $\gamma(t)\in \A$ for $t\in [0,1]$ small enough, thus $C_{\gamma(0^+)}\subset \A$, thus  $C_{\gamma(0^+)}=C_{\pi(0^+)}=v.C_0^+$ and hence $v=w'_+(0)\leq w$. Therefore: \[\supp(T_\lambda.f_w)\subset \bigcup_{v\leq w}B n_v K_I.\]

Suppose now that $v=w$. Then with the same notation as above, one has $w'_+(0)=w$. Therefore $w\leq w'_-(t) \leq w$ and $w\leq w'_+(t)\leq w$ for every $t\in [0,1]$ and hence $\pi$ is the line segment from $0$ to $w.\lambda$. Therefore if $g\in n_{w} K_In_\lambda K_I\cap B n_w K_I$, then $\rho_{+\infty}(g.C_0^+)=w.(\lambda+C_0^+)$. Consequently \[n_wK_In_\lambda K_I\cap B n_w K_I \subset U_+ n_{w.\lambda}n_w K_I,\] and   $n_{w.\lambda}\in T$. Thus
 \[ \begin{aligned} T_\lambda .f_w(n_w) & =\sum_{g\in K_In_\lambda K_I\cap n_w^{-1}Bn_w K_I/K_I}f_w(n_{w} g)  \\ &=|n_w K_I n_\lambda K_I\cap B n_w K_I/K_I|\tau\delta^{1/2}(w.\lambda).  \end{aligned}\]  Moreover $n_w n_\lambda\in  n_w K_In_\lambda K_I\cap B n_w K_I$, which proves that $T_\lambda .f_w(n_w)\neq 0$. 
\end{proof}

\begin{lemma}\label{lemFaithfullness_action}
Let $f\in I_{\tau,G^\epsilon}$. Suppose that for some $\mu\in Y\cap C^v_f$, $T_\mu.f=0$. Then $f=0$. 
\end{lemma}

\begin{proof}
Write $f=\sum_{w\in W^v} a_w f_w$, where $(a_w)\in \FC^{W^v}$ has finite support. Suppose that $f\neq 0$. Let $w\in \supp\big((a_v)\big)$ be maximal for the Bruhat order.  Then by Lemma~\ref{lemTriangular_support}, $T_\mu.f(n_w)=a_w T_\mu.f_w(n_w)\neq  0$. We reach a contradiction and thus $f=0$.
\end{proof}

\begin{lemma}\label{lemZlambda_weight_action}
Let $\lambda\in Y^+$. Then $Z^\lambda.f_1=\tau(\lambda).f_1$.
\end{lemma}

\begin{proof}
First assume that $\lambda\in Y^{++}$. Then $Z^\lambda=\delta^{-1/2}(\lambda) T_\lambda$, by \cite[5.7 and Theorem 5.5]{bardy2016iwahori}. By Lemma~\ref{lemTriangular_support}, $\supp(T_\lambda.f_1)=B K_I$ and thus $T_\lambda.f_1\in \FC f_1$. 

 We have $n_\lambda K_I\in K_I n_\lambda K_I\cap B K_I$. Let $g\in K_I n_\lambda K_I\cap B K_I$.  Let $C=g.C_0^+$. Then $\rho_{+\infty}(C)\in Y+C_0^+$ and $\dy(0,C)=\lambda$. Thus by Lemma~\ref{lemFiniteness_chamber_support_well_definedness}, $C=\lambda+C_0^+$. Hence $g\in n_\lambda K_I$ and $K_I n_\lambda K_I\cap B K_I=n_\lambda K_I$.
 Therefore $T_\lambda.f_1(1)=f_1(\lambda)=\delta^{1/2}\tau(\lambda)$. Hence $T_\lambda.f_1=\delta^{1/2}\tau(\lambda)f_1$ and $Z^\lambda.f_1=\tau(\lambda)f_1$.

Let now $\lambda\in Y^+$. Then by \cite[Theorem 5.5]{bardy2016iwahori} and the fact that $Z^\lambda=\delta^{-1/2}(\lambda) X^\lambda$, one has $T_\mu .Z^\lambda.f_1=\delta^{-1/2}(\lambda)T_{\lambda+\mu}.f_1=\tau(\lambda+\mu)\delta^{1/2}(\mu) f_1=T_\mu.(\tau(\lambda).f_1)$ for $\mu\in Y^{++}$ sufficiently dominant. Thus by Lemma~\ref{lemFaithfullness_action}, $Z^\lambda.f_1=\tau(\lambda).f_1$, which proves the lemma.
\end{proof}

\begin{proposition}\label{propIso_group_theoretic_algebraic}
Let $\epsilon\in \{+,\emptyset\}$. Let $\tau\in T_\FC^\epsilon$. Then the map $\phi:I_{\tau}^{\epsilon}\rightarrow I_{\tau,G^\epsilon}$ defined by $\phi(h.1\otimes_\tau 1)\mapsto h.f_1$  for $h\in \HC_\FC$ is well-defined and is an isomorphism of $\HC_{\FC}$-modules.
\end{proposition}

\begin{proof}
By Lemma~\ref{lemFrobenius_reciprocity} and Lemma~\ref{lemZlambda_weight_action}, $\phi$ is well-defined. Let $x\in I_\tau^\epsilon$ be such that $\phi(x)=0$. Write $x=\sum_{v\in W^v} a_v T_v\otimes_\tau 1$, with $(a_v)\in \FC^{W^v}$. Then $\phi(x)=\sum_{v\in W^v} a_v T_v.f_1$. Suppose that $x\neq 0$. Let $w\in W^v$ be such that $a_w\neq 0$ and such that $w$ is maximal for this property (for the Bruhat order). Then by Lemma~\ref{lemTriangular_support} and Lemma~\ref{lemSupport_T_w_f_1}, $\phi(x)(n_{w^{-1}})=a_w T_w.f_1(n_{w^{-1}})\neq 0$: a contradiction. Therefore $x=0$ and $\phi$ is injective. By Lemma~\ref{lemSupport_T_w_f_1} and Lemma~\ref{lemDecomposition_borel_G},  $(T_w.f_1)_{w\in W^v}$ is a basis of $I_{\tau,G^\epsilon}$. Consequently $\phi$ is surjective, which proves the proposition. 
\end{proof}

\subsection{Extendability of representations of $G^+$ and $\HC_\FC$}\label{subExtendability}

In this subsection, we study the extendability of $I_{\tau^+,G^+}$ (resp. $I_{\tau^+}^+$) to a representation of $G$ (resp. $\AC_\FC$), for $\tau\in T_\FC^+$. We obtain a criterion depending on the extendability of $\tau^+$ to an element of $T_\FC$ (see Proposition~\ref{propExtension_I_tau_G_+toG}). 

\subsubsection{Extendability of  elements of $T_\FC^+$}

Recall that if  $\tau:Y^+\rightarrow \FC$ is a monoid morphism $I_\tau^+=\mathrm{Ind}^{\HC_\FC}_{\FC[Y^+]}(\tau)=\HC_\FC\otimes_{\FC[Y^+]} \FC$ is a  representation of $\HC_\FC$. If $I_\tau^+$ is not the restriction of a representation of $\AC_\FC$ we call $I_\tau^+$ a \textbf{non-extendable principal series representation of $\HC_\FC$}. In this section we study the existence of non-extendable principal series representations of $\HC_\FC$. We prove that in some cases - for example when $\HC_\FC$ is associated with an affine root generating system or to a size $2$ Kac-Moody matrix - every  principal series representations of $\HC_\FC$ can be extended to a representation of $\AC_\FC$ (see Lemma~\ref{lemNonexistence_degenerate_representations}). We prove that there exist Kac-Moody matrices such that $\HC_\FC$ admits non-extendable principal series representations (see Lemma~\ref{lemExistence_degenerate_representations}).

   Let $\mathrm{res}_{Y^+}:\mathrm{Hom}_{\mathrm{Mon}}(Y,\FC)\rightarrow \mathrm{Hom}_\mathrm{Mon}(Y^+,\FC)$  be defined by $\mathrm{res}_{Y^+}(\tau)=\tau_{|Y^+}$ for all $\tau\in \mathrm{Hom}_\mathrm{Mon}(Y,\FC)$.

\begin{lemma}\label{lemRestriction_bijection}
The map $\mathrm{res}_{Y^+}:\mathrm{Hom}_\mathrm{Gr}(Y,\FC^*)=\mathrm{Hom}_\mathrm{Mon}(Y,\FC^*)\rightarrow \mathrm{Hom}_\mathrm{Mon}(Y^+,\FC^*)$ is a bijection.
\end{lemma}

\begin{proof}
Let $\tau\in \Hom_\mathrm{Mon}(Y,\FC^*)$. Let $\nu\in C^v_f$. Let $\lambda\in Y$ and $n\in \N$ be such that $\lambda+n\nu\in \T$. Then $\tau(\lambda)=\frac{\tau(\lambda+n\nu)}{\tau(n\nu)}$ and thus $\mathrm{res}_{|Y^+}$ is injective.

Let $\tau^+\in \mathrm{Hom}_{\mathrm{Mon}}(Y^+,\FC^*)$. Let $\lambda\in Y$. Write $\lambda=\lambda_+-\lambda_-$, with $\lambda_+,\lambda_-\in Y^+$. Set $\tau(\lambda)=\frac{\tau^+(\lambda_+)}{\tau^+(\lambda_-)}$, which does not depend on the choices of $\lambda_-$ and $\lambda_+$. Then $\tau\in \mathrm{Hom}_{\mathrm{Mon}}(Y,\FC^*)$ is well-defined and $\mathrm{res}_{|Y^+}(\tau)=\tau^+$, which finishes the proof.
\end{proof}

\begin{lemma}
Let $\tau\in \mathrm{Hom}_{\mathrm{Mon}}(Y^+,\FC)$ and $\chi\in T_\FC$. \begin{enumerate}

\item\label{itMorphismItauItau+} Suppose $\mathrm{Hom}_{\HC_\FC-\mathrm{mod}}(I_\tau^+,I_\chi)\neq \{0\}$.  Then there exists $w\in W^v$ such that $\tau=w.\chi_{|Y^+}$.

\item\label{itMorphismItau+Itau} Suppose $\mathrm{Hom}_{\HC_\FC-\mathrm{mod}}(I_\chi,I_\tau^+)\neq \{0\}$.  Then there exists $w\in W^v$ such that $\tau=w.\chi_{|Y^+}$.
\end{enumerate} 
\end{lemma}

\begin{proof}
(~\ref{itMorphismItauItau+}) Let $\phi\in \mathrm{Hom}_{\HC_\FC-\mathrm{mod}}(I_\tau^+,I_\chi)\setminus \{0\}$. Let $x=\phi(1\otimes_{\tau^+}1)$. Then $Z^\lambda.x=\tau(\lambda).x$ for all $\lambda\in Y^+$. By Lemma~\ref{lemCommutation relation}, $Z^\lambda.x\neq 0$ for all $\lambda\in Y^+$. Thus $\tau(\lambda)\neq 0$ for all $\lambda\in Y^+$.

 Let $\mu\in Y$. Let $\nu\in C^v_f\cap Y$ be such that $\mu+\nu\in Y^+$. Then $Z^\mu.x=\frac{\tau(\mu+\nu)}{\tau(\nu)}.x$. Therefore there exists $\chi'\in T_\FC$ such that $x\in I_\chi(\chi')$. By Lemma~\ref{lemDecomposition_submodules_I_tau}, $\chi'\in W^v.\chi$. Moreover, $\chi'_{|Y^+}=\tau$, which proves~(\ref{itMorphismItauItau+}).

(\ref{itMorphismItau+Itau}) Let $\phi\in \mathrm{Hom}_{\HC_\FC-\mathrm{mod}}(I_\chi,I_\tau^+)\setminus \{0\}$. Let $x=\phi(1\otimes_{\chi} 1)$. Then $Z^\lambda.x=\chi(\lambda).x$ for all $\lambda\in Y^+$. By a lemma similar to Lemma~\ref{lemDecomposition_submodules_I_tau} we deduce that $\chi_{|Y^+}\in W^v.\tau$, which proves the lemma.
\end{proof}

One has $\mathrm{Hom}_{\mathrm{Mon}}\big(Y,(\FC,.)\big)=\mathrm{Hom}_{\mathrm{Gr}}(Y,\FC^*)\cup \{0\}$. Set $\A_{in}=\bigcap_{s\in \SCC}\ker( \alpha_s)$. Let $\mathring{\T}$ be the interior of the Tits cone.

\begin{lemma}\label{lemNonexistence_degenerate_representations}
Let $\tau^+\in \mathrm{Hom}_{\mathrm{Mon}}\big(Y,(\FC,.)\big)$. Assume that there exists $\lambda\in Y^+$ such that $\tau^+(\lambda)=0$. Then $\tau^+(\mathring{\T}\cap Y)=\{0\}$. In particular, if $\T=\mathring{\T}\cup \A_{in}$, then $\mathrm{Hom}_{\mathrm{Mon}}\big(Y^+,(\FC,.)\big)=\mathrm{Hom}_{\mathrm{Mon}}(Y,\FC^*)\cup \{0\}$.
\end{lemma}

\begin{proof}
Let $\mu\in \mathring{\T}\cap Y$. Then for $n \gg 0$, $n\mu\in \lambda+ \T$. Indeed, $n\mu-\lambda=n(\mu-\frac{\lambda}{n})\in \T$ for $n\gg 0$. Hence $\tau^+(n\mu)=(\tau^+(\mu))^n=0$.
\end{proof}

A face $F^v\subset \T$ is called \textbf{spherical} if its fixator in $W^v$ is finite.

\begin{remark}
\begin{enumerate}
\item If $\A$ is associated to an affine Kac-Moody matrix, then $\T=\mathring{\T}\cup\A_{in}$ (see \cite[Corollary 2.3.8]{hebert2018study} for example).

\item  If $\A$ is associated to a size $2$ indefinite Kac-Moody matrix, then $\T=\mathring{\T}\cup\A_{in}$. Indeed, by \cite[Th{\'e}or{\`e}me 5.2.3 ]{remy2002groupes}, $\mathring{\T}$ is the union of the spherical vectorial faces. By \cite[1.3]{rousseau2011masures}, if $J\subset \SCC$ and $w\in W^v$, the fixator of $w.F^v$ is $w.W^v(J).w^{-1}$. Therefore the only non-spherical face of $\T$ is $\A_{in}$ and hence $\T=\mathring{\T}\cup \A_{in}$.

\item Let $A=(a_{i,j})_{i,j\in \llbracket 1,3\rrbracket}$ be a Kac-Moody matrix such that for all $i\neq j$, $a_{i,j}a_{j,i}\geq 4$. Then by \cite[Proposition 1.3.21]{kumar2002kac}, $W^v$ is the free group with 3 generators $s_1,s_2,s_3$ of order $2$. Thus for all $J\subset \SCC$ such that $|J|=2$, $F^v(J)$ is non-spherical. Hence $\T\supsetneq\mathring{\T}\cup \A_{in}$.
\end{enumerate}
\end{remark}

\subsubsection{Construction of an element of $\Hom_{\mathrm{Mon}}(Y^+,\FC)\setminus \Hom_{\mathrm{Mon}}(Y,\FC)$}

We now prove that there exist Kac-Moody matrices for which \[\Hom_{\mathrm{Mon}}(Y^+,\FC)\neq \Hom_{\mathrm{Mon}}(Y,\FC).\] Assume that $\A$ is associated to an invertible indefinite size 3 Kac-Moody  matrix (see \cite[Theorem 4.3]{kac1994infinite} for the definition of indefinite). Then one has $\A=\A'\oplus \A_{in}$, where $\A'=\bigoplus_{i\in I}\R\alpha_i^\vee$. Maybe considering $\A/\A_{in}$, we may assume that $\A_{in}=\{0\}$.

Recall that $\T$ is the disjoint union of the positive vectorial faces of $\A$. 

\begin{lemma}\label{lemExtremality_Nonspherical_face}
Assume that there exists a non-spherical vectorial face $F^v\neq \{0\}$. Let $x\in \T$ and $y\in \T\setminus \overline{F^v}$. Then $[x,y]\cap F^v\subset \{x\}$.
\end{lemma}

 \begin{proof}
 Assume  that $y\in \mathring{\T}$. Then $(x,y]\subset \mathring{\T}$ and thus $[x,y]\cap F^v\subset \{x\}$. 
 
 Assume that $y\notin \mathring{\T}$. For $a\in \T$, we denote by $F^v_a$ the vectorial face of $\T$ containing $a$. If $F^v_x=F^v_y$, then $[x,y]\subset F^v_x$. As $F^v_y\neq F^v$, we deduce that $[x,y]\cap F^v=\emptyset$. We now assume that $F^v_x\neq F^v_y$. As $W^v$ is countable, the number of positive vectorial faces is countable and thus there exist $u\neq u'\in [x,y]$ such that $F^v_u=F^v_{u'}$. Then the dimension of the vector space spanned by $F^v_{u}$ is at least $2$. Thus there exists $w\in W^v$ such that  $F^v_u=w.F^v(J)$, for some $J\subset \SCC$ such that $|J|\leq 1$. Then the fixator of $F^v_u$ is $w.W_J.w^{-1}$, where $W_J=\langle J\rangle $. Then $W_J$ is finite and thus $F^v_u$ is spherical.  Consequently, $(x,y)=(x,u]\cup [u,y)\subset \mathring{\T}$ and the lemma follows.
 \end{proof}

 \begin{lemma}\label{lemConvexity_Tits_cone_minus_nonspherical}
Assume that   there exists a non-spherical vectorial face $F^v\neq \{0\}$. Then  $\T\setminus \overline{F^v}$ and $\T\setminus\{0\}$ are convex. 
\end{lemma}

\begin{proof}
Let $x,y\in \T\setminus F^v$. Suppose that $ [x,y]\cap F^v\neq \emptyset$. By Lemma~\ref{lemExtremality_Nonspherical_face}, $y\in \overline{F^v}=F^v\cup\{0\}$ and hence $y=0$. Let $F^v_x$ be the vectorial face containing $x$. Then $[x,y)\subset F^v_x$ and hence $[x,y)\cap F^v=\emptyset$: a contradiction. Thus $\T\setminus F^v$ is convex.

By \cite[2.9 Lemma]{gaussent2014spherical}, there exists a basis $(\delta_s)_{s\in \SCC}$ of $\bigoplus_{s\in \SCC}\R\alpha_s^\vee$  such that $\delta_s(\T)\geq 0$ for all $s\in \SCC$. Thus $\T\setminus\{0\}$ is convex and hence $\T\setminus \overline{F^v}=\T\setminus F^v\cap \T\setminus \{0\}$ is convex.
\end{proof}

 \begin{lemma}\label{lemExistence_degenerate_representations}
 Assume that $\A$ is associated with an indefinite Kac-Moody matrix of size $3$ such that there exists a non-spherical face different from $\A_{in}$. Assume moreover that $(\alpha_s^\vee)_{s\in \SCC}$ is a basis of $\A$.  Then $\mathrm{Hom}_{\mathrm{Mon}}\big(Y^+,(\FC,.)\big)\supsetneq \mathrm{Hom}_{\mathrm{Mon}}\big(Y^+,\FC^*\big)\cup \{0\}$.
 \end{lemma}
 
 \begin{proof}
 Let $\tau^+=\mathds{1}_{\overline{F^v}}:\T\rightarrow  \FC$. Let us prove that $\tau^+\in \mathrm{Hom}_{\mathrm{Mon}}\big(\T,(\FC,.)\big)$. 
 
 Let $x,y\in \T$. If $x,y\in \T\setminus \overline{F^v}$, then $x+y=2.\frac{1}{2}(x+y)\in \T\setminus \overline{F^v}$ by Lemma~\ref{lemConvexity_Tits_cone_minus_nonspherical} and thus $\tau^+(x+y)=0=\tau^+(x)\tau^+(y)$. 
 
 Suppose $x\in F^v$ and $y\in \T\setminus \overline{F^v}$, then $x+y=2.\frac{1}{2}(x+y)\in \T\setminus \overline{F^v}$ by Lemma~\ref{lemExtremality_Nonspherical_face}. Thus $\tau^+(x+y)=0=\tau^+(x)\tau^+(y)$. 
 
 Suppose $x=\{0\}$ and $y\in \T\setminus \overline{F^v}$. Let $F^v_y$ be the vectorial face containing $y$. Then $(x,y]\subset F^v_y$ and hence $x+y\in F^v_y$: $\tau^+(x+y)=0=\tau^+(x)\tau^+(y)$. Consequently, $\tau^+\in \mathrm{Hom}_{\mathrm{Mon}}\big(\T,(\FC,.)\big)$.

Maybe considering $w.F^v$, for some $w\in W^v$, we can assume $F^v\subset \overline{C^v_f}$. Then there exist $s_1,s_2,s_3\in \SCC$ such that $\SCC=\{s_1,s_2,s_3\}$ and $F^v=\alpha_{s_1}^{-1}(\{0\})\cap \alpha_{s_2}^{-1}(\{0\})\cap \alpha_{s_3}^{-1}(\R^*_+)$. Let $\lambda\in \A$ be such that $\alpha_{s_1}(\lambda)=\alpha_{s_2}(\lambda)=0$ and $\alpha_{s_3}(\lambda)=1$. There exists $n\in \Ne$ such that $\lambda\in \frac{1}{n}Y$. Thus $\tau^+_{|Y^+}\in \mathrm{Hom}_{\mathrm{Mon}}\big(Y^+,(\FC,.)\big)\setminus(\mathrm{Hom}_{\mathrm{Mon}}\big(Y^+,\FC^*\big)\cup \{0\})$.
 \end{proof}

 \subsubsection{Extension of the representations from $G^+$ to $G$}

We now study under which condition the representation $I_{\tau,G^+}$ of $G^+$ extends to a representation of $G$, for $\tau\in T_\FC^+$. 

\begin{lemma}\label{lemG_+_spans_G}
Let $g\in G$. Then for $t\in T$ such that $t.0$ is sufficiently dominant, $tg\in G^+$.
\end{lemma}

\begin{proof}
Let $g\in G$ and $x=g.0$. There exists an apartment containing $-\infty$ and $x$, i.e there exists $g\in G$ such that $g.\A\cap \A$ contains $a-C^v_f$, for some $a\in \A$. For $q\in C^v_f$ sufficiently dominant, $a-q\leq x$. In particular, there  exists $y\in \A$ such that $y\leq x$. For $\lambda\in Y^{++}$ sufficiently dominant, $y+\lambda\geq 0$. Then $n_\lambda.y=y+\lambda\geq 0$. As $\leq$ is $G$-invariant, $n_\lambda.y\leq n_\lambda.x$ and thus $0\leq n_\lambda.x=n_\lambda g.0$. Therefore $n_\lambda g\in G^+$. 
\end{proof}

Let $x,y\in \I$. We write $x\mathring{<}y$ (resp. $x\mathring{\leq} y$) if there exists $g\in G$ such that $gx,g.y\in \A$ and $y-x\in \mathring{\T}$ (resp. $y-x\in \mathring{\T}\cup\{0\}$). This does not depend on the choice of $g$.

If $G$ is reductive, then $x\leq y$ for every $x,y\in \I$. We now assume that $G$ is not reductive. Then for every $x\in \A$, for every $y\in x+C^v_f$, one has $x\mathring{<} y$ and $y\not\leq x$.

\begin{lemma}\label{lemRefinement_transitivity_preorder}
Let $x,y,z\in \I$. Suppose that $x\leq y$,  $y\mathring{<} z$ and  $z\not\leq y$. Then $x\mathring{<} z$.
\end{lemma}

\begin{proof}
Let $A$ be an apartment containing $y$ and $z$. Let $F_y$ be a positive face of $A$ based at $y$ and containing $[y,y']$ for $y'\in [y,z]$ near $y$. Then by \cite[Theorem 4.4.17]{hebert2018study}, there exists an apartment $A'$ containing $F_y$ and $x$. Then $A'$ contains $[y,y']$ for some $y'\in [y,z]$ near $y$. In the apartment $A'$, one has $y\mathring{<} y'$ and $x\leq y$. Consequently $x \mathring{<} y'$ (because $\mathring{\T}+\T\subset \mathring{\T}$). We thus have $x\mathring{\leq} y'$ and $y'\mathring{\leq} z$. Using \cite[Théorème 5.9]{rousseau2011masures} we deduce that $x\mathring{\leq} z$. As $x\leq y$ and $z\not \leq y $, we have $x\neq z$, which proves the result.
\end{proof}

\begin{lemma}\label{lemExtensions_representations}
\begin{enumerate}
\item Let $\tau\in T_\FC^+$ be such that $\tau$ is the restriction of some element of $T_\FC$ (still denoted $\tau$). Then every element of $\widehat{I(\tau)^+}$ uniquely extends to an element of $\widehat{I(\tau)}$.

\item Let $\tau \in T_\FC^+$ be such that $\tau$ is not the restriction of some element of $T_\FC$. Then for every $f:G\rightarrow \FC$ such that for all $g\in G^+$ and $b\in B^+$, $f(bg)=(\delta^{1/2}\tau)(b)f(g)$, one has $f=0$.

\item Let $\tau\in T_\FC^+$ be such that $\tau$ is not the restriction of some element of $T_\FC$. Then there exists $t\in T$ such that for every $f\in I_{\tau,G^+}$, $t.f=0$. 
\end{enumerate}
\end{lemma}

\begin{proof}
(1) Let $f\in  \widehat{I(\tau)^+}$. Suppose that there exists $\tilde{f}\in  \widehat{I(\tau)}$ extending $f$. Let $g\in G$. Let $t\in T$ be such that $tg\in G^+$. Then $\tilde{f}(tg)=(\delta^{1/2}\tau)(t)\tilde{f}(g)=f(tg)$ and thus $\tilde{f}(g)=(\delta^{1/2}\tau(t))^{-1}f(tg)$. Thus $\tilde{f}$ is unique if it exists. 

We now set $f'(g)=(\delta^{1/2}\tau(t))^{-1}f(tg)$, for  $t\in T$ such that $t.0$ is dominant and such that $tg\in G^+$, which exists by Lemma~\ref{lemG_+_spans_G}. Let us prove that $f'$ is well-defined. Let $t,t'\in T$ be such that $tg,t'g\in G^+$ and such that $t.0,t'.0\in Y^{++}$. Then \[f(tt'g)=(\tau\delta^{1/2})(t')f(tg)=(\tau\delta^{1/2})(t)f(t'g)\] so that $f(t'g)\big(\tau\delta^{1/2}(t')\big)^{-1}=f(tg)\big(\tau\delta^{1/2}(t)\big)^{-1}$. This prove that $f'$ is well-defined. In particular, $f'$ extends $f$.

Let now $t\in T$ and $g\in G$. Let us prove that $f'(tg)=\tau\delta^{1/2}(t)f'(g)$. Let $t'\in T$ be such that $t'g,t'tg\in G^+$. Then \[f'(g)=f(tt'g)\big(\delta^{1/2}\tau(tt')\big)^{-1}=\tau\delta^{1/2}(t')f'(tg)\big(\tau\delta^{1/2}(tt')\big)^{-1}=f'(tg)\big(\tau\delta^{1/2}(t)\big)^{-1},\] which proves that  $f'(tg)=\tau\delta^{1/2}(t)f'(g)$.

Let now $g\in G^+$ and $u\in U_+$. Let $t\in T$ be such that $tg,tu\in G^+$. Then $f'(tug)=\tau\delta^{1/2}(t)f'(ug)$ and $f'(tug)=\tau\delta^{1/2}(tu)f'(g)=\tau\delta^{1/2}(t)f(g)$. Thus \[\tau\delta^{1/2}(t)f'(g)=\tau\delta^{1/2}(t)f'(ug)\] and hence $f'(ug)=f'(g)$ for every $u\in U_+$ and $g\in G^+$.

Let now $g\in G$ and $u\in U_+$. Let $t\in T$ be such that $tug,tg\in G^+$. As $t$ normalizes $U_+$, we can write $tu=u't$ for some $u'\in U_+$. Then \[f'(ug)=f'(tug)\big(\tau\delta^{1/2}(t)\big)^{-1}=f'(u'tg)\big(\tau\delta^{1/2}(t)\big)^{-1}=f'(tg)\big(\tau\delta^{1/2}(t)\big)^{-1}=f'(g).\]

Let $b\in B$ and $g\in G$. Write $b=tu$, with $t\in T$ and $u\in U_+$. Then  we have \[f'(bg)=f'(tug)=\tau\delta^{1/2}(t)f'(ug)=\tau\delta^{1/2}(t)f'(g)=\tau\delta^{1/2}(b)f'(g)\] and thus $f'\in \widehat{I(\tau)}$ and $f'$ extends $f$. This proves (1).

(2) Let $\tau \in T_\FC^+$ be such that $\tau$ is not the restriction of some element of $T_\FC$. Then by Lemma~\ref{lemRestriction_bijection}, there exists $t\in T$ such that $\tau(t)=0$. Let $f:G\rightarrow \FC$ be such that for all $g\in G^+$ and $b\in B^+$, $f(bg)=(\delta^{1/2}\tau)(b)f(g)$. Let $g\in G$. Then $f(g)=f(tt^{-1}g)=\tau\delta^{1/2}(t)f(t^{-1}g)=0$, which proves (2).

(3) By Lemma~\ref{lemNonexistence_degenerate_representations}, one has $\tau(t')=0$ for every $t'\in T$ such that $t'.0\in \mathring{\T}$. Let $t\in T$ be such that $t.0\in C^v_f$. Let $g\in G^+$ and $f\in I_{\tau,G^+}$. Then $t.0\mathring{>} 0$ and $t.0\not\leq 0$. Therefore $gt.0\mathring{>} g.0$ and $gt.0\not\leq g.0$. Moreover $g.0\geq 0$ and thus by Lemma~\ref{lemRefinement_transitivity_preorder} we have $gt.0\mathring{>} 0$.  Using Lemma~\ref{lemDecomposition_borel_G} we write $gt=bn_vk$, with $b\in B^+$, $v\in W^v$ and $k\in K_I$. Then $gt.0=b.0$, which proves that $b.0\mathring{>} 0$. Write $b=u't'$, with $u'\in U_+$ and $t'\in T$. Then by Theorem~\ref{thmHecke_path_retraction}, $\rho_{+\infty}(b.0)=t'.0\mathring{>} 0$ and thus $\tau(t')=0$. Therefore $f(gt)=t.f(g)=\tau\delta^{1/2}(t')f(n_vk)=0$, which proves (3).
\end{proof}

\begin{proposition}\label{propExtension_I_tau_G_+toG}

Let $\tau^+\in T_\FC^+$. 
\begin{enumerate}

\item Suppose that $\tau^+$ is not the restriction to $Y^+$ of an element of $T_\FC$. 

For every $f\in \widehat{I(\tau^+)}\setminus\{0\}$, for every $G$-module $M$, the restriction of $M$ to $G^+$ is not isomorphic to $G^+.f$. 

For every $x\in I_{\tau^+}^+\setminus\{0\}$, for every $\AC_\FC$-module $M$, the restriction of  $M$ to $\HC_\FC$ is not isomorphic to $\HC_\FC.x$.

\item Suppose that $\tau^+$ is the restriction to $Y^+$ of a (necessarily unique) element $\tau$ of $T_\FC$. 

Every element $f^+$ of $\widehat{I(\tau^+)^+}$ can be extended uniquely to an element $f$ of $\widehat{I(\tau)}$. Then $f^+\mapsto  f$ is an isomorphism of $G^+$-modules.

 The action of $\HC_\FC$ on  $I_{\tau^+}^+$ extends uniquely to an action of $\AC_\FC$ on $I_{\tau^+}^+$. Then $I_{\tau^+}^+$ is naturally isomorphic to $I_\tau$ as a $\AC_\FC$-module.
\end{enumerate}
\end{proposition}

\begin{proof}
(1) By Lemma~\ref{lemRestriction_bijection}, there exists $\lambda\in Y^+$ such that $\tau^+(\lambda)=0$. Then if $x\in I_{\tau^+}^+\setminus\{0\}$, $Z^\lambda.x=0$. If $M$ is a $\AC_\FC$-module, one has $Z^{-\lambda}.Z^\lambda.y=y\neq 0$ for every $y\in M\setminus\{0\}$.  The similar statement for $G^+$ is a consequence of Lemma~\ref{lemExtensions_representations}(3).

(2) The statement for $\widehat{I(\tau^+)^+}$ follows from Lemma~\ref{lemExtensions_representations}(1). The statement for $I_\tau$ follows from Proposition~\ref{propEquality_submodules_H+BLH}.  By Proposition~\ref{propIso_group_theoretic_algebraic},  the actions of $\HC_\FC$ on $I_{\tau,G^+}$ and $I_{\tau,G}$ extend to actions of $\AC_\FC$ on $I_{\tau,G^+}$ and $I_{\tau,G}$.  
\end{proof}

\begin{appendix}

\section{Existence of one dimensional representations of $^{\mathrm{BL}}\mathcal{H}_\C$}\label{appExistence_one_dimensional_representations}

In this section, we prove the existence of one dimensional representations of $\AC_\C$, when $\sigma_s=\sigma_s'=\sigma$, for all $s\in \SCC$.

 \begin{lemma}\label{lemOne_dimensional_representations}
Assume that $\FC=\C$ and that there exists $\sigma\in  \C$ such that  $\sigma_s=\sigma_s'=\sigma$ for all $s\in \SCC$ and such that $|\sigma|\neq 1$. Let $\epsilon\in \{-1,1\}$ and  $\tau\in T_\C$ be such that $\tau(\alpha_s^\vee)=\sigma^{2\epsilon}$ for all $s\in \SCC$. Then $I_\tau$ admits a unique maximal proper submodule $M$. Moreover, $I_\tau= M\oplus \C 1\otimes_\tau 1$ and if $x\in I_\tau/M$, then $Z^\lambda.x=\tau(\lambda).x$ and $H_w.x=(\epsilon \sigma^\epsilon)^{\ell(w)}.x$ for all $(w,\lambda)\in W^v\times Y$.
\end{lemma}

\begin{proof}
By Lemma~\ref{lemExistence_character_values_alphavee}, such a $\tau$ exists. Let $q=\sigma^2$. Let $\h:Y\rightarrow \Q$ be a $\Z$-linear map such that $\h(\alpha_s^\vee)=1$ for all $s\in \SCC$. Then one has $\tau(\alpha^\vee)=q^{\epsilon\h(\alpha^\vee)}$ for all $\alpha^\vee\in \Phi^\vee$.

 Let $s\in \SCC$. With the same notation as in Lemma~\ref{lemCondition on values for isomorphisms}, let $\phi_s=\phi(s.\tau,\tau):I_{s.\tau}\rightarrow I_{\tau}$. Then by Lemma~\ref{lemCondition on values for isomorphisms} $M_s:=\mathrm{Im}(\phi_s)$ is a proper submodule of $I_\tau$. Moreover,  $H_s-\epsilon\sigma^\epsilon\otimes_\tau 1\in M_s$. Let $M=\sum_{s\in \SCC} M_s$. Let $w\in W^v\setminus\{1\}$ and $w=s_1\ldots s_k$ be a reduced expression. Let $v=ws_k$. Then $H_v.(H_{s_k}-\epsilon \sigma^\epsilon)=H_w-\epsilon\sigma^\epsilon H_v\in M_{s_k}$. Therefore, for all $w\in W^v\setminus\{1\}$, there exists $x_w\in M$ such that $\pi^H_w(x_w)=1$ and $x_w\in M\cap I_\tau^{\leq w}$. By induction on $\ell(w)$ we deduce that $M+\C 1\otimes_\tau 1=I_\tau$.
 
 By \cite[Lemma 2.4 a)]{gaussent2014spherical}, $\tau\in T_\C^{\mathrm{reg}}$.  Moreover, by Proposition~\ref{propDecomposition_generalized_weight_spaces_Itau}~(\ref{itDecomposition_Itau_regular_case}), \[I_\tau=\bigoplus_{w\in W^v} I_{\tau}(w.\tau)\] and if we  choose  $\xi_v\in I_\tau(v.\tau)\setminus \{0\}$ for all $v\in W^v$, then $(\xi_v)_{v\in W^v}$ is a basis of $I_{\tau}$. For $w\in W^v$, let  $\pi^\xi_w:I_\tau\rightarrow \C$ be the linear map defined by $\pi^\xi_w(\xi_v)=\delta_{v,w}$ for all $v\in W^v$. As $\xi_1\in \C 1\otimes_\tau 1$, one has $\pi^\xi_1(M_s)=\{0\}$ for all $s\in \SCC$.
  Thus $I_\tau =M\oplus \C 1\otimes_\tau 1$. Moreover, $M\subset (\pi^{\xi}_1)^{-1}(\{0\})$ and by dimension $M=\pi^{\xi}_1(\{0\})$.  We deduce that $M$ is the unique maximal proper submodule of $I_\tau$ and the lemma follows.
\end{proof}

\begin{remark}
Actually, the representations constructed in Lemma~\ref{lemOne_dimensional_representations} generalize the well known trivial representation (when $\epsilon=1$) and Steinberg representation (when $\epsilon=-1$). For simplicity, we assumed all the $\sigma_s,\sigma_s'$ to be equal, but this is not necessary. We can also construct these representations directly by setting $\mathrm{triv}(H_s)=\sigma_s$, $\mathrm{triv}(Z^{\alpha_s^\vee})=\sigma_s\sigma'_s$, $\mathrm{St}(H_s)=-\sigma_s^{-1}$, $\mathrm{St}(Z^{\alpha_s^\vee}) = \sigma_s^{-1} \sigma_s'^{-1}$. Using the fact that the relations (BL1) to (BL4) are preserved by  $\mathrm{triv}$ and $\mathrm{St}$, we can extend them to representations of $\AC_\C$ over $\C$. 
\end{remark}

\section{Examples of possibilities for $W_\tau$ for size $2$ Kac-Moody matrices}\label{secExample_possibilities_Wtau}

In this section, we prove that there exist size $2$ Kac-Moody matrices such that for each subgroup $H$ of $W^v$, there exist $\tau\in T_\C$ such that $W_\tau$ is isomorphic to $H$.  We assume that $\alpha_s(Y)=\Z$ for all $s\in \SCC$ and thus $\Wta=W_\tau$.  We already proved the existence of regular elements in Lemma~\ref{lemDensity_regular_elements}. If $\tau\in T_\C$ is such that $\tau(\alpha_{s_1}^\vee)=1$ and $\tau(\alpha_{s_2}^\vee)$ is not a root of $1$, then $W_\tau=\{1,s_1\}$.

\begin{lemma}
Let $A=(a_{i,j})_{(i,j)\in \llbracket 1,2\rrbracket^2}$ be a Kac-Moody matrix. Assume that $a_{1,2}$ and $a_{2,1}$ are even and such that $a_{1,2}a_{2,1}$ is greater than $6$. Let $\gamma_2$ be a primitive $\frac12(a_{1,2}a_{2,1}-4)$-th root of $1$. Let $\gamma_1=\gamma_2^{\frac12a_{1,2}}$. Let $\tau:Y=\Z\alpha_1^\vee\oplus \Z\alpha_2^\vee \rightarrow \C^*$ be the group morphism defined by $\tau(\alpha_i^\vee)=\gamma_i$ for both $i\in \{1,2\}$. Then $W_\tau=\langle s_1s_2\rangle\simeq \Z$.
\end{lemma}

\begin{proof}
Let $\tau'\in T_\C$ and $\gamma_i'=\tau'(\alpha_i^\vee)$ for both $i\in \{1,2\}$. For $\lambda\in Y$, one has $(s_2-s_1).\lambda=\alpha_1(\lambda)\alpha_1^\vee-\alpha_2(\lambda)\alpha_2^\vee$. Thus \[
\begin{aligned} s_1.\tau'=s_2.\tau' & \Longleftrightarrow \forall \lambda\in Y, \tau'(\alpha_1(\lambda)\alpha_1^\vee-\alpha_2(\lambda)\alpha_2^\vee)=1\\ &\Longleftrightarrow \forall \lambda\in Y, \gamma_1'^{\alpha_1(\lambda)}= \gamma_2'^{\alpha_2(\lambda)}\\ &\Longleftrightarrow (\gamma_1')^{2}=(\gamma_2')^{a_{1,2}}\mathrm{\ and\ }(\gamma_2')^{2}=(\gamma_1')^{a_{2,1}}.\end{aligned}\] Thus $s_1.s_2.\tau=\tau$. Moreover $s_2.\tau\neq \tau$ and hence $W_\tau=\langle s_1s_2\rangle$.
\end{proof}

If $\tau=\mathds{1}:Y\rightarrow \{1\}$, then $W_\tau =1$. The following lemma proves that $W_\tau$ can be a proper subgroup of $W^v$ isomorphic to the infinite dihedral group.

\begin{lemma}\label{lemTau_fixer_2_generators}
Let $A=(a_{i,j})_{(i,j)\in  \llbracket 1,2\rrbracket^2}$ be an irreducible Kac-Moody matrix which is not a Cartan matrix. One has $a_{1,2}a_{2,1}\geq 4$ and maybe considering ${^t\!}A$, one may assume $a_{1,2}\leq -2$. Write $W^v=\langle s_1,s_2\rangle$. Let $\gamma_2$ be an $a_{1,2}$-th primitive root of $1$ and  $\tau\in T_\C$ be defined by $\tau(\alpha_{s_1}^\vee)=1$ and $\tau(\alpha_{s_2}^\vee)=\gamma_2$. Then $W_\tau=\langle s_1,s_2s_1s_2\rangle$. 
\end{lemma}

\begin{proof}
Let $\widetilde{\tau}=s_2.\tau$. Let us prove that $s_1.\widetilde{\tau}=\widetilde{\tau}$, i.e that $\widetilde{\tau}(\alpha_{s_1}^\vee)=1$.  One has $\widetilde{\tau}(\alpha_{s_1}^\vee)=\tau(s_2.\alpha_{s_1}^\vee)=\tau(\alpha_{s_1}^\vee-\alpha_{s_2}(\alpha_{s_1}^\vee)\alpha_{s_2}^\vee)=\tau(\alpha_{s_2}^\vee)^{-a_{1,2}}=1$. Thus $W_\tau\ni\{s_1,s_2s_1s_2\}$. Therefore $W^v/W_\tau=\{W_\tau,t.W_\tau\}$. Moreover $t\notin W_\tau$, thus $[W^v:W_\tau]=2$ and hence $W_\tau=\langle s_1,s_2s_1s_2\rangle$.\end{proof}

\end{appendix}

\printindex

\bibliography{/home/auguste_pro/Documents/Projets/bibliographie.bib}
\bibliographystyle{plain}
\end{document}